\documentclass[11pt,reqno]{amsart}

\usepackage[a4paper, margin=2cm]{geometry}
\usepackage{amsmath,amssymb,amsthm}
\usepackage{hyperref}

\theoremstyle{plain}
\newtheorem{theorem}{Theorem}[section]
\newtheorem{proposition}[theorem]{Proposition}
\newtheorem{lemma}[theorem]{Lemma}
\newtheorem{corollary}[theorem]{Corollary}

\theoremstyle{definition}
\newtheorem{definition}[theorem]{Definition}
\newtheorem{convention}[theorem]{Convention}

\theoremstyle{remark}
\newtheorem{remark}[theorem]{Remark}

\newcommand{\Q}{\mathbb{Q}}
\newcommand{\Z}{\mathbb{Z}}
\newcommand{\R}{\mathbb{R}}
\newcommand{\Ok}{\mathcal{O}}
\newcommand{\pp}{\mathfrak{p}}
\newcommand{\mmm}{\mathfrak{m}}
\newcommand{\PP}{\mathfrak{P}}
\newcommand{\qq}{\mathfrak{q}}
\newcommand{\aaa}{\mathfrak{a}}
\newcommand{\bb}{\mathfrak{b}}
\newcommand{\hh}{\mathfrak{h}}
\newcommand{\nn}{\mathfrak{n}}
\newcommand{\ssc}{\mathfrak{s}}
\newcommand{\vv}{\mathfrak{v}}
\newcommand{\ord}{\operatorname{ord}}
\newcommand{\rank}{\operatorname{rank}}
\newcommand{\gen}{\operatorname{gen}}
\newcommand{\spn}{\operatorname{spn}}
\newcommand{\cls}{\operatorname{cls}}
\newcommand{\covol}{\operatorname{covol}}
\newcommand{\Nm}{\mathrm{N}}
\newcommand{\NFQ}{\Nm_{F/\Q}}
\newcommand{\dual}[1]{#1^{\#}}
\newcommand{\Kv}{K(v)}
\newcommand{\bigperp}{\mathop{\perp}}
\newcommand{\eps}{\varepsilon}
\newcommand{\Dodd}{D_{\mathrm{odd}}}
\newcommand{\Dp}{D_{2}}
\newcommand{\pref}[2]{\hyperref[#1:#2]{\ref*{#1}(#2)}}

\begin{document}

\title[Explicit Tartakowsky and criterion sets]
{Explicit Tartakowsky theorem and universality criterion sets over totally real number fields}

\author{V\'{\i}t\v{e}zslav Kala}
\author{Siu Hang Man}
\author{Robin Visser}
\author{Pavlo Yatsyna}

\address{Charles University, Faculty of Mathematics and Physics, Department of Algebra,
Sokolovsk\'a 83, 18600 Praha 8, Czech Republic}
\email{vitezslav.kala@matfyz.cuni.cz}
\email{siu-hang.man@matfyz.cuni.cz}
\email{robin.visser@matfyz.cuni.cz}
\email{p.yatsyna@matfyz.cuni.cz}

\thanks{V.K., S.H.M., and P.Y. were supported by Czech Science Foundation (GA\v{C}R) grant
26-20514S, S.H.M. by Charles University programme PRIMUS/25/SCI/008, R.V. and P.Y. by
Charles University programme PRIMUS/24/SCI/010, and R.V. by Charles University programme
UNCE/24/SCI/022.}


\date{November 16, 2026}

\begin{abstract}
We prove an explicit version of Tartakowsky's theorem over all totally real number fields. Namely, a positive definite integral lattice of rank at
least five represents every totally positive integer that it represents locally, provided that the
norm of the integer exceeds an explicit bound, exponential in powers of the volume of the lattice
and of the discriminant of the field. As an application, we give explicit criterion
sets for universal quadratic forms over all but finitely many real quadratic fields. Finally, we determine the minimal criterion sets for an infinite family
of real quadratic fields.
\end{abstract}

\maketitle


\section{Introduction}\label{sec:intro}

Which integers are represented by a given quadratic form is one of the oldest questions in number
theory. Lagrange's four-square theorem asserts that the form $x^{2}+y^{2}+z^{2}+w^{2}$ represents
every positive integer, and Ramanujan \cite{Ra} determined the diagonal quaternary forms with this
property, up to one error in his list. Positive definite forms representing every positive integer are called
\emph{universal}, and the $15$-Theorem of Conway--Schneeberger \cite{Bh, Co} and the
$290$-Theorem of Bhargava--Hanke \cite{BH} show that universality is decided by the
representation of finitely many integers.

A basic tool is the comparison with local representations. By the Hasse--Minkowski theorem and
the theory of genera, a positive integer that is represented by a positive definite form $Q$ over
$\Z_p$ for every prime $p$ is represented by some form in the genus of $Q$, but not necessarily by
$Q$ itself \cite{Si1}, \cite[\S102]{OM2}. Tartakowsky \cite{T} proved that the difference
disappears for large integers when $Q$ has at least five variables, i.e., such a form represents
every sufficiently large integer that it represents locally. For quaternary forms the analogous
statement holds for integers with bounded divisibility by the anisotropic primes \cite{Kl, T}, see
also \cite{BO}. For ternary forms, Duke--Schulze-Pillot \cite{DSP} proved a statement of this kind
for integers outside finitely many square classes, using bounds for Fourier coefficients of modular
forms of half-integral weight, but their bound is not effective in general. See \cite{Cog} for an
extension to totally real fields, and \cite{Kit, SP2} for accounts of these developments.

Effective versions of Tartakowsky's theorem were obtained by Watson \cite{W}, with unspecified
implied constants, and by Hsia--Icaza \cite{HI}, who produced the first explicit bound by
arithmetic means, following \cite{Kn} and \cite{HKK}. Analytic methods also yield explicit bounds,
see, e.g., \cite{Ha}, \cite{Bl}, and \cite{BD}. Over totally real number fields,
Hsia--Kitaoka--Kneser \cite{HKK} proved an asymptotic local--global principle for representations
of lattices of rank $n$ by lattices of rank at least $2n+3$. Over $\Z$, Ellenberg--Venkatesh
\cite{EV} proved such a principle in codimension at least $5$, under an additional local
condition, by ergodic methods. Hsia--Icaza \cite[\S6]{HI} extended their explicit bound to totally
real number fields $F$, under the standing hypothesis, stated at \cite[p. 250]{HI}, that $2$ does
not ramify in $F$.

Over totally real number fields, explicit local--global principles are particularly relevant for
universal quadratic forms, i.e., positive definite forms that represent every totally
positive integer (see the surveys \cite{Kim, Ka}). Maa\ss{} \cite{Ma} showed that the sum of three
squares is universal over $\Q(\sqrt5)$, while Siegel \cite{Si2} proved that a sum of squares is
universal only over $\Q$ and $\Q(\sqrt5)$. Universal forms exist over every totally real field by
\cite{HKK}, but their minimal rank is often large \cite{BK1, BK2, Ya, KT, Man, KYZ, KMa}, and
quadratic forms with coefficients in a proper subfield are rarely universal \cite{KY1}. Kitaoka
conjectured that only finitely many totally real fields admit a universal ternary form. See
\cite{CKR, EK, KTZ, KKP, KY2, KK+} for results towards this conjecture.

Universal forms over totally real fields also admit an analogue of the $15$- and $290$-Theorems.
Chan--Oh \cite{CO} proved that finite \emph{universality criterion sets} exist over every totally
real field, extending a finiteness theorem of Kim--Kim--Oh \cite{KKO2} over $\Z$. Over $\Z$, the
minimal criterion sets for the representation of all integral quadratic forms of rank $2$,
respectively $8$, are unique \cite{KKO1, Oh, Kom1, Kom2}, while Elkies--Kane--Kominers \cite{EKK}
showed that for other families of quadratic forms the minimal criterion sets need not be unique and
may even vary in size. For universality over totally real fields, Kala--Krásenský--Romeo
\cite{KKR} showed that there is a unique minimal criterion set, up to multiplication by squares of
units, and Krásenský--Romeo \cite{KR} developed escalations over real quadratic fields. However,
criterion sets are known explicitly only in very few cases, e.g., over $\Q(\sqrt5)$ by Lee
\cite{Lee}, and \cite[Question 6.2]{KKR} asks for an explicit bound for the norms of their
elements.

In this paper we prove an explicit version of Tartakowsky's theorem over totally real number
fields with no assumption on the ramification of $2$ (\S\ref{subsec:tartintro}), and we apply it to
give explicit universality criterion sets over real quadratic fields, answering
\cite[Question 6.2]{KKR} for real quadratic fields of large discriminant
(\S\ref{subsec:critintro}). For the fields $\Q(\sqrt{n^{2}-1})$, we also determine the minimal
criterion sets, without Tartakowsky's theorem.

\subsection{An explicit Tartakowsky theorem}\label{subsec:tartintro}

Throughout, $F$ is a totally real number field of degree $l:=[F:\Q]$, with ring of integers $\Ok$,
discriminant $d_F$, and real embeddings $\tau_1,\dots,\tau_l$. All lattices are quadratic lattices,
i.e., an $\Ok$-lattice is a finitely generated $\Ok$-submodule of a quadratic space over $F$, with
the bilinear form $B$ and the quadratic form $Q(x)=B(x,x)$ of that space. For an $\Ok$-lattice $L$
we write $\ssc L$, $\nn L$, and $\vv L$ for its scale, norm, and volume ideals,
$Q(L):=\{Q(x)\mid x\in L\}$, and $L_\PP$ for the localization of $L$ at a prime $\PP$ of $F$. These
notions are recalled in \S\ref{sec:notation}. The lattice $L$ is \emph{integral} if
$\ssc L\subset\Ok$, and it is \emph{positive definite} if $Q(x)$ is totally positive for every
nonzero $x\in L$. Throughout the paper, all $\Ok$-lattices are positive definite. This convention
does not concern the lattices over the local rings $\Ok_\PP$. We write $\Nm\aaa$ for the absolute norm of an ideal $\aaa$ of $\Ok$, and for
an integer $m\ge5$ we put
\begin{equation}\label{eq:rho}
\rho_5:=10,\qquad \rho_m:=\frac{m(m-2)}{m-3}\quad(m\ge6).
\end{equation}

\begin{theorem}\label{cor:crude}
Let $M$ be an integral positive definite $\Ok$-lattice of rank $m\ge5$, put
$D:=\Nm(\vv M)$, and let $A\in\Ok$ be totally positive and represented by $M_\PP$ for every
prime $\PP$. If
\[
\NFQ(A)\ >\ \exp\Bigl(\bigl(2lm\,d_F\bigr)^{\,l^{2}m^{2}+4lm^{2}}\,D^{\rho_m}
\bigl(\log(4D)\bigr)^{2lm+m-1}\Bigr),
\]
then $A$ is represented by $M$.
\end{theorem}

We will deduce Theorem \ref{cor:crude} from a more precise result. Let $\qq$ be an odd prime of
$F$ that does not divide $\vv M$. In Theorem \ref{thm:main} and Corollary \ref{cor:exponent} we
will show that $A$ is represented by $M$ as soon as
\begin{equation}\label{eq:G}
\NFQ(A)\ >\ G_{m,F}\,(\Nm\qq)^{4(1+r)}\,D^{\,1+\frac{2}{m-3}+\frac{2}{m-4}},\qquad\text{where}\quad
G_{m,F}:=\bigl(2^{12}l^{5}m^{m+2}\bigr)^{l}\,d_F^{\,m+2l}
\end{equation}
and $r\ge0$ is an integer that depends only on $F$, $m$, $\qq$, and $D$ (see \eqref{eq:r}). The
existence of $r$ follows from a finiteness argument. We will then show that $r$ is smaller than the
number of isometry classes of integral lattices of rank $m-1$ and suitably bounded volume, which is
bounded by a polynomial in $D$ and $\Nm\qq$ (Propositions \ref{prop:classes} and \ref{prop:rH}). Taking for
$\qq$ a prime above the least rational prime not dividing $2D$, we will obtain
Theorem \ref{thm:explicitintro}, and Theorem \ref{cor:crude} is a simplified form of it. 

Corollary \ref{cor:exponent} also gives separate exponents for the dyadic and the odd part of $D$,
and Corollary \ref{cor:nonclassical} treats lattices with $\nn M\subset\Ok$, which correspond to
quadratic forms $\sum_{i\le j}a_{ij}x_ix_j$ with $a_{ij}\in\Ok$.

\subsection{Universality criterion sets}\label{subsec:critintro}

An $\Ok$-lattice $L$ is \emph{integer-valued} if $\nn L\subset\Ok$, i.e., if
$Q(L)\subset\Ok$. Every integral lattice is integer-valued, as $\nn L\subset\ssc L$. Let $\Ok^{+}$ be
the set of totally positive elements of $\Ok$. For $S\subset\Ok^{+}$, an integer-valued lattice $L$
is \emph{$S$-universal} if $S\subset Q(L)$, and it is \emph{universal} if it is
$\Ok^{+}$-universal. In \cite{KKR}, integer-valued lattices are simply called lattices, and integral
lattices are called classical.

\begin{definition}\label{def:crit}
A set $\mathcal C\subset\Ok^{+}$ is a \emph{universality criterion set for integer-valued lattices} if every
$\mathcal C$-universal integer-valued lattice is universal. It is a \emph{universality criterion set
for integral lattices} if every $\mathcal C$-universal integral lattice is universal.
\end{definition}

For $F=\Q$, the $15$-Theorem \cite{Bh} and the $290$-Theorem \cite{BH}
state that the sets
\begin{align*}
\mathcal C_{15}&:=\{1,2,3,5,6,7,10,14,15\},\\
\mathcal C_{290}&:=\{1,2,3,5,6,7,10,13,14,15,17,19,21,22,23,26,29,30,\\
&\phantom{:=\{}31,34,35,37,42,58,93,110,145,203,290\}
\end{align*}
are universality criterion sets over $\Q$, for integral lattices and for integer-valued lattices,
respectively, and that no element of these sets can be omitted.
Over real quadratic fields of large discriminant, Theorem \ref{cor:crude} yields the following
explicit criterion sets.

\begin{theorem}\label{cor:critexplicit}
Let $F$ be a real quadratic field.
\begin{enumerate}
\item\label{cor:critexplicit:1} If $d_F\ge1641$, then
$\bigl\{A\in\Ok^{+}\mid\NFQ(A)\le\exp\bigl((24\,d_F)^{476}\bigr)\bigr\}$ is a universality
criterion set for integral lattices.
\item\label{cor:critexplicit:2} If $d_F\ge25\,761$, then
$\bigl\{A\in\Ok^{+}\mid\NFQ(A)\le\exp\bigl((24\,d_F)^{486}\bigr)\bigr\}$ is a universality
criterion set for integer-valued lattices.
\end{enumerate}
\end{theorem}

Combining part \hyperref[cor:critexplicit:1]{(1)} with further computations, Kr\'asensk\'y \cite{Kr} has extended it to all but 12 real quadratic fields.

Let $\mathcal U_F$ be the unit group of $\Ok$, and $\mathcal U_F^{+}$ its subgroup of totally
positive units. A lattice represents $A$ if and only if it represents $A\eta^{2}$ for
$\eta\in\mathcal U_F$, so only the class of an element of a criterion set modulo
$\mathcal U_F^{2}$ matters. The two sets in Theorem \ref{cor:critexplicit} are infinite, but each of them is the union of finitely many
classes modulo $\mathcal U_F^{2}$, because $\Ok^{+}$ contains only finitely many elements of bounded
norm up to multiplication by totally positive units, and $\mathcal U_F^{2}$ has finite index in
$\mathcal U_F^{+}$. Choosing one element in each class gives a finite universality
criterion set of the same kind.

Theorem \ref{cor:critexplicit} is a consequence of Theorem \ref{thm:tuniv}, which concerns
$t\Ok^{+}$-universal lattices for a positive integer $t$. There the norms are bounded in terms of
the number $\mathcal T_{F,6}(W)$ of Definition \ref{def:T}, i.e., the least bound beyond which
the conclusion of Theorem \ref{cor:crude} holds for all integral lattices $L$ of rank $5$ or $6$
with $\Nm(\vv L)\le W$. Hence every improvement of Theorem \ref{cor:crude} directly
improves Theorem \ref{cor:critexplicit}. The case $t=1$ of Theorem \ref{thm:tuniv} is stated
separately as Corollary \ref{cor:tone}.

The proof of Theorem \ref{thm:tuniv} starts from five linearly independent vectors of $L$ of small
norm, which we take from \cite{KK+}. They span a lattice $L_5$ that represents all local integers
at all but finitely many primes. To represent a given element $A$, we then add one further vector
of $L$, chosen in terms of $A$, such that the resulting lattice $L_A$ represents $A$ at every prime
and has bounded volume. Finally, Tartakowsky's theorem is applied to $L_A$, or to a rescaling of
$L_A$ if $L$ is not integral.

For the fields $\Q(\sqrt{n^{2}-1})$ with $n$ large we determine minimal universality criterion
sets, independently of Tartakowsky's theorem.

\begin{theorem}\label{thm:pell}
Let $n$ be an integer such that $d:=n^{2}-1$ is squarefree, let $F=\Q(\sqrt d)$ and
$\eps:=n+\sqrt d$.
\begin{enumerate}
\item\label{thm:pell:1} If $n\ge114$, then
\[
S_n:=\mathcal C_{15}\cup\eps\,\mathcal C_{15}
=\{1,2,3,5,6,7,10,14,15,\ \eps,2\eps,3\eps,5\eps,6\eps,7\eps,10\eps,14\eps,15\eps\}
\]
is a universality criterion set for integral lattices, and no element of $S_n$ can be omitted,
i.e., for every $s\in S_n$ the set $S_n\smallsetminus\{s\}$ is not a universality criterion set for
integral lattices.
\item\label{thm:pell:2} If $n\ge168\,202$, then $S'_n:=\mathcal C_{290}\cup\eps\,\mathcal C_{290}$ is a universality
criterion set for integer-valued lattices, and no element of $S'_n$ can be omitted, i.e., for every
$s\in S'_n$ the set $S'_n\smallsetminus\{s\}$ is not a universality criterion set for
integer-valued lattices.
\end{enumerate}
\end{theorem}

The sets $S_n$ and $S'_n$ have $18$ and $58$ elements. This makes precise the result announced in \cite[\S6.1]{KKR}, namely that the
sizes of these sets are bounded. To our knowledge, apart from
$\Q$ and $\Q(\sqrt5)$, these are the first fields for which a minimal criterion set is completely
known. Over $\Q(\sqrt2)$ and $\Q(\sqrt3)$, Krásenský--Romeo \cite[Theorem 5.2]{KR} determined the elements of
norm below $250\,000$ of the minimal criterion set for integral lattices.

There are infinitely many fields to which Theorem \ref{thm:pell} applies, as $n^{2}-1$ is squarefree for a positive proportion of integers $n$. The proof that $S_n$ and $S'_n$ are
criterion sets uses the $15$- and $290$-Theorems and Kim's description
\cite[Lemma 4]{Ki} of $\Ok^{+}$ in terms of the powers of $\eps$. That no element can be omitted
follows from the results of \cite{KKR} on truants.

\subsection{The proof of Theorem \ref{cor:crude}}\label{subsec:proofintro}

Our proof follows the arithmetic method of Hsia--Icaza \cite{HI}. Let $M$, $A$, and $\qq$ be as in
\S\ref{subsec:tartintro}. We will find a vector $v\in M$ and an element $x\in\Ok$ such that
$\beta:=A-Q(v)x^{2}$ is totally positive and is represented by the lattice $K:=v^{\perp}\cap M$ of
rank $m-1$. Then $A=Q(xv)+\beta$ is represented by $M$.

To represent $\beta$ by $K$, we will prescribe $v$ modulo an ideal $\aaa$ and $x$ modulo an ideal
$\hh$ in such a way that $\beta$ is represented by the sublattice $\qq^{r}K$ of $K$ locally at
every prime. Then $\beta$ is represented by a lattice in the genus of $\qq^{r}K$, and even in its
proper spinor genus, by a theorem of Hsia \cite{Hs}. As $K$ is isotropic at $\qq$, strong
approximation shows that every lattice in this proper spinor genus is isometric to a sublattice of
$K$, provided that $r$ is large enough. By the geometry of numbers, we can choose $v$ and $x$ so
small that $\beta$ is totally positive whenever $\NFQ(A)$ satisfies \eqref{eq:G}, and that the
volume of $K$ is bounded in terms of $D$ and $\Nm\qq$. Hence only finitely many lattices $K$ occur
up to isometry, and one integer $r$ suffices for all of them. Finally, in \S\ref{sec:r} we will
bound $r$ by a class number, using Kneser's neighbors at $\qq$.

The main new points concern the dyadic primes $\pp$, where \cite{HI} uses results that are special
to $\ord_\pp(2)=1$. We will bound the value set of $M_\pp$ from below with the help of maximal
lattices (Proposition \ref{prop:lowerbound}), and we will obtain the congruence conditions on $v$
from the local square theorem (Proposition \ref{prop:covering}). Both arguments are valid at every
dyadic prime. Two further differences concern the odd primes $\PP$ with $\ord_\PP\vv M=m-4$, where
the argument of \cite{HI} has to be modified for $m=5$ (Remark \ref{rem:HIgap}), and the real
places, where we use Minkowski's theorems for a weighted Euclidean norm instead of reduction
theory. As a result, the exponent of $D$ in \eqref{eq:G} does not depend on $l$, and for $F=\Q$ it
is the exponent of \cite[Corollary 3.1]{HI}. 

\section{Preliminaries}\label{sec:notation}

Let us recall the definitions used in this paper.

\subsection{Number fields}
Throughout the paper, $F$ is a totally real number field of degree $l=[F:\Q]$, with ring of
integers $\Ok$, discriminant $d_F$, and real embeddings $\tau_1,\dots,\tau_l$. We write $\NFQ$ for
the norm from $F$ to $\Q$ and $\mathcal U_F$ for the unit group of $\Ok$.

An element $\alpha\in F$ is \emph{totally positive} if $\tau_j(\alpha)>0$ for $j=1,\dots,l$, and
$\Ok^{+}$ is the set of totally positive elements of $\Ok$. For $\alpha,\beta\in F$ we write
$\alpha\prec\beta$ if $\beta-\alpha$ is totally positive, and $\alpha\preceq\beta$ if
$\alpha\prec\beta$ or $\alpha=\beta$. For a nonzero fractional ideal $\aaa$ of $\Ok$ we write $\Nm\aaa$ for its
absolute norm, so that $\Nm(\alpha\Ok)=|\NFQ(\alpha)|$ for $\alpha\in F^{\times}$.

A \emph{prime} of $F$ is a nonzero prime ideal of $\Ok$. For a prime $\PP$ we write $F_\PP$ and
$\Ok_\PP$ for the completions of $F$ and $\Ok$ at $\PP$, $\pi_\PP$ for a uniformizer, $\ord_\PP$ for
the valuation of $F_\PP$ with $\ord_\PP(\pi_\PP)=1$ (also applied to fractional ideals), $e_\PP:=\ord_\PP(2)$, and
$f_\PP$ for the residue degree of $\PP$ over $\Q$. Thus $\Nm\PP=p^{f_\PP}$ if $\PP\mid p$, and
$\sum_{\pp\mid2}e_\pp f_\pp=l$. A prime is \emph{dyadic} if it divides $2$ and \emph{odd}
otherwise. The letter $\pp$ always denotes a dyadic prime and $\qq$ an odd prime. The \emph{places} of $F$ are its primes and its real
embeddings, and $F_v$ is the completion of $F$ at the place $v$. We write $J_F$ for the idele group
of $F$.

\subsection{Quadratic spaces}
Let $E$ be one of the fields $F$, $F_\PP$, or $\R$. A \emph{quadratic space} over $E$ is a finite
dimensional $E$-vector space $V$ with a symmetric bilinear form $B\colon V\times V\to E$. We put
$Q(x):=B(x,x)$, so that $2B(x,y)=Q(x+y)-Q(x)-Q(y)$. All quadratic spaces in this paper are
nondegenerate. The \emph{Gram matrix} of $x_1,\dots,x_n\in V$ is $\bigl(B(x_i,x_j)\bigr)_{i,j}$.
We write $V=U\perp W$ if $V=U\oplus W$ and $B(U,W)=0$, and
$x^{\perp}:=\{y\in V\mid B(x,y)=0\}$. A nonzero $x\in V$ is \emph{isotropic} if $Q(x)=0$ and
\emph{anisotropic} otherwise, and $V$ is \emph{isotropic} if it contains an isotropic vector. An
element $\alpha\in E$ is \emph{represented} by $V$ if $\alpha=Q(x)$ for some $x\in V$.

A quadratic space $V$ over $F$ is \emph{positive definite} if $Q(x)$ is totally positive for
every nonzero $x\in V$. For a place $v$ of $F$ we put $V_v:=V\otimes_FF_v$, a quadratic space over
$F_v$.

We write $O(V)$ for the orthogonal group of $V$ and $O^{+}(V)$ for its subgroup of rotations, the
elements of determinant $1$. For an anisotropic $u\in V$, the \emph{symmetry}
$S_u\in O(V)$ is $x\mapsto x-2B(x,u)Q(u)^{-1}u$. Every rotation is a product
$S_{u_1}\cdots S_{u_{2k}}$ of an even number of symmetries, and its \emph{spinor norm} is the class
of $Q(u_1)\cdots Q(u_{2k})$ in $E^{\times}/(E^{\times})^{2}$. This defines a homomorphism
$\theta\colon O^{+}(V)\to E^{\times}/(E^{\times})^{2}$ \cite[\S55]{OM2}. We identify a subgroup of
$E^{\times}/(E^{\times})^{2}$ with its preimage in $E^{\times}$.

\subsection{Lattices}
Let $R$ be $\Ok$ or $\Ok_\PP$, and $E$ its field of fractions. An \emph{$R$-lattice} is a finitely
generated $R$-submodule $L$ of a quadratic space $V$ over $E$, together with the restrictions of
$B$ and $Q$ to $L$. Thus all lattices in this paper are quadratic lattices. We say that $L$ is a
lattice \emph{on} $V$ if $EL=V$. In general $L$ is a lattice on the subspace $EL$, and
$\rank L:=\dim EL$. An $\Ok$-lattice $L$ is \emph{positive definite} if $FL$ is
positive definite. Recall that all $\Ok$-lattices in this paper are positive definite.

For an $\Ok$-lattice $L$ on $V$ and a prime $\PP$, the \emph{localization}
$L_\PP:=\Ok_\PP L\subset V_\PP$ is an $\Ok_\PP$-lattice on $V_\PP$. A lattice on $V$ is determined
by its localizations, and if $L$ is a lattice on $V$ and $(Y_\PP)_\PP$ is a family of
$\Ok_\PP$-lattices on $V_\PP$ with $Y_\PP=L_\PP$ for almost all $\PP$, then there is a lattice $X$
on $V$ with $X_\PP=Y_\PP$ for every $\PP$ \cite[81:14]{OM2}. Every $\Ok_\PP$-lattice is free,
because $\Ok_\PP$ is a principal ideal domain.

For an $R$-lattice $L$ we use the following notation.
\begin{itemize}
\item $Q(L):=\{Q(x)\mid x\in L\}$ is the \emph{value set} of $L$. In general it is not an ideal.
\item The \emph{scale} $\ssc L$ and the \emph{norm} $\nn L$ are the fractional ideals of $R$
generated by $B(L,L)$ and by $Q(L)$ respectively. $L$ is \emph{integral} if $\ssc L\subset R$.
\item The \emph{volume} $\vv L$ of a free lattice $L$ with basis $x_1,\dots,x_n$ is the
fractional ideal generated by $\det\bigl(B(x_i,x_j)\bigr)$. For an $\Ok$-lattice $L$, $\vv L$ is
the fractional ideal with $(\vv L)_\PP=\vv(L_\PP)$ for every $\PP$.
\item The \emph{dual} of $L$ is $\dual L:=\{x\in EL\mid B(x,L)\subset R\}$.
\item If $L'\subset L$ are $R$-lattices of the same rank, then $L/L'\cong\bigoplus_iR/\aaa_i$ for
nonzero integral ideals $\aaa_i$, and the \emph{index ideal} is $[L:L']:=\prod_i\aaa_i$. We put
$(L:L'):=\#(L/L')=\Nm[L:L']$, and recall $\vv L'=[L:L']^{2}\,\vv L$ \cite[82:11]{OM2}.
\item A nonzero $x\in L$ is \emph{primitive} in $L$ if $L\cap Ex=Rx$. For $R=\Ok_\PP$ this means
$x\notin\PP L$.
\end{itemize}
If $L_1,L_2\subset V$ are lattices with $B(L_1,L_2)=0$ and $EL_1\cap EL_2=0$, we write
$L_1\perp L_2$ for $L_1+L_2$. A sublattice $N$ of $L$ \emph{splits} $L$ if $L=N\perp N'$ for some
sublattice $N'$ of $L$. For $\alpha_1,\dots,\alpha_n\in R$, $\langle\alpha_1,\dots,\alpha_n\rangle$
denotes a lattice with a basis $x_1,\dots,x_n$ such that $B(x_i,x_j)=0$ for $i\ne j$ and
$Q(x_i)=\alpha_i$. A lattice with such a basis is \emph{diagonalizable}. Lattices of rank $2$ and
$3$ are called \emph{binary} and \emph{ternary}, and $L$ is \emph{isotropic} if $EL$ is. A
\emph{hyperbolic plane} is a lattice with a basis $e,f$ such that $Q(e)=Q(f)=0$ and $B(e,f)=1$, and
$L$ \emph{splits a hyperbolic plane} if some hyperbolic plane $H\subset L$ splits $L$.

Now let $R=\Ok_\PP$. For $\sigma\in\Z$, an $\Ok_\PP$-lattice $L$ is \emph{$\PP^{\sigma}$-modular} if
$\dual L=\PP^{-\sigma}L$, or equivalently, if $B(x,L)=\PP^{\sigma}$ for every $x\in L\smallsetminus\PP L$
\cite[82:14, 82:17]{OM2}. A $\PP^{0}$-modular lattice is \emph{unimodular}. A modular lattice $L$ is
\emph{proper} if $\nn L=\ssc L$ and \emph{improper} otherwise. A \emph{Jordan splitting} of $L$ is a
splitting $L=L_1\perp\dots\perp L_t$ in which each $L_i$ is $\PP^{\sigma_i}$-modular and
$\sigma_1<\dots<\sigma_t$. Every $\Ok_\PP$-lattice has one, and $t$, the $\sigma_i$, and the ranks of
the $L_i$ do not depend on the choice of the Jordan splitting \cite[91:9]{OM2}. The
\emph{unimodular component} of the Jordan splitting is $L_1$ if $\sigma_1=0$, and $0$ otherwise.

\begin{convention}\label{conv:scaling}
We will use two different operations of scaling. For
$c\in E^{\times}$ and an $R$-lattice $L$ we write $cL$ for the lattice $\{cx\mid x\in L\}$, so that
$Q(cL)=c^{2}Q(L)$ and $\ssc(cL)=c^{2}\ssc L$. We write $L^{c}$ for the $R$-module $L$ equipped with
the rescaled form $cB$, which is a lattice on the quadratic space $(EL,cB)$, so that
$Q(L^{c})=cQ(L)$, $\ssc(L^{c})=c\,\ssc L$, and $\vv(L^{c})=c^{\rank L}\vv L$. For an integer $j$ and a prime $\PP$ we write $\PP^{j}L$ for the
product of the fractional ideal $\PP^{j}$ with $L$. For an $\Ok_\PP$-lattice $L$ this is $\pi_\PP^{j}L$, so
that $Q(\PP^{j}L)=\pi_\PP^{2j}Q(L)$. For an $\Ok$-lattice $K$ we have
$(\PP^{j}K)_\PP=\pi_\PP^{j}K_\PP$ and $(\PP^{j}K)_{\PP'}=K_{\PP'}$ for every prime $\PP'\ne\PP$.
\end{convention}

\subsection{Isometries, genera, and spinor genera}
Lattices $L$ on $V$ and $L'$ on $V'$ are \emph{isometric}, written $L\cong L'$, if there is an
isometry $\sigma\colon V\to V'$ of quadratic spaces with $\sigma L=L'$. The \emph{isometry class}
of $L$ consists of all lattices isometric to $L$. A lattice $K$ is \emph{represented} by $L$ if $L$
contains a sublattice isometric to $K$, and $\alpha\in E$ is represented by $L$ if
$\alpha\in Q(L)$. A lattice or an element is represented by a set of lattices if it is represented
by some member of that set.

Let $L$ be an $\Ok$-lattice on a quadratic space $V$ over $F$. The \emph{adelic orthogonal group}
$O_{\mathbb A}(V)$ consists of the families $\Sigma=(\sigma_v)_v$, where $v$ runs through the places
of $F$, such that $\sigma_v\in O(V_v)$ for every $v$ and $\sigma_\PP L_\PP=L_\PP$ for almost all
primes $\PP$. This group does not depend on $L$. We write $O^{+}_{\mathbb A}(V)$ for the subgroup of
families of rotations, and $\theta(\Sigma)\in J_F/J_F^{2}$ for the class of the
idele $(\theta(\sigma_v))_v$ \cite[\S\S101--102]{OM2}. The kernel of $\theta$ on
$O^{+}_{\mathbb A}(V)$ is denoted by $O'_{\mathbb A}(V)$. For
$\Sigma\in O_{\mathbb A}(V)$, $\Sigma L$ is the lattice on $V$ with $(\Sigma L)_\PP=\sigma_\PP L_\PP$
for every $\PP$. The \emph{genus}, \emph{class}, \emph{proper class}, \emph{spinor genus}, and
\emph{proper spinor genus} of $L$ are
\[
\gen(L)=\{\Sigma L\mid\Sigma\in O_{\mathbb A}(V)\},\quad
\cls(L)=\{\sigma L\mid\sigma\in O(V)\},\quad
\cls^{+}(L)=\{\sigma L\mid\sigma\in O^{+}(V)\},
\]
\[
\spn(L)=\{\sigma\Sigma L\mid\sigma\in O(V),\ \Sigma\in O'_{\mathbb A}(V)\},\qquad
\spn^{+}(L)=\{\sigma\Sigma L\mid\sigma\in O^{+}(V),\ \Sigma\in O'_{\mathbb A}(V)\}
\]
\cite[\S102]{OM2}. A lattice $L'$ on $V$ lies in $\gen(L)$ if and only if $L'_\PP\cong L_\PP$ for
every prime $\PP$. Indeed, if $L'_\PP\cong L_\PP$ for every $\PP$, we can choose
$\sigma_\PP\in O(V_\PP)$ with $\sigma_\PP L_\PP=L'_\PP$, and we can take $\sigma_\PP=1$ for the primes
with $L'_\PP=L_\PP$, i.e., for almost all $\PP$. Then $L'=\Sigma L$ for the family $\Sigma$ with
these components at the primes and with the identity at the real places. The genus
$\gen(L)$ is the union of finitely many classes \cite[103:4]{OM2}, and we write $h(\gen L)$ for
their number.

For a finite set $S$ of places of $F$ containing the real places, we put
$\Ok_S:=\{x\in F\mid\ord_\PP x\ge0\text{ for every prime }\PP\notin S\}$, and for a lattice $L$ on $V$
we write $\Ok_SL$ for the $\Ok_S$-module generated by $L$. The class, proper class, and proper
spinor genus of $\Ok_SL$ are defined by the formulas above, in which the localizations are now taken
only at the primes $\PP\notin S$ \cite[\S\S101--102]{OM2}. We call $S$ \emph{indefinite} for $V$ if $V_v$
is isotropic for some $v\in S$.

\subsection{Local preliminaries}\label{sec:prelim}

Throughout this subsection $\PP$ is a prime of $F$, $e:=e_\PP$, and $\ord:=\ord_\PP$.

\begin{lemma}[Local square theorem {\cite[63:1]{OM2}}]\label{lem:lst}
$1+\PP^{2e+1}\subset(\Ok_\PP^{\times})^{2}$.
\end{lemma}

\begin{lemma}\label{lem:normrange}
For every lattice $L$ we have $2\,\ssc L\subset\nn L\subset\ssc L$. In particular a
$\PP^{\sigma}$-modular lattice $L$ has $\nn L=\PP^{a}$ with $\sigma\le a\le\sigma+e$, and if $\PP$ is odd
then $L$ is proper.
\end{lemma}

\begin{proof}
The inclusions are proved in \cite[\S82E]{OM2}. A $\PP^{\sigma}$-modular lattice $L$ has scale
$\PP^{\sigma}$ \cite[\S82G]{OM2}, so the inclusions give $\PP^{\sigma+e}\subset\nn L\subset\PP^{\sigma}$.
If $\PP$ is odd, then $2$ is a unit, and so $\nn L=\ssc L$.
\end{proof}

\begin{lemma}[{\cite[\S\S91C, 92]{OM2}}]\label{lem:split}
Every integral $\Ok_\PP$-lattice $M$ splits as $M=N_1\perp\dots\perp N_k$, where each $N_i$ is a
modular lattice of rank at most $2$. If $\PP$ is odd, all the $N_i$ can be chosen of rank $1$.
\end{lemma}

\begin{definition}\label{def:exponents}
Let $M$ be an integral $\Ok_\PP$-lattice of rank $n$, and fix a splitting
$M=N_1\perp\dots\perp N_k$ as in Lemma \ref{lem:split}, where $N_i$ is $\PP^{\sigma_i}$-modular and
$\sigma_1\le\sigma_2\le\dots\le\sigma_k$. The
\emph{scale exponents} of this splitting are the integers $s_1\le s_2\le\dots\le s_n$ obtained by
listing each $\sigma_i$ with multiplicity $\rank N_i$.
\end{definition}

The scale exponents do not depend on the choice of the splitting. Indeed, collecting the lattices
$N_i$ with the same $\sigma_i$ gives a Jordan splitting of $M$, and the scales and the ranks of the
components of a Jordan splitting are invariants of $M$ \cite[91:9]{OM2}. We will therefore speak of
the scale exponents of $M$. They satisfy $s_1\ge0$, $\sum_i s_i=\ord\vv M$, and $\PP^{s_1}=\ssc M$.
If $n\ge4$, then
$(n-3)s_4\le s_4+s_5+\dots+s_n\le\ord\vv M$, and so
\begin{equation}\label{eq:s4}
s_4\le\bigl\lfloor\ord\vv M/(n-3)\bigr\rfloor .
\end{equation}

\begin{lemma}[cf. {\cite[Lemma 2.4]{HI}}]\label{lem:vol}
Let $M$ be an integral $\Ok$-lattice, let $v\in M$ satisfy $Q(v)\ne0$, and put
$\Kv=v^{\perp}\cap M$. Let $\PP$ be a prime, and let $\bb$ be the ideal of $\Ok_\PP$ generated by
$B(v,M_\PP)$. Then
\[
\ord_\PP\vv\Kv=\ord_\PP\vv M+\ord_\PP Q(v)-2\ord_\PP\bb .
\]
In particular $\ord_{\PP}\vv \Kv\le\ord_{\PP}\vv M+\ord_{\PP}Q(v)$.
\end{lemma}

\begin{proof}
Put $N:=\Ok_\PP v\perp\Kv_\PP\subset M_\PP$. The $\Ok_\PP$-linear map $x\mapsto B(x,v)/Q(v)$ maps
$M_\PP$ onto $Q(v)^{-1}\bb$, and $N$ is the preimage of $\Ok_\PP$. Hence $[M_\PP:N]=Q(v)\bb^{-1}$,
and \cite[82:11]{OM2} gives $Q(v)\,\vv\Kv_\PP=\vv N=Q(v)^{2}\bb^{-2}\,\vv M_\PP$. Taking orders, we
obtain the equality. The inequality follows, because $\ord_\PP\bb\ge0$ as $M$ is integral.
\end{proof}

\section{A lower bound for local value sets}\label{sec:lower}

In this section we will show that the value set of an integral $\Ok_\PP$-lattice $M$ of rank
$\ge4$ contains an ideal that is explicitly given by the volume of $M$
(Proposition \ref{prop:lowerbound}). The proof will use maximal lattices, and it requires no
hypothesis on $e$. We keep the notation of \S\ref{sec:prelim}. A lattice $L$ on a quadratic space
$W$ over $F_\PP$ is
\emph{$\aaa$-maximal} if $Q(L)\subset\aaa$ and no lattice on $W$ properly containing $L$ has this
property \cite[\S82H]{OM2}.

\begin{lemma}\label{lem:maximal}
Let $W$ be a quadratic space over $F_\PP$ with $\dim W\ge4$, let $\aaa$ be a fractional ideal, and
let $L$ be an $\aaa$-maximal lattice on $W$. Then $Q(L)=\aaa$.
\end{lemma}

\begin{proof}
By \cite[91:3]{OM2}, $Q(L)=Q(W)\cap\aaa$, and $Q(W)=F_\PP$ because $\dim W\ge4$
\cite[63:18]{OM2}.
\end{proof}

\begin{lemma}[{\cite[82:18]{OM2}}]\label{lem:exist}
Let $N$ be a lattice and $\aaa:=\nn N$. Then there is an $\aaa$-maximal lattice $L\supset N$ on
$F_\PP N$.
\end{lemma}

\begin{lemma}\label{lem:conductor}
Let $N=N_1\perp\dots\perp N_k$ with $N_i$ a $\PP^{\sigma_i}$-modular lattice and
$\rank N\ge4$. Write $\nn N=\PP^{a}$, $\sigma_{\min}=\min_i\sigma_i$,
$\sigma_{\max}=\max_i\sigma_i$. Then $\sigma_{\min}\le a\le\sigma_{\min}+e$ and
\[
Q(N)\ \supset\ \PP^{\,2\sigma_{\max}+2e-a}.
\]
\end{lemma}

\begin{proof}
We have $\ssc N=\PP^{\sigma_{\min}}$, so the bounds for $a$ follow from Lemma \ref{lem:normrange}.

By Lemma \ref{lem:exist}, there is a $\PP^{a}$-maximal lattice $L\supset N$ on $F_\PP N$. Since
$Q(L)\subset\PP^{a}$, Lemma \ref{lem:normrange} gives $2B(L,L)\subset\PP^{a}$, i.e.,
$B(L,L)\subset\PP^{a-e}$. In particular $B(x,N)\subset\PP^{a-e}$ for every $x\in L$, which means
that $L\subset\PP^{a-e}\dual N$. As $\dual N=\bigperp_i\dual{N_i}$ and
$\dual{N_i}=\PP^{-\sigma_i}N_i$, we get
\[
L\ \subset\ \PP^{a-e}\dual N=\bigperp_i \PP^{\,a-e-\sigma_i}N_i\ \subset\
\PP^{\,a-e-\sigma_{\max}}N .
\]
It follows that $\PP^{\,\sigma_{\max}+e-a}L\subset N$. The exponent $\sigma_{\max}+e-a$ is
nonnegative, because $a\le\sigma_{\min}+e\le\sigma_{\max}+e$. Finally $Q(L)=\PP^{a}$ by
Lemma \ref{lem:maximal}, and so Convention \ref{conv:scaling} gives
\[
Q(N)\supset Q\bigl(\PP^{\sigma_{\max}+e-a}L\bigr)
=\PP^{\,2(\sigma_{\max}+e-a)}\,\PP^{a}=\PP^{\,2\sigma_{\max}+2e-a}. \qedhere
\]
\end{proof}

\begin{proposition}\label{prop:lowerbound}
Let $M$ be an integral $\Ok_\PP$-lattice of rank $n\ge4$ with scale exponents
$0\le s_1\le\dots\le s_n$. Then
\[
Q(M)\ \supset\ \PP^{\,s_4+2e+1}\ \supset\ \PP^{\,\lfloor \ord\vv M/(n-3)\rfloor+2e+1}.
\]
\end{proposition}

\begin{proof}
Let $M=N_1\perp\dots\perp N_t$ be the splitting of Definition \ref{def:exponents}, in which $N_i$
is $\PP^{\sigma_i}$-modular of rank at most $2$ and $\sigma_1\le\dots\le\sigma_t$. Let $k$ be minimal
such that the lattice $N:=N_1\perp\dots\perp N_k$ has rank at least $4$. Then
$N_1\perp\dots\perp N_{k-1}$ has rank at most $3$, so the fourth scale exponent $s_4$ comes from
the block $N_k$. Thus $\sigma_k=s_4$ and $\sigma_i\le s_4$ for every $i\le k$.

We will now replace $N$ by a sublattice whose blocks have almost equal scales. Put
$c_i:=\lfloor(s_4-\sigma_i)/2\rfloor\ge0$ for $i\le k$ and
$N':=\bigperp_{i\le k}\pi_\PP^{c_i}N_i\subset N$. The lattice $\pi_\PP^{c_i}N_i$ is
$\PP^{\sigma_i'}$-modular with $\sigma_i':=\sigma_i+2c_i\in\{s_4-1,s_4\}$, and $\sigma'_k=s_4$. Hence
the largest of the exponents $\sigma'_i$ is $s_4$ and the smallest is at least $s_4-1$. By
Lemma \ref{lem:normrange}, the integer $a':=\ord\nn N'$ satisfies $a'\ge s_4-1$. So
Lemma \ref{lem:conductor}, applied to $N'$, gives
\[
Q(M)\supset Q(N')\supset\PP^{\,2s_4+2e-a'}\supset\PP^{\,s_4+2e+1}.
\]
The last inclusion of the statement follows from \eqref{eq:s4}.
\end{proof}

At odd primes we will need a sharper bound. Recall from Lemma \ref{lem:split} that for odd $\PP$
every integral $\Ok_\PP$-lattice is diagonalizable.

\begin{lemma}\label{lem:odduni}
Let $\PP$ be odd and $L$ a unimodular $\Ok_\PP$-lattice.
\begin{enumerate}
\item\label{lem:odduni:1} If $\rank L=2$ then $Q(L)\supset\{x\in\Ok_\PP\mid\ord x\ \text{even}\}\cup\{0\}$.
\item\label{lem:odduni:2} If $\rank L\ge3$ then $Q(L)=\Ok_\PP$ and $L$ splits a hyperbolic plane.
\end{enumerate}
\end{lemma}

\begin{proof}
By \cite[92:1b]{OM2}, we have $Q(L)\supset\Ok_\PP^{\times}$ if $\rank L=2$, and $Q(L)=\Ok_\PP$ if
$\rank L\ge3$. In the first case, $Q(L)$ also contains all nonzero elements of even order, because
$Q(\pi_\PP^{k}x)=\pi_\PP^{2k}Q(x)$ for $x\in L$ and $k\ge0$.

It remains to show that $L$ splits a hyperbolic plane if $\rank L\ge3$. By \cite[92:1]{OM2}, $L$ has
an orthogonal basis $x_1,\dots,x_n$ with $Q(x_i)\in\Ok_\PP^{\times}$ for every $i$. The unimodular
lattice $\Ok_\PP x_1+\Ok_\PP x_2$ represents the unit $-Q(x_3)$ by \cite[92:1b]{OM2}, so $L$ contains
a nonzero isotropic vector, and we can choose it primitive in $L$. A primitive isotropic vector of
$L$ lies in a binary unimodular sublattice that splits $L$ \cite[82:16]{OM2}. As $\PP$ is odd, this
sublattice is a hyperbolic plane.
\end{proof}

\begin{proposition}\label{prop:oddlower}
Let $\PP$ be odd and $M$ an integral $\Ok_\PP$-lattice of rank $n\ge4$ with scale exponents
$0\le s_1\le\dots\le s_n$. Then
\[
Q(M)\ \supset\ \PP^{\,s_4}\qquad\text{and}\qquad
Q(M)\ \supset\ \PP^{\,\max\{0,\,\lfloor(\ord\vv M-1)/(n-3)\rfloor\}}.
\]
\end{proposition}

\begin{proof}
Let us prove the first inclusion. By Lemma \ref{lem:split}, the lattice $M$ has an orthogonal
basis $x_1,\dots,x_n$ with $\ord Q(x_i)=s_i$. Put $\alpha_i:=Q(x_i)$ and
$c_i:=\lfloor(s_4-s_i)/2\rfloor\ge0$ for $i\le4$, and let $N\subset M$ be the lattice with the basis
$\pi_\PP^{c_1}x_1,\dots,\pi_\PP^{c_4}x_4$. Thus
$N=\langle\pi_\PP^{2c_1}\alpha_1,\dots,\pi_\PP^{2c_4}\alpha_4\rangle$, the four diagonal entries of
$N$ have orders in $\{s_4-1,s_4\}$, and the fourth of them has order $s_4$.

If $s_4=0$, then $N$ is unimodular of rank $4$, and $Q(M)\supset Q(N)=\Ok_\PP$ by
Lemma \pref{lem:odduni}{2}. Assume now that $s_4\ge1$, and rescale the form of $N$ by
$\pi_\PP^{-(s_4-1)}$. The diagonal entries of the rescaled lattice have orders $0$ and $1$, so this
lattice equals $U\perp V^{\pi_\PP}$, where $U$ and $V$ are unimodular, $\rank V\ge1$, and
$\rank U+\rank V=4$. By Lemma \ref{lem:odduni}, $Q(U)$ contains every element of $\Ok_\PP$ of even
order if $\rank U\ge2$, and $Q(U)=\Ok_\PP$ if $\rank U\ge3$. In the same way,
$Q(V^{\pi_\PP})=\pi_\PP Q(V)$ contains every element of $\Ok_\PP$ of odd order if $\rank V\ge2$, and
$Q(V^{\pi_\PP})=\PP$ if $\rank V\ge3$. One of the three cases $\rank U\ge3$, $\rank V\ge3$, and
$\rank U=\rank V=2$ occurs, and in each of them $Q(U)\cup Q(V^{\pi_\PP})\supset\PP$. Rescaling back,
we get $Q(M)\supset Q(N)\supset\pi_\PP^{s_4-1}\PP=\PP^{\,s_4}$.

Let us prove the second inclusion. If $s_3=0$, then
$\Ok_\PP x_1+\Ok_\PP x_2+\Ok_\PP x_3\subset M$ is unimodular of rank $3$, and $Q(M)\supset\Ok_\PP$ by
Lemma \pref{lem:odduni}{2}. Otherwise $s_1+s_2+s_3\ge1$, so that
$(n-3)s_4\le s_4+\dots+s_n\le\ord\vv M-1$. Hence
$s_4\le\lfloor(\ord\vv M-1)/(n-3)\rfloor$, and the claim follows from the first inclusion.
\end{proof}

Over $\Z_p$, the weaker inclusion $Q(M)\supset\PP^{\lfloor\ord\vv M/(n-3)\rfloor}$, which follows
from Proposition \ref{prop:oddlower}, is the non-dyadic case of \cite[Lemma 2.2]{HI}.

\section{The dyadic primes}\label{sec:covering}

Throughout this section $\pp$ is a dyadic prime, $e:=e_{\pp}$, and $\ord:=\ord_{\pp}$. We will show
that every element of $Q(M_\pp)$ has the form $Q(v)c^{2}$ with $c\in\Ok_\pp$, simultaneously for
all vectors $v$ in a suitable class modulo a power of $\pp$. This will give the congruence
conditions at the dyadic primes on the vector $v$ of \S\ref{subsec:proofintro}.

\begin{proposition}\label{prop:covering}
Let $M_{\pp}$ be an integral $\Ok_\pp$-lattice of rank $n\ge4$, put $w:=\ord_{\pp}\vv M_\pp$,
\[
\kappa:=\Bigl\lfloor\frac{w}{n-3}\Bigr\rfloor+2e+1,\qquad n_0:=\kappa+e+2 ,
\]
and
\[
\mathcal Y(\pp):=\bigl\{\,y+\pp^{n_0}M_{\pp}\mid y\in M_{\pp},\ \ord(Q(y))\le\kappa+1\,\bigr\}.
\]
Then:
\begin{enumerate}
\item\label{prop:covering:1} $\#\mathcal Y(\pp)\le(\Nm\pp)^{n_0\,n}$.
\item\label{prop:covering:2} The integer $\ord(Q(y))$ depends only on the class of $y$ in $\mathcal Y(\pp)$.
\item\label{prop:covering:3} For every $\alpha\in Q(M_{\pp})$ there is a class $[y]\in \mathcal Y(\pp)$ such that every $v\in[y]$
satisfies $\alpha=Q(v)c_v^{2}$ for some $c_v\in\Ok_\pp$.
\end{enumerate}
\end{proposition}

\begin{proof}
Assertion (1) holds because $M_\pp/\pp^{n_0}M_\pp$ has $(\Nm\pp)^{n_0n}$ elements. In the rest of
the proof we will use several times that $\pp^{\kappa}\subset Q(M_{\pp})$, which holds by
Proposition \ref{prop:lowerbound}.

Let $y,v\in M_\pp$ with $v\equiv y\pmod{\pp^{n_0}M_{\pp}}$. Then
\begin{equation}\label{eq:Qdiff}
Q(v)-Q(y)=Q(v-y)+2B(y,v-y)\in\pp^{2n_0}+\pp^{n_0+e}=\pp^{n_0+e},
\end{equation}
where we used that $\nn M_\pp\subset\ssc M_{\pp}\subset\Ok_\pp$ and that $n_0\ge e$. If
$\ord(Q(y))\le\kappa+1$, then $\ord(Q(y))<\kappa+2e+2=n_0+e$, and \eqref{eq:Qdiff} shows that
$\ord(Q(v))=\ord(Q(y))$. This proves (2). The set $\mathcal Y(\pp)$ is not empty, because
$\pp^{\kappa}\subset Q(M_{\pp})$ gives a vector $y\in M_\pp$ with $\ord(Q(y))=\kappa$.

Let us prove (3). Let $\alpha\in Q(M_{\pp})$. If $\alpha=0$, we can take an arbitrary class
$[y]\in\mathcal Y(\pp)$ and $c_v=0$ for every $v\in[y]$. So assume that $\alpha\ne0$.

\emph{Case 1: $\ord(\alpha)\le\kappa+1$.} Write $\alpha=Q(y)$ with $y\in M_{\pp}$. Then
$[y]\in \mathcal Y(\pp)$. For every $v\in[y]$, \eqref{eq:Qdiff} gives
\[
\frac{Q(v)}{\alpha}\in 1+\pp^{\,n_0+e-\ord(\alpha)}\subset 1+\pp^{\,2e+1}\subset(\Ok_\pp^{\times})^{2},
\]
by Lemma \ref{lem:lst}, since $n_0+e-\ord(\alpha)\ge (\kappa+e+2)+e-(\kappa+1)=2e+1$. Hence
$\alpha=Q(v)c_v^{2}$ for a unit $c_v$.

\emph{Case 2: $\ord(\alpha)\ge\kappa+2$.} Put $j:=\lfloor(\ord(\alpha)-\kappa)/2\rfloor\ge1$ and
$\alpha_0:=\pi_\pp^{-2j}\alpha$, so that $\ord(\alpha_0)\in\{\kappa,\kappa+1\}$. Then
$\alpha_0\in\pp^{\kappa}\subset Q(M_{\pp})$, and Case 1 applies to $\alpha_0$. It gives a class
$[y]\in \mathcal Y(\pp)$ such that every $v\in[y]$ satisfies $\alpha_0=Q(v)c_v'^{2}$ for some unit
$c_v'$. Therefore $\alpha=Q(v)c_v^{2}$ with $c_v:=\pi_\pp^{j}c_v'\in\Ok_\pp$.
\end{proof}

The next lemma will give the congruence conditions at $\pp$ on the element $x$ of
\S\ref{subsec:proofintro}.

\begin{lemma}\label{lem:dyadscalar}
Let $M$ be an integral $\Ok$-lattice of rank $m\ge5$, put $w:=\ord_{\pp}\vv M$ and
$b:=\lfloor w/(m-4)\rfloor+e+1$, let $v\in M$ satisfy $Q(v)\ne0$, and put $\Kv:=v^{\perp}\cap M$.
If $c,x\in\Ok_\pp$ satisfy $x\equiv c\pmod{\pp^{b}}$, then $Q(v)(c^{2}-x^{2})\in Q(\Kv_\pp)$.
\end{lemma}

\begin{proof}
Write $x=c+t$ with $\ord(t)\ge b$. Then $c^{2}-x^{2}=-2ct-t^{2}$, where $\ord(2ct)\ge e+b$ and
$\ord(t^{2})\ge2b\ge b+e+1$, so $\ord(c^{2}-x^{2})\ge b+e$. With $s:=\ord(Q(v))$ this gives
\[
\ord\bigl(Q(v)(c^{2}-x^{2})\bigr)\ \ge\ s+\Bigl\lfloor\frac{w}{m-4}\Bigr\rfloor+2e+1\ \ge\
\Bigl\lfloor\frac{w+s}{m-4}\Bigr\rfloor+2e+1 .
\]
By Lemma \ref{lem:vol}, $\ord_\pp\vv\Kv\le w+s$. As $\Kv_\pp$ is integral of rank $m-1\ge4$,
Proposition \ref{prop:lowerbound} gives
$Q(\Kv_\pp)\supset\pp^{\lfloor\ord_\pp\vv\Kv/(m-4)\rfloor+2e+1}\supset
\pp^{\lfloor(w+s)/(m-4)\rfloor+2e+1}$.
\end{proof}

\section{The odd primes}\label{sec:otherprimes}

Throughout this section $M$ is an integral $\Ok$-lattice of rank $m\ge5$, and for $v\in M$ we
write $\Kv:=v^{\perp}\cap M$. We will keep this notation in \S\S\ref{sec:global}--\ref{sec:r}.

In this section we will study the value sets $Q(\Kv_\PP)$ at the odd primes $\PP$. If
$\ord_\PP\vv M$ is small, then $Q(\Kv_\PP)=\Ok_\PP$ for every nonzero $v$
(Lemmas \ref{lem:complement} and \ref{lem:complement6}). At the other odd primes, and at one
auxiliary odd prime $\qq$, we will impose congruence conditions on $v$
(Lemmas \ref{lem:oddT} and \ref{lem:q}).

\begin{lemma}\label{lem:complement}
Let $\PP$ be an odd prime with $\ord_{\PP}\vv M\le m-5$ and let $v\in M$ be nonzero. Then
$\Kv_\PP$ contains a unimodular lattice of rank $\ge3$, and $Q(\Kv_{\PP})=\Ok_{\PP}$.
\end{lemma}

\begin{proof}
By Lemma \ref{lem:split}, the lattice $M_\PP$ has an orthogonal basis $x_1,\dots,x_m$. Let
$M(0)$ be the lattice spanned by those $x_i$ for which $Q(x_i)$ is a unit, and $X$ the lattice
spanned by the other $x_i$. Then $M_\PP=M(0)\perp X$ and $M(0)$ is unimodular. As
$\ord_\PP\vv M=\sum_i\ord_\PP Q(x_i)$, at most $\ord_\PP\vv M$ of the elements $Q(x_i)$ are not units,
and so $\rank M(0)\ge m-\ord_{\PP}\vv M\ge5$. Write $v=v_0+v_1$ with $v_0\in M(0)$ and $v_1\in X$.

If $v_0=0$, then $M(0)\subset\Kv_\PP$, and the first assertion holds. Otherwise write
$v_0=\pi_{\PP}^{i}v_0'$ with $v_0'$ primitive in $M(0)$. Since $B(v,z)=B(v_0,z)$ for $z\in M(0)$,
we have $\Kv_{\PP}\supset (v_0')^{\perp}\cap M(0)$. We distinguish two cases.

If $Q(v_0')\in\Ok_{\PP}^{\times}$, then $\Ok_\PP v_0'$ is unimodular and therefore splits $M(0)$
\cite[82:15]{OM2}. So $M(0)=\Ok_\PP v_0'\perp M'$, where $M'$ is unimodular of rank at least $4$ and
$M'\subset(v_0')^{\perp}\cap M(0)$.

If $Q(v_0')\in\PP$, choose $y\in M(0)$ with $B(v_0',y)=1$. This is possible because $M(0)$ is
unimodular and $v_0'$ is primitive. The Gram matrix of $v_0',y$ has determinant
$Q(v_0')Q(y)-1\in\Ok_{\PP}^{\times}$, so the lattice $J:=\Ok_\PP v_0'+\Ok_\PP y$ is unimodular and
splits $M(0)$ \cite[82:15]{OM2}. Thus $M(0)=J\perp R$, where $R$ is unimodular of rank at least $3$
and $R\subset(v_0')^{\perp}\cap M(0)$.

This proves the first assertion. A unimodular lattice of rank at least $3$ has value set $\Ok_\PP$
by Lemma \pref{lem:odduni}{2}, and $Q(\Kv_\PP)\subset\Ok_\PP$ because $M$ is integral.
\end{proof}

\begin{lemma}\label{lem:complement6}
Let $m\ge6$, let $\PP$ be an odd prime with $\ord_{\PP}\vv M\le m-4$, and let $v\in M$ be nonzero.
Then $Q(\Kv_{\PP})=\Ok_{\PP}$.
\end{lemma}

\begin{proof}
Put $\pi:=\pi_\PP$, and let $x_1,\dots,x_m$ and $M_\PP=M(0)\perp X$ be as in the proof of
Lemma \ref{lem:complement}. If $\rank M(0)\ge5$, that proof applies without change. So assume that
$\rank M(0)\le4$. As $\rank M(0)\ge m-\ord_\PP\vv M\ge4$, we then have $\rank M(0)=4$. The lattice
$X$ is spanned by $m-4$ of the vectors $x_i$, the orders of the corresponding elements $Q(x_i)$ are
positive, and their sum is $\ord_\PP\vv M\le m-4$. So each of these orders equals $1$. Hence
$M_\PP=U\perp V^{\pi}$, where $U:=M(0)$ and $V$ are unimodular, $\rank U=4$, and $\rank V=m-4\ge2$.
Write $B_U$ and $B_V$ for the bilinear forms of $U$ and $V$, so that
\[
B(t+b,t'+b')=B_U(t,t')+\pi B_V(b,b')\qquad\text{for }t,t'\in U\text{ and }b,b'\in V.
\]
The lattice $\Kv_\PP=v^{\perp}\cap M_\PP$ depends only on the line $F_\PP v$, so we may assume that
$v$ is primitive in $M_\PP$. Write $v=u+z$ with $u\in U$ and $z\in V$. We distinguish three cases.

\emph{Case 1: $u\in\pi U$.} Then $z$ is primitive in $V$, and since $V$ is unimodular, there is
$y\in V$ with $B_V(z,y)=1$. For $t\in U$ put $\varphi(t):=t-\pi^{-1}B_U(u,t)\,y$. This vector lies in
$M_\PP$, because $B_U(u,t)\in\PP$, and it is orthogonal to $v$, because
\[
B(v,\varphi(t))=B_U(u,t)-\pi^{-1}B_U(u,t)\,\pi B_V(z,y)=0 .
\]
Moreover, for $t,t'\in U$ we have
\[
B(\varphi(t),\varphi(t'))=B_U(t,t')+\pi^{-1}B_U(u,t)B_U(u,t')B_V(y,y)\equiv B_U(t,t')\pmod{\PP}.
\]
Hence $\varphi(U)$ is a unimodular sublattice of $\Kv_\PP$ of rank $4$, and
$Q(\Kv_\PP)\supset\Ok_\PP$ by Lemma \pref{lem:odduni}{2}.

\emph{Case 2: $u\notin\pi U$ and $Q(u)\in\Ok_\PP^{\times}$.} Then $\Ok_\PP u$ splits $U$
\cite[82:15]{OM2}, so $u^{\perp}\cap U$ is a unimodular lattice of rank $3$. It is contained in
$\Kv_\PP$, and again $Q(\Kv_\PP)\supset\Ok_\PP$ by Lemma \pref{lem:odduni}{2}.

\emph{Case 3: $u\notin\pi U$ and $Q(u)\in\PP$.} Choose $y\in U$ with $B_U(u,y)=1$. The Gram matrix
of $u,y$ has the unit determinant $Q(u)Q(y)-1$, so $\Ok_\PP u+\Ok_\PP y$ is unimodular and splits
$U$. Thus $U=(\Ok_\PP u+\Ok_\PP y)\perp R_0$, where $R_0$ is unimodular of rank $2$ and
$R_0\subset\Kv_\PP$. For $b\in V$ put $\psi(b):=b-\pi B_V(z,b)\,y\in M_\PP$. Then
\[
B(v,\psi(b))=\pi B_V(z,b)-\pi B_V(z,b)B_U(u,y)=0 ,
\]
the vector $\psi(b)$ is orthogonal to $R_0$, and for $b,b'\in V$ we have
\[
B(\psi(b),\psi(b'))=\pi\bigl(B_V(b,b')+\pi Q(y)B_V(z,b)B_V(z,b')\bigr).
\]
Hence $\psi(V)=W^{\pi}$ for a unimodular lattice $W$ of rank $m-4\ge2$, and
$R_0\perp\psi(V)\subset\Kv_\PP$. By Lemma \pref{lem:odduni}{1}, $Q(R_0)$ contains every element of
$\Ok_\PP$ of even order, and $Q(\psi(V))=\pi Q(W)$ contains every element of $\Ok_\PP$ of odd
order. So $Q(\Kv_\PP)\supset\Ok_\PP$ in this case, too.

Finally $Q(\Kv_\PP)\subset\Ok_\PP$, because $M$ is integral.
\end{proof}

\begin{remark}\label{rem:HIgap}
The hypothesis $m\ge6$ of Lemma \ref{lem:complement6} cannot be dropped. Let $F=\Q$ and
$M=(\Z e+\Z f)\perp\langle1,1\rangle\perp\langle3\rangle$, where $Q(e)=Q(f)=0$ and $B(e,f)=1$, and
let $v:=2e+9f$. Then $m=5$, $\ord_3\vv M=1=m-4$, and $Q(v)=36$. We have
$v^{\perp}\cap(\Z_3e+\Z_3f)=\Z_3w$ with $w:=2e-9f$ and $Q(w)=-36$, and so
\[
\Kv_3=\langle-36\rangle\perp\langle1,1\rangle\perp\langle3\rangle .
\]
This lattice does not represent $6$. Indeed, suppose that $-36a^{2}+b^{2}+c^{2}+3d^{2}=6$ with
$a,b,c,d\in\Z_3$. Then $b^{2}+c^{2}\equiv0\pmod3$, and so $3\mid b$ and $3\mid c$, as $-1$ is a
quadratic nonresidue modulo $3$. Dividing the equation by $3$, we get $d^{2}\equiv2\pmod3$, which
is impossible.
\end{remark}

\begin{lemma}\label{lem:oddT}
Let $\PP$ be odd and put $t_{\PP}:=\ord_{\PP}\ssc M_{\PP}$. There is $v(\PP)\in M_{\PP}$ with
$\ord_{\PP}Q(v(\PP))=t_{\PP}$, and every $v\in M$ with $v\equiv v(\PP)\pmod{\PP M_{\PP}}$ satisfies
$\ord_{\PP}Q(v)=t_{\PP}$, $M_{\PP}=\Ok_\PP v\perp \Kv_{\PP}$, and
$\ord_{\PP}\vv \Kv=\ord_{\PP}\vv M-t_{\PP}$. Moreover, with
$h_{\PP}:=\max\bigl\{0,\lfloor(\ord_{\PP}\vv M-t_{\PP}-1)/(m-4)\rfloor\bigr\}$ we have
$\PP^{h_{\PP}}\subset Q(\Kv_{\PP})$.
\end{lemma}

\begin{proof}
As $\PP$ is odd, $\nn M_{\PP}=\ssc M_{\PP}=\PP^{t_\PP}$ by Lemma \ref{lem:normrange}. Since
$\nn M_\PP$ is generated by the values of $Q$ on $M_\PP$, some $v(\PP)\in M_\PP$ has
$\ord_\PP Q(v(\PP))=t_\PP$. For $v\equiv v(\PP)\pmod{\PP M_\PP}$,
\[
Q(v)-Q(v(\PP))=Q\bigl(v-v(\PP)\bigr)+2B\bigl(v(\PP),v-v(\PP)\bigr)
\in\PP^{2}\nn M_{\PP}+\PP\,\ssc M_{\PP}=\PP^{\,t_\PP+1},
\]
and so $\ord_\PP Q(v)=t_\PP$. The lattice $\Ok_\PP v$ of rank one is therefore
$\PP^{t_\PP}$-modular. As $B(v,M_\PP)\subset\ssc M_\PP=\PP^{t_\PP}$, it splits $M_\PP$
\cite[82:15a]{OM2}, i.e., $M_{\PP}=\Ok_\PP v\perp\Kv_\PP$. Taking volumes, we get
$\vv M_\PP=Q(v)\,\vv\Kv_\PP$ and $\ord_{\PP}\vv\Kv=\ord_{\PP}\vv M-t_\PP$. The last statement is
the second inclusion of Proposition \ref{prop:oddlower}, applied to the integral lattice
$\Kv_{\PP}$ of rank $m-1\ge4$.
\end{proof}

\begin{lemma}\label{lem:q}
Let $\qq$ be an odd prime with $\qq\nmid\vv M$. Then $M_{\qq}$ is unimodular. Fix a nonsquare unit
$\delta_\qq$ of $\Ok_\qq$, put $\Sigma_\qq:=\{1,\delta_\qq,\pi_\qq,\pi_\qq\delta_\qq\}$, and for
$c\in\Sigma_\qq$ choose $x_c\in M_{\qq}$ with $Q(x_c)=c$. Then:
\begin{enumerate}
\item\label{lem:q:1} Every $v\in M$ with $v\equiv x_c\pmod{\qq^{2}M_{\qq}}$ satisfies $Q(v)=c\gamma^{2}$ with
$\gamma\in\Ok_{\qq}^{\times}$, and $\ord_{\qq}\vv \Kv\le1$.
\item\label{lem:q:2} For every nonzero $v\in M$, $Q(\Kv_{\qq})=\Ok_{\qq}$.
\item\label{lem:q:3} For every $\alpha\in\Ok_{\qq}$ there is $c\in\Sigma_\qq$ such that every $v$ as in
\textup{(1)} for that $c$ satisfies $\alpha=Q(v)\beta_v^{2}$ for some $\beta_v\in\Ok_\qq$.
\end{enumerate}
\end{lemma}

\begin{proof}
As $M$ is integral, $\vv M_\qq\subset(\ssc M_\qq)^{m}\subset\Ok_\qq$ \cite[82:10]{OM2}, while
$\qq\nmid\vv M$ gives $\vv M_\qq=\Ok_\qq$. Hence $\ssc M_\qq=\Ok_\qq$, and $M_\qq=\dual{(M_\qq)}$ is
unimodular. In particular $Q(M_\qq)=\Ok_\qq$ by Lemma \pref{lem:odduni}{2}, so the vectors $x_c$
exist.

(1) We have
\[
Q(v)-c=Q(v-x_c)+2B(x_c,v-x_c)\in\qq^{4}\nn M_\qq+\qq^{2}\ssc M_\qq=\qq^{2},
\]
and so $Q(v)/c\in1+\qq^{2-\ord_\qq c}\subset1+\qq$. As $\qq$ is odd, every element of $1+\qq$ is the
square of a unit by Lemma \ref{lem:lst}. Hence $Q(v)=c\gamma^{2}$ with $\gamma\in\Ok_\qq^{\times}$.
In particular $\ord_\qq Q(v)=\ord_\qq c\le1$, and Lemma \ref{lem:vol} gives
$\ord_\qq\vv\Kv\le\ord_\qq\vv M+\ord_\qq Q(v)\le1$.

(2) This is Lemma \ref{lem:complement}, as $\ord_\qq\vv M=0\le m-5$.

(3) If $\alpha=0$, we can take an arbitrary $c\in\Sigma_\qq$ and $\beta_v=0$. Let $\alpha\ne0$.
The four elements of $\Sigma_\qq$ represent the four classes of
$F_\qq^{\times}/(F_\qq^{\times})^{2}$, so $\alpha=c\,\beta_0^{2}$ for some $c\in\Sigma_\qq$ and
$\beta_0\in F_\qq^\times$. Here $\beta_0\in\Ok_\qq$, because $\ord_\qq\alpha\ge0$ and
$\ord_\qq c\le1$. Now let $v$ be as in (1) for this $c$, and write $Q(v)=c\gamma^{2}$ with a unit
$\gamma$, which depends on $v$. Then $\alpha=Q(v)\beta_v^{2}$ with
$\beta_v:=\beta_0\gamma^{-1}\in\Ok_\qq$.
\end{proof}
\section{Genera and spinor genera}\label{sec:genera}

In this section we will collect the global results on genera and spinor genera that the proof of
Theorem \ref{thm:main} uses, and we will define the integers $r_\qq(K)$ that enter its statement.

\begin{proposition}[{\cite[p. 148]{Hs}}]\label{prop:gentospn}
Let $V$ be a quadratic space over $F$, let $L$ be a lattice on $V$, and let $K$ be a lattice whose
ambient space $FK$ is a nondegenerate subspace of $V$ with $\dim V-\dim FK\ge3$. If $K$ is
represented by $\gen(L)$, then $K$ is represented by $\spn^{+}(L)$.
\end{proposition}

Proposition \ref{prop:gentospn} requires no hypothesis at the dyadic primes. Indeed, let $U$ be the
orthogonal complement of $FK$ in $V$. The proof in \cite{Hs} uses spinor norms only through the
equality $\theta\bigl(O^{+}(U_\PP)\bigr)=F_\PP^{\times}$, which holds at every prime $\PP$ because
$\dim U\ge3$ \cite[91:6]{OM2}.

The situation is different in codimension $2$. There spinor exceptional representations occur, and
their description requires integral spinor norm groups. At ramified dyadic primes, the groups
$\theta\bigl(O^{+}(L_\PP)\bigr)$ are hard to compute, and they were determined only in \cite{Be}.
Compare \cite{BH2}, which treats lattices of rank $m-2$ represented by a lattice of rank $m$, and
which assumes the discriminant of the field to be odd where the local spinor norm groups have to be
determined.

\begin{lemma}\label{lem:localtogen}
Let $L$ be an $\Ok$-lattice on a quadratic space $V$ over $F$, and let
$\beta\in\Ok$ be totally positive with $\beta\in Q(L_\PP)$ for every prime $\PP$. Then
$\beta$ is represented by $\gen(L)$.
\end{lemma}

\begin{proof}
The space $V_\PP$ represents $\beta$ at every prime $\PP$ by hypothesis. At every real place $v$,
the space $V_v$ represents $\beta$ because $V$ is positive definite and $\beta$ is totally
positive. By the Hasse--Minkowski theorem, $V$ therefore represents $\beta$, i.e., there is
$y\in V$ with $Q(y)=\beta$. For every prime $\PP$, the lattice $\Ok_\PP y$ is represented by
$L_\PP$, because $\beta\in Q(L_\PP)$. By \cite[102:5]{OM2}, it follows that some lattice in
$\gen(L)$ represents $\Ok y$, and in particular $\beta$.
\end{proof}

\begin{lemma}\label{lem:spinor}
Let $K$ be an $\Ok$-lattice on a quadratic space $V$ over $F$ with $\dim V\ge3$, and let $\qq$ be
a prime of $F$ such that $V_{\qq}$ is isotropic. Then:
\begin{enumerate}
\item\label{lem:spinor:1} Every lattice in $\spn^{+}(K)$ is isometric to a lattice $N\in\spn^{+}(K)$ with
$N_\PP=K_\PP$ for every prime $\PP\ne\qq$.
\item\label{lem:spinor:2} There is an integer $s\ge0$ such that every lattice in $\spn^{+}(K)$ is isometric to a
lattice $N'$ with $\qq^{s}N'\subset K$.
\end{enumerate}
\end{lemma}

\begin{proof}
(1) Let $S$ be the set that consists of the real places of $F$ and of $\qq$. The primes outside
$S$ are the primes $\PP\ne\qq$, and so two lattices $N,N'$ on $V$ satisfy $\Ok_SN=\Ok_SN'$ if and
only if $N_\PP=N'_\PP$ for every prime $\PP\ne\qq$. As $V_{\qq}$ is isotropic, the set $S$ is
indefinite for $V$. Since $\dim V\ge3$, \cite[104:5]{OM2} gives
\[
\cls^{+}(\Ok_SK)=\spn^{+}(\Ok_SK).
\]

Let $N_0\in\spn^{+}(K)$, i.e., $N_0=\sigma_0\Sigma K$ with $\sigma_0\in O^{+}(V)$ and
$\Sigma\in O'_{\mathbb A}(V)$. Then $\Ok_SN_0=\sigma_0\Sigma(\Ok_SK)$, where only the components of
$\Sigma$ at the primes $\PP\ne\qq$ are used. So $\Ok_SN_0$ lies in $\spn^{+}(\Ok_SK)$, and hence in
$\cls^{+}(\Ok_SK)$. This means that there is $\sigma\in O^{+}(V)$ with $\sigma(\Ok_SN_0)=\Ok_SK$. The
lattice $N:=\sigma N_0$ is isometric to $N_0$, it lies in $\spn^{+}(K)$, and it satisfies
$\Ok_SN=\Ok_SK$, i.e., $N_{\PP}=K_{\PP}$ for every prime $\PP\ne\qq$.

(2) By \cite[103:4]{OM2}, the proper spinor genus $\spn^{+}(K)$ is the union of finitely many proper
classes. By (1), we can choose representatives $K_1,\dots,K_g$ of these classes such that
$(K_i)_\PP=K_\PP$ for every $i$ and every prime $\PP\ne\qq$. As $(K_i)_{\qq}$ and $K_{\qq}$ are
lattices on the same space $V_{\qq}$, there is an integer $s\ge0$ with
$\pi_\qq^{s}(K_i)_\qq\subset K_\qq$ for $i=1,\dots,g$. By Convention \ref{conv:scaling}, the lattices
$\qq^{s}K_i$ and $K$ have the same localizations at the primes $\PP\ne\qq$, and therefore
$\qq^{s}K_i\subset K$. This proves (2), because every lattice in $\spn^{+}(K)$ is isometric to one
of the lattices $K_i$.
\end{proof}

\begin{lemma}\label{lem:spinorcrit}
Let $K$ be an $\Ok$-lattice of rank $\ge4$ on a quadratic space $V$ over $F$, let $\qq$ be a prime
of $F$, and let $s\ge0$ be an integer such that every lattice in $\spn^{+}(K)$ is isometric to a
lattice $N'$ with $\qq^{s}N'\subset K$. Then, for every integer $r\ge s$, every nonzero
$\beta\in F$ represented by $\gen(\qq^{r}K)$ is represented by $K$.
\end{lemma}

\begin{proof}
Let $r\ge s$, and let $\beta\in F$ be nonzero and represented by $\gen(\qq^{r}K)$. Thus
$\beta=Q(x)$ for some vector $x$ in a lattice $L'\in\gen(\qq^{r}K)$. The lattice $\Ok x$ of rank one
is contained in $L'$, so it is represented by $\gen(\qq^{r}K)$. Its ambient space $Fx$ is
nondegenerate, because $\beta\ne0$, and $\dim V-\dim Fx\ge3$. So Proposition \ref{prop:gentospn},
applied to the lattice $\qq^{r}K$ on $V$ and to $\Ok x$, shows that some lattice
$L_0\in\spn^{+}(\qq^{r}K)$ represents $\Ok x$. In particular $\beta\in Q(L_0)$.

Write $L_0=\sigma\Sigma(\qq^{r}K)$ with $\sigma\in O^{+}(V)$ and $\Sigma\in O'_{\mathbb A}(V)$. The
lattices $\sigma\Sigma(\qq^{r}K)$ and $\qq^{r}\sigma\Sigma K$ have the same localization at every
prime, so they are equal. Hence $L_0=\qq^{r}N$ with $N:=\sigma\Sigma K\in\spn^{+}(K)$. By the
hypothesis on $s$, there is $\tau\in O(V)$ such that $N':=\tau N$ satisfies $\qq^{s}N'\subset K$. Then
$\qq^{r}N'=\tau(\qq^{r}N)$, and so $\beta\in Q(\qq^{r}N)=Q(\qq^{r}N')$. Moreover
$\qq^{r}N'\subset\qq^{s}N'\subset K$, because $r\ge s$. Therefore $\beta\in Q(K)$.
\end{proof}

Let $K$ be an $\Ok$-lattice of rank $\ge4$, and let $\qq$ be a prime of
$F$ such that $K_\qq$ is isotropic. By Lemmas \pref{lem:spinor}{2} and \ref{lem:spinorcrit}, there
is an integer $r_0\ge0$ such that, for every integer $r\ge r_0$, every totally positive
$\beta\in\Ok$ that is represented by $\gen(\qq^{r}K)$ is represented by $K$.

\begin{definition}\label{def:rK}
For $K$ and $\qq$ as above, $r_\qq(K)$ is the least integer $r_0\ge0$ with this property.
\end{definition}

The integer $r_\qq(K)$ depends only on $\qq$ and on the isometry class of $K$.
Lemmas \pref{lem:spinor}{2} and \ref{lem:spinorcrit} are the analogue of \cite[Lemma 1.2]{HKK},
which assumes $\gen L=\spn L$ in place of Proposition \ref{prop:gentospn}. This proposition is the
reason for the hypothesis $\rank K\ge4$. The proof of Lemma \pref{lem:spinor}{2} gives no bound for
$s$, and we will bound $r_\qq(K)$ by a different argument in \S\ref{sec:r}.

\section{The main theorem}\label{sec:global}

In this section we will prove our main result, Theorem \ref{thm:main}. We begin with the estimates
from the geometry of numbers that we will use to choose the vector $v$ and the element $x$ of
\S\ref{subsec:proofintro}, and to count the lattices $v^{\perp}\cap M$ up to isometry.

\subsection{Geometry of numbers}\label{subsec:gon}
A \emph{lattice in $\R^{n}$} is a discrete subgroup $\Lambda$ of $\R^{n}$ of rank $n$, and
$\covol(\Lambda)$ is the volume of $\R^{n}/\Lambda$. For a norm $\|\cdot\|$ on $\R^{n}$, the
\emph{successive minima} $\lambda_1\le\dots\le\lambda_n$ of $\Lambda$ are defined by letting
$\lambda_k$ be the least $\lambda>0$ such that $\Lambda$ contains $k$ linearly independent vectors of
norm at most $\lambda$.

\begin{theorem}[Minkowski's second theorem {\cite[Ch. VIII]{Ca}}]\label{thm:minkowski2}
Let $\Lambda$ be a lattice of rank $n$ in $\R^{n}$, let $\|\cdot\|$ be a norm on $\R^{n}$ with
unit ball $\mathcal B$, and let $\lambda_1\le\dots\le\lambda_n$ be the successive minima of
$\Lambda$ for $\|\cdot\|$. Then
$\lambda_1\cdots\lambda_n\operatorname{vol}(\mathcal B)\le2^{n}\covol(\Lambda)$.
\end{theorem}

\begin{lemma}\label{lem:covering}
Let $\Lambda$ be a lattice of rank $n$ in $\R^{n}$, let $\|\cdot\|$ be a norm on $\R^{n}$, and let
$b_1,\dots,b_n\in\Lambda$ be linearly independent. Then for every $x\in\R^{n}$ there is
$y\in\Lambda$ with $\|x-y\|\le\tfrac12\sum_{i=1}^{n}\|b_i\|$. In particular, if
$\lambda_1\le\dots\le\lambda_n$ are the successive minima of $\Lambda$ for $\|\cdot\|$, there is
$y\in\Lambda$ with $\|x-y\|\le\tfrac12\sum_{i=1}^{n}\lambda_i$.
\end{lemma}

\begin{proof}
Write $x=\sum_it_ib_i$ with $t_i\in\R$, choose $k_i\in\Z$ with $|t_i-k_i|\le\frac12$, and put
$y:=\sum_ik_ib_i\in\Lambda$. Then $\|x-y\|\le\sum_i|t_i-k_i|\,\|b_i\|\le\frac12\sum_i\|b_i\|$. For the
second assertion, recall that $\Lambda$ contains linearly independent vectors $b_1,\dots,b_n$ with
$\|b_i\|=\lambda_i$, and apply the first assertion to them.
\end{proof}

\begin{lemma}\label{lem:idealcover}
Let $\bb$ be a nonzero integral ideal of $\Ok$. For every $y\in\R^{l}$ there is $b\in\bb$ with
$|y_j-\tau_j(b)|\le\frac l2d_F^{1/2}\Nm(\bb)^{1/l}$ for every $j$. In particular, every coset of
$\bb$ in $\Ok$ contains a representative $x$ with $|\tau_j(x)|\le\frac l2d_F^{1/2}\Nm(\bb)^{1/l}$
for every $j$.
\end{lemma}

\begin{proof}
The map $x\mapsto(\tau_j(x))_j$ embeds $\bb$ into $\R^{l}$ as a lattice of covolume
$\Nm\bb\cdot d_F^{1/2}$. Let $\lambda_1\le\dots\le\lambda_l$ be the successive minima of this
lattice for the sup-norm $|\cdot|_\infty$ on $\R^{l}$. Every nonzero $x\in\bb$ satisfies
$\Nm\bb\le|\NFQ(x)|\le|x|_{\infty}^{l}$, and so $\lambda_1\ge(\Nm\bb)^{1/l}$. The unit ball of the
sup-norm has volume $2^{l}$, so Theorem \ref{thm:minkowski2} gives
$\lambda_1\cdots\lambda_l\le\covol(\bb)$. Therefore
\[
\lambda_l\le\frac{\covol(\bb)}{\lambda_1^{l-1}}\le d_F^{1/2}(\Nm\bb)^{1/l}.
\]
By Lemma \ref{lem:covering}, for every $y\in\R^{l}$ there is $b\in\bb$ whose image in $\R^{l}$ has
sup-distance at most $\tfrac12\sum_i\lambda_i\le\tfrac l2\lambda_l$ from $y$. This is the first
assertion.

For the second assertion, let $x_0$ be an element of the given coset, and let $b\in\bb$ be as in
the first assertion for $y:=(-\tau_j(x_0))_j$. Then $x:=x_0+b$ lies in the coset and satisfies
$|\tau_j(x)|\le\frac l2d_F^{1/2}\Nm(\bb)^{1/l}$ for every $j$.
\end{proof}

\begin{lemma}\label{lem:weighted}
Let $L$ be an integral positive definite $\Ok$-lattice of rank $n$, put
$Y:=\Nm(\vv L)$, let $w_1,\dots,w_l>0$ satisfy $\prod_jw_j=1$, and for $x\in FL$ put
\[
\|x\|_w:=\Bigl(\sum_{j=1}^{l}w_j\tau_j\bigl(Q(x)\bigr)\Bigr)^{1/2}.
\]
Then there are $a_1,\dots,a_n\in L$ such that, for $1\le i\le n$ and
$W_{i-1}:=Fa_1+\dots+Fa_{i-1}$, we have $a_i\notin W_{i-1}$ and $\|a_i\|_w\le\|x\|_w$ for every
$x\in L\smallsetminus W_{i-1}$. Every such sequence $a_1,\dots,a_n$ is a basis of $FL$, and the
lattice $L_0:=\Ok a_1+\dots+\Ok a_n$ has the following properties:
\begin{enumerate}
\item\label{lem:weighted:1} $l\le\|a_1\|_w^{2}\le\dots\le\|a_n\|_w^{2}$.
\item\label{lem:weighted:2} $\prod_{i=1}^{n}\|a_i\|_w^{2}\le(ln)^{n}d_F^{n/l}Y^{1/l}$.
\item\label{lem:weighted:3} $\|a_n\|_w^{2}\le l\,n^{n}d_F^{n/l}Y^{1/l}$.
\item\label{lem:weighted:4} $(L:L_0)\le d_F^{n/2}n^{ln}$.
\end{enumerate}
\end{lemma}

\begin{proof}
For $1\le j\le l$ put $V_j:=FL\otimes_{F,\tau_j}\R$, and let $Q^{\tau_j}$ be the real quadratic
form on $V_j$ obtained from $Q$ by extension of scalars along $\tau_j$. It is positive definite,
because $L$ is positive definite. Put
\[
V_\R:=FL\otimes_\Q\R=\bigoplus_{j=1}^{l}V_j .
\]
The form $\sum_jw_jQ^{\tau_j}$ makes $V_\R$ a Euclidean space of dimension $N:=ln$, in which the
subspaces $V_j$ are mutually orthogonal. We embed $FL$ into $V_\R$ diagonally. The Euclidean norm of
$V_\R$ then restricts to $\|\cdot\|_w$ on $FL$, and we denote it by $\|\cdot\|_w$ as well. For
$\xi\in F$ and $x\in FL$ we have
\begin{equation}\label{eq:xinorm}
\|\xi x\|_w^{2}=\sum_{j=1}^{l}w_j\,\tau_j(\xi)^{2}\,\tau_j(Q(x))\le\max_j|\tau_j(\xi)|^{2}\cdot\|x\|_w^{2}.
\end{equation}

\emph{Existence of $a_1,\dots,a_n$.} The image of $L$ in $V_\R$ is a lattice of rank $N$. In
particular, $L$ contains only finitely many vectors of bounded norm. So for each $i\le n$ the
minimum of $\|x\|_w$ over $x\in L\smallsetminus W_{i-1}$ is attained, and this gives the vectors
$a_1,\dots,a_n$ one after another. They are linearly independent over $F$ by construction, and
therefore they form a basis of $FL$.

\emph{The covolume of $L$.} Suppose first that $L$ is free with basis $e_1,\dots,e_n$, and let
$\omega_1,\dots,\omega_l$ be a $\Z$-basis of $\Ok$. Then $(\omega_ke_i)_{k,i}$ is a $\Z$-basis of
$L$, and its Gram matrix in the Euclidean space $V_\R$ has determinant
\[
\det(\tau_j\omega_k)_{j,k}^{2n}\,\prod_{j=1}^{l}\Bigl(w_j^{n}\det\bigl(\tau_jB(e_i,e_{i'})\bigr)_{i,i'}\Bigr)
=d_F^{n}\,\NFQ\bigl(\det(B(e_i,e_{i'}))_{i,i'}\bigr)=d_F^{n}Y .
\]
Hence $\covol(L)=d_F^{n/2}Y^{1/2}$. In general, $L$ contains a free sublattice $L'$ of finite
index. Then $\covol(L')=(L:L')\covol(L)$ and $\vv L'=[L:L']^{2}\vv L$, and so the formula
$\covol(L)=d_F^{n/2}Y^{1/2}$ holds for every $L$.

(1) Let $x\in L$ be nonzero. Then $Q(x)$ is a nonzero totally positive element of $\Ok$, and the
inequality between the arithmetic and the geometric mean gives
\[
\|x\|_w^{2}\ge l\Bigl(\prod_jw_j\tau_j(Q(x))\Bigr)^{1/l}=l\,\NFQ(Q(x))^{1/l}\ge l .
\]
The sets $L\smallsetminus W_{i-1}$ decrease as $i$ increases, and so
$\|a_1\|_w\le\dots\le\|a_n\|_w$.

(2) Let $\lambda_1\le\dots\le\lambda_N$ be the successive minima of $L$ for $\|\cdot\|_w$, and let
$1\le i\le n$. The real span of $W_{i-1}$ in $V_\R$ has dimension $l(i-1)$, and its intersection
with $L$ is $L\cap W_{i-1}$. Hence, among $l(i-1)+1$ linearly independent vectors of $L$ of norm
at most $\lambda_{l(i-1)+1}$, at least one lies outside $W_{i-1}$. By the choice of $a_i$, this
gives $\|a_i\|_w\le\lambda_{l(i-1)+1}\le\lambda_{l(i-1)+t}$ for $1\le t\le l$. Taking the product
over all $i$ and $t$, we get $\prod_i\|a_i\|_w^{l}\le\prod_k\lambda_k$. The unit ball of
$\|\cdot\|_w$ contains a cube of side $2/\sqrt N$, so its volume is at least $(2/\sqrt N)^{N}$, and
Theorem \ref{thm:minkowski2} gives
\[
\prod_k\lambda_k\le N^{N/2}\covol(L)=N^{N/2}d_F^{n/2}Y^{1/2}.
\]
Raising $\prod_i\|a_i\|_w^{l}\le N^{N/2}d_F^{n/2}Y^{1/2}$ to the power $2/l$, we obtain (2).

(3) By (1), we have $\|a_i\|_w^{2}\ge l$ for $i<n$, and so
$\|a_n\|_w^{2}\le l^{-(n-1)}\prod_i\|a_i\|_w^{2}$. Now (2) gives the claim.

(4) As $a_1,\dots,a_n$ is a basis of $FL$, there is an isomorphism of real vector spaces
$V_\R\to\R^{N}$ that sends $\sum_i\xi_ia_i$ to $(\tau_j(\xi_i))_{i,j}$ for $\xi_1,\dots,\xi_n\in F$.
It maps $L_0$ onto the image of $\Ok^{n}$, which is a lattice of covolume $d_F^{n/2}$. As $L_0$ has
index $(L:L_0)$ in $L$, it maps $L$ onto a lattice $\Lambda$ of covolume $d_F^{n/2}/(L:L_0)$.

Let $x\in L$ be nonzero, and write $x=\sum_{i\le s}\xi_ia_i$ with $\xi_i\in F$ and $\xi_s\ne0$. Then
$x\notin W_{s-1}$, and so $\|x\|_w\ge\|a_s\|_w$ by the choice of $a_s$. On the other hand,
\eqref{eq:xinorm} and (1) give
\[
\|x\|_w\le\sum_{i\le s}\max_j|\tau_j(\xi_i)|\,\|a_i\|_w\le s\,\max_{i,j}|\tau_j(\xi_i)|\,\|a_s\|_w .
\]
Hence $\max_{i,j}|\tau_j(\xi_i)|\ge1/s\ge1/n$. In other words, every nonzero point of $\Lambda$ has
sup-norm at least $1/n$, so all successive minima of $\Lambda$ for the norm $n|\cdot|_\infty$ are
at least $1$. The unit ball of this norm has volume $(2/n)^{N}$, and Theorem \ref{thm:minkowski2}
gives $(2/n)^{N}\le2^{N}d_F^{n/2}/(L:L_0)$. This is (4), as $N=ln$.
\end{proof}

\begin{proposition}\label{prop:coset}
Let $L$ be an integral positive definite $\Ok$-lattice of rank $n$, put
$Y:=\Nm(\vv L)$, let $w_1,\dots,w_l>0$ satisfy $\prod_jw_j=1$, let $\bb$ be a nonzero integral
ideal of $\Ok$, and put
\[
C:=\tfrac14\,l^{3}n^{n+2}d_F^{1+n/l}Y^{1/l}\Nm(\bb)^{2/l}.
\]
Every coset $v_0+\bb L$ with $v_0\in L$ contains a vector $v$ with
$\sum_jw_j\tau_j(Q(v))\le C$. In particular
$w_j\tau_j(Q(v))\le C$ for every $j$, and $\NFQ(Q(v))\le(C/l)^{l}$.
\end{proposition}

\begin{proof}
Let $a_1,\dots,a_n$ be as in Lemma \ref{lem:weighted}, write $v_0=\sum_ic_ia_i$ with $c_i\in F$,
and choose by Lemma \ref{lem:idealcover} elements $b_i\in\bb$ with
$|\tau_j(c_i-b_i)|\le\frac l2d_F^{1/2}\Nm(\bb)^{1/l}$ for all $i$ and $j$. The vector
$v:=v_0-\sum_ib_ia_i=\sum_i(c_i-b_i)a_i$ lies in $v_0+\bb L$. By \eqref{eq:xinorm} and
Lemma \pref{lem:weighted}{1}, we have
\[
\|v\|_w\le\sum_i\max_j|\tau_j(c_i-b_i)|\,\|a_i\|_w\le\tfrac l2d_F^{1/2}\Nm(\bb)^{1/l}\,n\,\|a_n\|_w ,
\]
and so $\|v\|_w^{2}\le C$ by Lemma \pref{lem:weighted}{3}. This is the first assertion, and it implies
that $w_j\tau_j(Q(v))\le C$ for every $j$. Finally, the inequality between the arithmetic and the
geometric mean gives $\NFQ(Q(v))=\prod_jw_j\tau_j(Q(v))\le(C/l)^{l}$, as $\prod_jw_j=1$.
\end{proof}

\begin{lemma}\label{lem:box}
For $\Upsilon>0$, $\#\{x\in\Ok\mid|\tau_j(x)|\le \Upsilon\text{ for all }j\}\le(2\Upsilon+1)^{l}$.
\end{lemma}

\begin{proof}
We identify $x\in\Ok$ with the point $(\tau_j(x))_j\in\R^{l}$. A nonzero $x\in\Ok$ satisfies
$|\NFQ(x)|\ge1$, and therefore $|x|_{\infty}\ge1$. So two distinct points of $\Ok$ have sup-distance
at least $1$, and the open cubes of side $1$ centered at the points that we want to count are
pairwise disjoint. Each of these cubes has volume $1$, and all of them lie in the cube
$[-\Upsilon-\tfrac12,\Upsilon+\tfrac12]^{l}$, whose volume is $(2\Upsilon+1)^{l}$.
\end{proof}

\begin{lemma}\label{lem:overlattices}
Let $N\ge1$ and $t\ge1$ be integers. Then $\Q^{N}$ has at most $t^{N}$ subgroups $\Lambda$ with
$\Z^{N}\subset\Lambda$ and $(\Lambda:\Z^{N})=t$.
\end{lemma}

\begin{proof}
Let $b_N(t)$ be the number of these subgroups. We will show by induction on $N$ that
$b_N(t)\le t^{N}$. For $N=1$ the only such subgroup is $t^{-1}\Z$, and so $b_1(t)=1$.

Let $N\ge2$, and let $\Lambda$ be a subgroup of $\Q^{N}$ with $\Z^{N}\subset\Lambda$ and
$(\Lambda:\Z^{N})=t$. Let $p\colon\Q^{N}\to\Q$ be the projection to the last coordinate, and identify
$\Q^{N-1}$ with the kernel of $p$. The group $\Lambda$ is finitely generated, because $\Z^{N}$ has finite index in it. So
$p(\Lambda)$ is a finitely generated subgroup of $\Q$ that contains $\Z$, and therefore
$p(\Lambda)=d^{-1}\Z$ for an integer $d\ge1$. Put
$\Lambda':=\Lambda\cap\Q^{N-1}$. Then $\Z^{N-1}\subset\Lambda'$ and
\[
t=(\Lambda:\Z^{N})=\bigl(p(\Lambda):p(\Z^{N})\bigr)\cdot(\Lambda':\Z^{N-1})=d\cdot(\Lambda':\Z^{N-1}).
\]
So $d$ divides $t$, and $\Lambda'$ is one of the $b_{N-1}(t/d)$ subgroups of $\Q^{N-1}$ that contain
$\Z^{N-1}$ with index $t/d$.

Choose $\gamma\in\Lambda$ with $p(\gamma)=1/d$, and write $\gamma=(u,1/d)$ with $u\in\Q^{N-1}$. Then
$\Lambda=\Lambda'+\Z\gamma$. Moreover, the element $d\gamma-(0,\dots,0,1)$ lies in $\Lambda$ and in
the kernel of $p$, and so $du\in\Lambda'$. Hence $\Lambda$ is determined by $d$, by $\Lambda'$, and
by the class of $u$ in $d^{-1}\Lambda'/\Lambda'$, which is a group of order $d^{N-1}$. Together with
the induction hypothesis, this gives
\[
b_N(t)\le\sum_{d\mid t}d^{N-1}\,b_{N-1}(t/d)\le\sum_{d\mid t}d^{N-1}(t/d)^{N-1}=\tau(t)\,t^{N-1}\le t^{N},
\]
where $\tau(t)\le t$ is the number of positive divisors of $t$.
\end{proof}

\begin{proposition}\label{prop:classes}
Let $n\ge1$ and $Y\ge1$, and put $I_n:=d_F^{n/2}n^{ln}$, $U_n(Y):=n^{n}d_F^{n/l}Y^{1/l}$, and
\[
\mathcal H_n(Y):=2\,I_n^{\,ln+1}\,(6l)^{ln(n+1)/2}\,U_n(Y)^{l(n+1)/2}
\bigl(1+\log_2U_n(Y)\bigr)^{n-1}.
\]
The integral positive definite $\Ok$-lattices of rank $n$ whose volume ideal has norm at
most $Y$ fall into at most $\mathcal H_n(Y)$ isometry classes.
\end{proposition}

\begin{proof}
Let $L$ be such a lattice. We apply Lemma \ref{lem:weighted} with $w_1=\dots=w_l=1$, and we let
$a_1,\dots,a_n$ and $L_0$ be as in that lemma. Put
$b_i:=\|a_i\|_w^{2}=\operatorname{Tr}_{F/\Q}Q(a_i)$. By Lemma \pref{lem:weighted}{1} and \hyperref[lem:weighted:2]{(2)}, we have
$b_i\ge l$ and $\prod_i(b_i/l)\le U_n(Y)$. For each $i$ let $k_i\ge0$ be the integer with
$2^{k_i}\le b_i/l<2^{k_i+1}$. Then $2^{k_1+\dots+k_n}\le U_n(Y)$, i.e.,
$k_1+\dots+k_n\le u:=\lfloor\log_2U_n(Y)\rfloor$.

We will first count the possible Gram matrices of $a_1,\dots,a_n$. Their entries are
$g_{ii'}:=B(a_i,a_{i'})\in\Ok$. As $L$ is positive definite, the Cauchy--Schwarz inequality
holds in every real embedding, and it gives
\[
|\tau_j(g_{ii'})|\le\bigl(\tau_j(Q(a_i))\,\tau_j(Q(a_{i'}))\bigr)^{1/2}\le(b_ib_{i'})^{1/2}
<R_{ii'}:=2l\cdot2^{(k_i+k_{i'})/2}
\]
for every $j$. For a given tuple $(k_1,\dots,k_n)$, Lemma \ref{lem:box} therefore bounds the number
of these symmetric matrices by
\[
\prod_{i\le i'}(2R_{ii'}+1)^{l}\le\prod_{i\le i'}(3R_{ii'})^{l}=(6l)^{ln(n+1)/2}\,2^{c(k_1+\dots+k_n)},
\qquad\text{where}\quad c:=l(n+1)/2\ge1 .
\]
For the last equality, note that each $k_i$ occurs $n+1$ times in the sum
$\sum_{i\le i'}(k_i+k_{i'})$. There are at most $(s+1)^{n-1}$ tuples $(k_1,\dots,k_n)$ of nonnegative
integers with $k_1+\dots+k_n=s$, and $\sum_{s=0}^{u}2^{cs}\le2\cdot2^{cu}$ because $c\ge1$. So the
number of possible Gram matrices is at most
\begin{align*}
(6l)^{ln(n+1)/2}\sum_{s=0}^{u}(s+1)^{n-1}2^{cs}
&\le2\,(6l)^{ln(n+1)/2}(u+1)^{n-1}2^{cu}\\
&\le2\,(6l)^{ln(n+1)/2}\bigl(1+\log_2U_n(Y)\bigr)^{n-1}U_n(Y)^{c}.
\end{align*}

Next we count the lattices $L$ that belong to a given Gram matrix. The Gram matrix determines $L_0$
up to isometry, and $L$ is a lattice with $L_0\subset L\subset FL_0$ whose index $t:=(L:L_0)$ is at
most $I_n$ by Lemma \pref{lem:weighted}{4}. As a $\Z$-module, $L_0$ is free of rank $ln$, and $L$
is a subgroup of $FL_0=\Q L_0$ that contains $L_0$ with index $t$. By
Lemma \ref{lem:overlattices}, applied with $N=ln$, there are at most $t^{ln}$ such subgroups when
$L_0$ and $t$ are given. Summing over $t\le I_n$, we see that for a given $L_0$ there are at most
$I_n^{ln+1}$ possibilities for $L$. The product of this number and the bound for the number of Gram
matrices is $\mathcal H_n(Y)$.
\end{proof}

\subsection{Setup}\label{subsec:setup}
For the rest of this section, $M$ is a positive definite integral $\Ok$-lattice of rank
$m\ge5$, $D:=\Nm(\vv M)$, and $\qq$ is an odd prime of $F$ with $\qq\nmid\vv M$. As in
\S\ref{sec:otherprimes}, $\Kv:=v^{\perp}\cap M$ for $v\in M$.

Let $A\in\Ok$ be totally positive and represented by $M_\PP$ for every prime $\PP$. To represent
$A$ by $M$, we will choose $v\in M$ in a prescribed class modulo an ideal $\aaa$ and $x\in\Ok$ in a
prescribed class modulo an ideal $\hh$ such that $\beta:=A-Q(v)x^{2}$ is totally positive and is
represented by $\qq^{r}\Kv$ at every prime, for a suitable integer $r\ge0$. It will follow that
$\beta$ is represented by $\Kv$, and so $A=Q(xv)+\beta$ is represented by $M$. Let us define
$\aaa$, $r$, and $\hh$.

\emph{The parts of $D$ and the set $T$.} For a prime $\PP$ put $w_\PP:=\ord_\PP\vv M$. We write
$\Dp:=\prod_{\pp\mid2}(\Nm\pp)^{w_\pp}$ for the dyadic part of $D$, and $\Dodd:=D/\Dp$ for its odd
part. Thus $\Dp$ is the largest power of $2$ that divides the integer $D$, and $\Dp$ and $\Dodd$
are determined by $D$. Put $d_5:=1$ and $d_m:=m-3$ for $m\ge6$, and let $T$ be the set of odd primes $\PP\ne\qq$
with $w_\PP\ge d_m$.

\emph{Local exponents.} For $\pp\mid2$ put
\[
\kappa_\pp:=\Bigl\lfloor\frac{w_\pp}{m-3}\Bigr\rfloor+2e_\pp+1,\qquad n_\pp:=\kappa_\pp+e_\pp+2,
\qquad b_\pp:=\Bigl\lfloor\frac{w_\pp}{m-4}\Bigr\rfloor+e_\pp+1.
\]
These are the numbers $\kappa$ and $n_0$ of Proposition \ref{prop:covering} for the lattice
$M_\pp$ of rank $n=m$, and the number $b$ of Lemma \ref{lem:dyadscalar} for the prime $\pp$. For $\PP\in T$ put
$t_\PP:=\ord_\PP\ssc M_\PP$ and $h_\PP:=\max\{0,\lfloor(w_\PP-t_\PP-1)/(m-4)\rfloor\}$, as in
Lemma \ref{lem:oddT}.

\emph{The ideal $\aaa$ and the bound $X$.} We put
\[
\aaa:=\prod_{\pp\mid2}\pp^{n_\pp}\cdot\qq^{2}\cdot\prod_{\PP\in T}\PP ,\qquad
X:=V_{m,F}\,(\Nm\qq)^{4}D^{2}\Dp^{2/(m-3)}\Dodd^{2/d_m},
\]
where
\begin{equation}\label{eq:V}
V_{m,F}:=\bigl(2^{10}l^{2}m^{m+2}\bigr)^{l}d_F^{m+l}.
\end{equation}

\emph{The integer $r$ and the ideal $\hh$.} Let $\mathcal K$ be the set of isometry classes of
integral positive definite $\Ok$-lattices $K$ of rank $m-1$ with $\Nm(\vv K)\le X$ and
$\ord_\qq\vv K\le1$. For such $K$, the unimodular component of a Jordan splitting of $K_\qq$ has
rank at least $m-2\ge3$, because $\ord_\qq\vv K\le1$. As $\qq$ is odd, $K_\qq$ is therefore isotropic
by Lemma \pref{lem:odduni}{2}, and so $r_\qq(K)$ is defined (Definition \ref{def:rK}). By
Proposition \ref{prop:classes}, $\mathcal K$ is a finite set with at most $\mathcal H_{m-1}(X)$
elements. We put
\begin{equation}\label{eq:r}
r:=\max\bigl(\{0\}\cup\{r_\qq(K)\mid K\in\mathcal K\}\bigr)
\end{equation}
and, finally,
\[
\hh:=\prod_{\pp\mid2}\pp^{b_\pp}\cdot\qq^{2r}\cdot\prod_{\PP\in T}\PP^{h_\PP}.
\]

\subsection{The main theorem}

\begin{theorem}\label{thm:main}
Let $M$, $D$, $\qq$, $\aaa$, $r$, and $\hh$ be as in \S\textup{\ref{subsec:setup}}, put
$\Gamma_{m,F}:=\bigl(l^{5}m^{m+2}/16\bigr)^{l}d_F^{m+2l}$, and let $A\in\Ok$ be totally positive
and represented by $M_{\PP}$ for every prime $\PP$. If
$\NFQ(A)>\Gamma_{m,F}\,D\,\Nm(\aaa)^{2}\Nm(\hh)^{2}$, then $A$ is represented by $M$.
\end{theorem}

\begin{proof}
\emph{Choice of $v$.} We apply Proposition \ref{prop:covering} at each $\pp\mid2$ to
$\alpha:=A\in Q(M_\pp)$ and obtain classes $[y_\pp]\in\mathcal Y(\pp)$. We also apply
Lemma \pref{lem:q}{3} to $\alpha:=A$ and obtain an element $c_0\in\Sigma_\qq$. By the Chinese
remainder theorem there is $v_0\in M$ with $v_0\equiv y_\pp\pmod{\pp^{n_\pp}M_\pp}$ for $\pp\mid2$,
$v_0\equiv x_{c_0}\pmod{\qq^{2}M_\qq}$, and $v_0\equiv v(\PP)\pmod{\PP M_\PP}$ for $\PP\in T$, where
$v(\PP)$ is the vector of Lemma \ref{lem:oddT}. By Proposition \ref{prop:coset}, applied to $L=M$,
$\bb=\aaa$, and the weights $\NFQ(A)^{1/l}/\tau_j(A)$, whose product is $1$, there is
$v\in v_0+\aaa M$ with
\begin{equation}\label{eq:vsmall}
\begin{gathered}
\tau_j(Q(v))\le\frac{\tau_j(A)}{\NFQ(A)^{1/l}}\cdot\frac{l^{3}m^{m+2}}{4}\,d_F^{1+m/l}D^{1/l}
\Nm(\aaa)^{2/l}\qquad(1\le j\le l),\\
\NFQ(Q(v))\le\Bigl(\frac{l^{2}m^{m+2}}{4}\Bigr)^{l}d_F^{m+l}D\,\Nm(\aaa)^{2}.
\end{gathered}
\end{equation}

Since $v\equiv v_0\pmod{\aaa M}$, the vector $v$ satisfies the local conditions imposed on $v_0$.
Hence $A=Q(v)c_\pp^{2}$ for some $c_\pp\in\Ok_\pp$ at each $\pp\mid2$, by
Proposition \ref{prop:covering}. At $\qq$ we have $A=Q(v)c_\qq^{2}$ for some $c_\qq\in\Ok_\qq$ by
Lemma \pref{lem:q}{3}, and $Q(v)\ne0$ by Lemma \pref{lem:q}{1}. At each $\PP\in T$, we have
$M_\PP=\Ok_\PP v\perp\Kv_\PP$ by Lemma \ref{lem:oddT}, so that $A=Q(v)c_\PP^{2}+Q(y_\PP)$ for some
$c_\PP\in\Ok_\PP$ and $y_\PP\in\Kv_\PP$.

\emph{The lattice $K:=\Kv$.} It is integral and positive definite of rank $m-1$, and
$\ord_\qq\vv K\le1$ by Lemma \pref{lem:q}{1}. By Lemma \ref{lem:vol},
$\Nm(\vv K)\le D\,\NFQ(Q(v))$. To bound the right-hand side, we estimate
\[
\Nm(\aaa)=2^{\sum_{\pp\mid2}f_\pp n_\pp}\,(\Nm\qq)^{2}\prod_{\PP\in T}\Nm\PP .
\]
As $n_\pp=\lfloor w_\pp/(m-3)\rfloor+3e_\pp+3$, $\sum_{\pp\mid2}f_\pp e_\pp=l$, and
$\Dp=2^{\sum_{\pp\mid2}f_\pp w_\pp}$, we have
\[
\sum_{\pp\mid2}f_\pp n_\pp\le\frac{1}{m-3}\sum_{\pp\mid2}f_\pp w_\pp+3\sum_{\pp\mid2}f_\pp e_\pp
+3\sum_{\pp\mid2}f_\pp\le\frac{\log_2\Dp}{m-3}+6l .
\]
Moreover $\Nm\PP\le\Nm\PP^{\,w_\PP/d_m}$ for $\PP\in T$, because $w_\PP\ge d_m$, and so
$\prod_{\PP\in T}\Nm\PP\le\Dodd^{1/d_m}$. Hence
\begin{equation}\label{eq:anorm}
\Nm(\aaa)^{2}\le2^{12l}(\Nm\qq)^{4}\Dp^{2/(m-3)}\Dodd^{2/d_m},
\end{equation}
and \eqref{eq:vsmall} gives $\Nm(\vv K)\le X$. So the class of $K$ lies in $\mathcal K$, and
$r_\qq(K)\le r$.

\emph{Choice of $x$.} By the Chinese remainder theorem, the conditions
$x\equiv c_\pp\pmod{\pp^{b_\pp}}$ for $\pp\mid2$, $x\equiv c_\qq\pmod{\qq^{2r}}$, and
$x\equiv c_\PP\pmod{\PP^{h_\PP}}$ for $\PP\in T$ define a coset of $\hh$ in $\Ok$. By
Lemma \ref{lem:idealcover}, we can choose $x$ in this coset with
$|\tau_j(x)|\le\frac l2d_F^{1/2}\Nm(\hh)^{1/l}$ for every $j$. Put $\beta:=A-Q(v)x^{2}$. We will
show that $\beta$ is totally positive and that $\beta\in Q\bigl((\qq^{r}K)_\PP\bigr)$ for every
prime $\PP$. Recall that $(\qq^{r}K)_\qq=\pi_\qq^{r}K_\qq$ and $(\qq^{r}K)_\PP=K_\PP$ for
$\PP\ne\qq$.

\emph{The real places.} By \eqref{eq:vsmall}, the choice of $x$, and the definition of
$\Gamma_{m,F}$, we have
\[
\tau_j(Q(v)x^{2})\le\tau_j(A)\Bigl(\frac{\Gamma_{m,F}\,D\,\Nm(\aaa)^{2}\Nm(\hh)^{2}}{\NFQ(A)}\Bigr)^{1/l}
<\tau_j(A)
\]
for every $j$, where the last inequality is the hypothesis on $\NFQ(A)$. So $\beta$ is totally
positive.

\emph{The primes $\pp\mid2$.} We have $\beta=Q(v)(c_\pp^{2}-x^{2})$ and
$x\equiv c_\pp\pmod{\pp^{b_\pp}}$, and so $\beta\in Q(K_\pp)$ by Lemma \ref{lem:dyadscalar}.

\emph{The prime $\qq$.} We have $\beta=Q(v)(c_\qq-x)(c_\qq+x)\in\pi_\qq^{2r}\Ok_\qq$, because
$x\equiv c_\qq\pmod{\qq^{2r}}$. As $Q(K_\qq)=\Ok_\qq$ by Lemma \pref{lem:q}{2}, we have
$\pi_\qq^{2r}\Ok_\qq=Q(\pi_\qq^{r}K_\qq)=Q\bigl((\qq^{r}K)_\qq\bigr)$.

\emph{The primes $\PP\in T$.} We have $\beta=Q(y_\PP)+Q(v)(c_\PP^{2}-x^{2})$, and the second term
lies in $\PP^{h_\PP}$, because $x\equiv c_\PP\pmod{\PP^{h_\PP}}$. If $\ord_\PP Q(y_\PP)<h_\PP$, then
$\beta=Q(y_\PP)(1+\eta)$ with $\eta\in\PP$. By Lemma \ref{lem:lst}, $1+\eta=\gamma^{2}$ for some
$\gamma\in\Ok_\PP$, and so $\beta=Q(\gamma y_\PP)\in Q(K_\PP)$. If $\ord_\PP Q(y_\PP)\ge h_\PP$, then
$\beta\in\PP^{h_\PP}\subset Q(K_\PP)$ by Lemma \ref{lem:oddT}.

\emph{The odd primes $\PP\notin T$ with $\PP\ne\qq$.} For these primes $w_\PP<d_m$, i.e.,
$w_\PP\le m-5$ if $m=5$ and $w_\PP\le m-4$ if $m\ge6$. So $Q(K_\PP)=\Ok_\PP$ by
Lemma \ref{lem:complement} if $m=5$ and by Lemma \ref{lem:complement6} if $m\ge6$, and
$\beta\in\Ok\subset Q(K_\PP)$.

\emph{Conclusion.} The lattice $\qq^{r}K$ is positive definite, and $\beta$ is totally
positive and lies in $Q\bigl((\qq^{r}K)_\PP\bigr)$ for every prime $\PP$. By
Lemma \ref{lem:localtogen}, $\beta$ is represented by $\gen(\qq^{r}K)$. As $r\ge r_\qq(K)$, it
follows that $\beta\in Q(K)$. Therefore $A=Q(xv)+\beta\in Q(\Ok v\perp K)\subset Q(M)$.
\end{proof}

\begin{corollary}\label{cor:exponent}
In the situation of Theorem \textup{\ref{thm:main}}, $A$ is represented by $M$ if
\[
\NFQ(A)\ >\ G_{m,F}\,(\Nm\qq)^{4(1+r)}\,D\,\Dodd^{\,4/(m-3)}\,\Dp^{\,\frac{2}{m-3}+\frac{2}{m-4}},
\]
with $G_{m,F}$ as in \eqref{eq:G}, and in particular if
$\NFQ(A)>G_{m,F}\,(\Nm\qq)^{4(1+r)}D^{\,1+\frac{2}{m-3}+\frac{2}{m-4}}$.
\end{corollary}

\begin{proof}
The $\qq$-part of $\Nm(\aaa)^{2}\Nm(\hh)^{2}$ is $(\Nm\qq)^{4(1+r)}$. Its dyadic part is
$2^{2\sum_\pp f_\pp(n_\pp+b_\pp)}$, and
\[
\sum_{\pp\mid2}f_\pp b_\pp\le\frac{\log_2\Dp}{m-4}+\sum_{\pp\mid2}f_\pp e_\pp+\sum_{\pp\mid2}f_\pp
\le\frac{\log_2\Dp}{m-4}+2l .
\]
With the estimate preceding \eqref{eq:anorm}, the dyadic part is therefore at most
$2^{16l}\Dp^{2/(m-3)+2/(m-4)}$.

A prime $\PP\in T$ contributes the factor $\Nm\PP^{\,2+2h_\PP}$ to
$\Nm(\aaa)^{2}\Nm(\hh)^{2}$. Put $w:=w_\PP$. Then $w\ge d_m\ge1$ and
$h_\PP\le\lfloor(w-1)/(m-4)\rfloor$. We claim that
\[
2+2\Bigl\lfloor\frac{w-1}{m-4}\Bigr\rfloor\le\frac{4w}{m-3}.
\]
For $m=5$ both sides equal $2w$. For $m\ge6$ put $a:=m-4\ge2$, so that $w\ge d_m=a+1$. The
inequality $2+2(w-1)/a\le4w/(a+1)$ is equivalent to $(a-1)(a+1)\le(a-1)w$, and this holds. By the
claim, the primes in $T$ together contribute at most
$\prod_{\PP\in T}\Nm\PP^{\,4w_\PP/(m-3)}\le\Dodd^{4/(m-3)}$.

As $2^{16l}\Gamma_{m,F}=G_{m,F}$, we conclude that $\Gamma_{m,F}D\,\Nm(\aaa)^{2}\Nm(\hh)^{2}$ is at
most the right-hand side of the first inequality. The second inequality implies the first, since
$\frac{4}{m-3}\le\frac{2}{m-3}+\frac{2}{m-4}$.
\end{proof}

\begin{corollary}\label{cor:nonclassical}
Let $M$ be a positive definite $\Ok$-lattice of rank $m\ge5$ with $\nn M\subset\Ok$, put
$M':=M^{2}$ and $\widehat D:=\Nm(\vv M')=\Nm\bigl(2^{m}\vv M\bigr)$, let $\qq$ be an odd prime of $F$
with $\qq\nmid\vv M'$, let $r$ be the integer \eqref{eq:r} for the lattice $M'$ and the prime $\qq$,
and let $G_{m,F}$ be as in \eqref{eq:G}. If $A\in\Ok$ is totally positive, is represented by
$M_{\PP}$ for every prime $\PP$, and
\[
\NFQ(A)\ >\ 2^{-l}\,G_{m,F}\,(\Nm\qq)^{4(1+r)}\,\widehat D^{\,1+\frac{2}{m-3}+\frac{2}{m-4}},
\]
then $A$ is represented by $M$.
\end{corollary}

\begin{proof}
By Convention \ref{conv:scaling} and Lemma \ref{lem:normrange},
$\ssc M'=2\,\ssc M\subset\nn M\subset\Ok$, so $M'$ is an integral positive definite
$\Ok$-lattice of rank $m\ge5$, and $\vv M'=2^{m}\vv M$. Writing $Q$ for the form of $M$, the form
of $M'$ is $2Q$, so the value set of $M'$ is $2Q(M)$ and that of $M'_{\PP}$ is $2Q(M_{\PP})$ at every
prime $\PP$. Hence $2A$ is totally positive and represented by $M'_{\PP}$ at every prime, and $2A$
is represented by $M'$ if and only if $A$ is represented by $M$. Corollary \ref{cor:exponent},
applied to $M'$ and $2A$, therefore gives $A\in Q(M)$ as soon as $\NFQ(2A)=2^{l}\NFQ(A)$ exceeds
$G_{m,F}(\Nm\qq)^{4(1+r)}\widehat D^{\,1+\frac{2}{m-3}+\frac{2}{m-4}}$.
\end{proof}

\section{An explicit bound for \texorpdfstring{$r$}{r}}\label{sec:r}

The only inexplicit ingredient of Theorem \ref{thm:main} is the integer $r$ of \eqref{eq:r}. In
this section we will bound $r_\qq(K)$, for $K$ in a class of $\mathcal K$, by the class number of
$\gen(K)$ (Lemma \ref{lem:rbound}). The proof will use Kneser's neighbors at the prime $\qq$,
which we will develop in \S\S\ref{subsec:nbr} and \ref{subsec:kneser} in the form that we need. At the
end of the section we will deduce Theorem \ref{cor:crude}.

\subsection{Neighbors along a hyperbolic plane}\label{subsec:nbr}

The prime $\qq$ may divide the volume of the lattice $\Kv$ in the proof of
Theorem \ref{thm:main}. Indeed, if $\ord_\qq A$ is odd, then the vector $v$ of that proof satisfies
$\ord_\qq Q(v)=1$. Then $v$ is primitive in the unimodular lattice $M_\qq$, so
$B(v,M_\qq)=\Ok_\qq$, and Lemma \ref{lem:vol} gives $\ord_\qq\vv\Kv=1$. At a prime
that divides the volume, the definition of a $\qq$-neighbor by index conditions alone does not
preserve the genus in general, and
\cite[Theorem 5.8]{GV} assumes $\vv K$ squarefree, $\#\mathrm{Cl}(\Ok)$ odd, and $\qq\nmid2\vv K$.
We therefore use the following variant, in which a neighbor lies in the same genus by definition.
It is the adaptation indicated in \cite[\S4]{Vo} and carried out over $\Z$ in \cite[(v)]{SP}.

In this section a hyperbolic plane $H$ in an integral $\Ok_\qq$-lattice $L$ always comes with an
ordered basis $(e,f)$ as in \S\ref{sec:notation}, i.e., $Q(e)=Q(f)=0$ and $B(e,f)=1$, and we
write $H=\Ok_\qq e\oplus\Ok_\qq f$. Being unimodular, $H$ splits $L$ \cite[82:15]{OM2}, and we
write $L=H\perp W$ with $W=H^{\perp}\cap L$.

\begin{definition}\label{def:nbr}
\begin{enumerate}
\item\label{def:nbr:1} Let $L$ be an integral $\Ok_\qq$-lattice and $H=\Ok_\qq e\oplus\Ok_\qq f$ a hyperbolic plane
in $L$, with $L=H\perp W$. The \emph{$\qq$-neighbor of $L$ along $H$} is
$L(\qq,H):=\Ok_\qq\,\pi_\qq^{-1}e\oplus\Ok_\qq\,\pi_\qq f\perp W$.
\item\label{def:nbr:2} Let $K$ be an integral $\Ok$-lattice on a quadratic space $V$ over $F$ and $H$ a hyperbolic
plane in $K_\qq$. The \emph{$\qq$-neighbor of $K$ along $H$} is the lattice $K(\qq,H)$ on $V$
with $K(\qq,H)_\qq=K_\qq(\qq,H)$ and $K(\qq,H)_\PP=K_\PP$ for $\PP\ne\qq$.
\end{enumerate}
\end{definition}

The lattice $K(\qq,H)$ of Definition \pref{def:nbr}{2} exists by \cite[81:14]{OM2}, since its family of localizations differs from
that of $K$ at one prime only.

\begin{lemma}\label{lem:nbrbasic}
Let $K'=K(\qq,H)$ be as in Definition \textup{\ref{def:nbr}}, with
$H=\Ok_\qq e\oplus\Ok_\qq f$. Then:
\begin{enumerate}
\item\label{lem:nbrbasic:1} $K'$ is integral, $K'_{\qq}\cong K_{\qq}$, and $K'_\PP=K_\PP$ for $\PP\ne\qq$. Hence
$K'\in\gen(K)$ and $\vv K'=\vv K$.
\item\label{lem:nbrbasic:2} $\qq K'\subset K$ and $\qq K\subset K'$, and
$[K:K\cap K']=[K':K\cap K']=\qq$.
\item\label{lem:nbrbasic:3} If $\sigma\in O(V)$ then $\sigma K'=(\sigma K)(\qq,\sigma H)$.
\item\label{lem:nbrbasic:4} $K=K'(\qq,H')$ for the hyperbolic plane $H':=\Ok_\qq(\pi_\qq f)\oplus\Ok_\qq(\pi_\qq^{-1}e)$
of $K'_\qq$, with the ordered basis $(\pi_\qq f,\pi_\qq^{-1}e)$.
\end{enumerate}
\end{lemma}

\begin{proof}
(1) Put $e':=\pi_\qq^{-1}e$ and $f':=\pi_\qq f$. Then $Q(e')=Q(f')=0$ and $B(e',f')=B(e,f)=1$, so
$\Ok_\qq e'\oplus\Ok_\qq f'$ is again a hyperbolic plane, and
$K'_\qq=(\Ok_\qq e'\oplus\Ok_\qq f')\perp W\cong H\perp W=K_\qq$. Hence $K'$ is integral and lies
in $\gen(K)$ (see \S\ref{sec:notation}).

(2) Every element of $K_\qq$ can be written as $ae+bf+w$ with $a,b\in\Ok_\qq$ and $w\in W$, and
such an element lies in $K'_\qq$ if and only if $b\in\qq$. So
$K_\qq\cap K'_\qq=\Ok_\qq e\oplus\qq f\perp W$. This lattice has index $\qq$ in $K_\qq$. It also has
index $\qq$ in $K'_\qq$, because $K'_\qq/(K_\qq\cap K'_\qq)$ is generated by the image of $e'$, and
$\pi_\qq e'=e$. As $K_\PP=K'_\PP$ for $\PP\ne\qq$, this gives the two index statements. Moreover
\[
\qq K'_\qq=\Ok_\qq e\oplus\qq^{2}f\perp\qq W\subset K_\qq\qquad\text{and}\qquad
\qq K_\qq=\qq e\oplus\qq f\perp\qq W\subset K'_\qq ,
\]
and for $\PP\ne\qq$ we have $(\qq K')_\PP=K'_\PP=K_\PP=(\qq K)_\PP$. Hence $\qq K'\subset K$ and
$\qq K\subset K'$.

(3) An isometry $\sigma\in O(V)$ maps the hyperbolic plane $H$ of $K_\qq$ to the hyperbolic plane
$\sigma H$ of $(\sigma K)_\qq$, with the ordered basis $(\sigma e,\sigma f)$, and the assertion
follows from Definition \ref{def:nbr}.

(4) We have $K'_\qq=H'\perp W$, and $H'$ is a hyperbolic plane with the ordered basis
$(\pi_\qq f,\pi_\qq^{-1}e)$. By Definition \ref{def:nbr}, the $\qq$-neighbor of $K'_\qq$ along $H'$
is
\[
\Ok_\qq\pi_\qq^{-1}(\pi_\qq f)\oplus\Ok_\qq\pi_\qq(\pi_\qq^{-1}e)\perp W=\Ok_\qq f\oplus\Ok_\qq e\perp W=K_\qq .
\]
\end{proof}

The next two lemmas will produce such neighbors from a vector, under the hypothesis
$\ord_\qq\vv L\le1$. Recall that $\qq$ is odd, so that integral lattices over $\Ok_\qq$ are diagonalizable
(Lemma \ref{lem:split}).

\begin{lemma}\label{lem:qprim}
Let $L$ be an integral $\Ok_\qq$-lattice with $\ord_\qq\vv L\le1$. If
$x\in L\smallsetminus\qq L$ satisfies $\ord_\qq Q(x)\ge2$, then $B(x,L)=\Ok_\qq$.
\end{lemma}

\begin{proof}
Write $L=\bigperp_{i=1}^{n}\Ok_\qq u_i$ with $Q(u_i)=\alpha_i\in\Ok_\qq$, so that
$\sum_i\ord_\qq\alpha_i=\ord_\qq\vv L\le1$. After renumbering, $\alpha_i\in\Ok_\qq^{\times}$ for
$i\le n-1$ and $\ord_\qq\alpha_n\le1$. Write $x=\sum_ic_iu_i$. Then $B(x,u_i)=c_i\alpha_i$, so it
suffices to find $i\le n-1$ with $c_i\in\Ok_\qq^{\times}$. If there were none, then $c_n$ would be
a unit, since $x$ is primitive, and
\[
\ord_\qq Q(x)=\ord_\qq\Bigl(\sum_i c_i^{2}\alpha_i\Bigr)=\ord_\qq\alpha_n\le1,
\]
because $\ord_\qq(c_i^{2}\alpha_i)\ge2$ for $i\le n-1$ while $\ord_\qq(c_n^{2}\alpha_n)\le1$. This
contradicts the hypothesis.
\end{proof}

The hypothesis $\ord_\qq\vv L\le1$ of Lemma \ref{lem:qprim} cannot be dropped. For
$L=\langle1,1,1,\pi_\qq^{2}\rangle$ and $x$ the last basis vector, we have $B(x,L)=\qq^{2}$.

\begin{lemma}\label{lem:nbrvector}
Let $L$ be an integral $\Ok_\qq$-lattice on a quadratic space $U$ over $F_\qq$
with $\ord_\qq\vv L\le1$, and let $x\in L\smallsetminus\qq L$ satisfy $\ord_\qq Q(x)\ge2$. Put
\[
L_x:=\{w\in L\mid B(x,w)\in\qq\},\qquad L(\qq,x):=\Ok_\qq\,\pi_\qq^{-1}x+L_x .
\]
Then there is a hyperbolic plane $H$ in $L$ with $L(\qq,x)=L(\qq,H)$. In particular $L(\qq,x)$ is
integral and isometric to $L$.
\end{lemma}

\begin{proof}
Write $\pi:=\pi_\qq$. By Lemma \ref{lem:qprim} there is $z\in L$ with $B(x,z)=1$. Put
$\alpha:=Q(x)$ and $\beta:=Q(z)$, so $\ord_\qq\alpha\ge2$ and $\beta\in\Ok_\qq$. As $\qq$ is
odd, $1-\alpha\beta\in1+\qq^{2}\subset1+\qq$ is a square by Lemma \ref{lem:lst}. If $s$ is a
square root of $1-\alpha\beta$, then $s^{2}\equiv1\pmod{\qq^{2}}$, so one of $s-1$ and $s+1$ lies in
$\qq$, and not both, as their difference $2$ is a unit. Replacing $s$ by $-s$ if necessary, we
obtain $s\in\Ok_\qq^{\times}$ with $s^{2}=1-\alpha\beta$ and $s\equiv1\pmod\qq$, so that $s+1$ is a
unit. Put
\[
c:=\frac{-\alpha}{s+1},\qquad\text{so}\qquad \ord_\qq c=\ord_\qq\alpha\ge2 .
\]
From $(s-1)(s+1)=-\alpha\beta$ we get $c\beta=s-1$, whence
\[
Q(x+cz)=\alpha+2c+c^{2}\beta=\alpha+2c+c(s-1)=\alpha+c(s+1)=0 .
\]
Set $e:=x+cz$. Then $Q(e)=0$ and $e\equiv x\pmod{\qq^{2}L}$. Moreover $B(e,z)=1+c\beta=s$ is a
unit, so $B(e,L)=\Ok_\qq$.

Both $L_x$ and $L(\qq,x)$ are unchanged by this replacement. Indeed, for $w\in L$ we have
$B(e,w)-B(x,w)=cB(z,w)\in\qq^{2}$, so $B(e,w)\in\qq$ if and only if $B(x,w)\in\qq$. Thus
$L_e=L_x$. Moreover $\pi^{-1}e-\pi^{-1}x=(\pi^{-1}c)z$ with $\pi^{-1}c\in\qq$, and
$B\bigl(x,(\pi^{-1}c)z\bigr)=\pi^{-1}c\in\qq$, so $(\pi^{-1}c)z\in L_x$. Hence
$L(\qq,e)=L(\qq,x)$.

Choose $z'\in L$ with $B(e,z')=1$ and put $f:=z'-\tfrac12Q(z')e$, which is possible as $2$ is a
unit. Then $Q(f)=Q(z')-Q(z')B(z',e)+\tfrac14Q(z')^{2}Q(e)=0$ and $B(e,f)=B(e,z')=1$, so
$H:=\Ok_\qq e\oplus\Ok_\qq f$ is a hyperbolic plane. Being unimodular, it splits $L$, say
$L=H\perp W$.
For $w=ae+bf+w_0$ with $w_0\in W$ we have $B(e,w)=b$, so $L_e=\Ok_\qq e\oplus\qq f\perp W$ and
\[
L(\qq,e)=\Ok_\qq\pi^{-1}e+\Ok_\qq e\oplus\qq f\perp W
=\Ok_\qq\pi^{-1}e\oplus\Ok_\qq\pi f\perp W ,
\]
which is $L(\qq,H)$. As in the proof of Lemma \pref{lem:nbrbasic}{1}, we have
$L(\qq,H)\cong H\perp W=L$, and so $L(\qq,x)$ is integral.
\end{proof}

\begin{lemma}\label{lem:localchain}
Let $L,L'$ be integral $\Ok_\qq$-lattices on a quadratic space $U$ over
$F_\qq$ with $\vv L=\vv L'$ and $\ord_\qq\vv L\le1$. Then there are lattices
$L=L_0,L_1,\dots,L_s=L'$ on $U$ such that, for each $i$, $L_i$ is the $\qq$-neighbor of
$L_{i-1}$ along some hyperbolic plane of $L_{i-1}$.
\end{lemma}

\begin{proof}
We induct on the index $(L:L\cap L')$. Since $\vv L=\vv L'$ we have
$[L:L\cap L']^{2}\vv L=\vv(L\cap L')=[L':L\cap L']^{2}\vv L'$, so
$[L:L\cap L']=[L':L\cap L']$. In particular the index is $1$ exactly when $L=L'$, which is the
base case $s=0$.

Suppose $(L:L\cap L')>1$ and choose $y\in L'\smallsetminus L$. Let $j\ge1$ be minimal with
$\pi_\qq^{j}y\in L$ and put $x:=\pi_\qq^{j}y$. Then $x\notin\qq L$, by minimality of $j$, and
$Q(x)=\pi_\qq^{2j}Q(y)$ with $Q(y)\in\Ok_\qq$, so $\ord_\qq Q(x)\ge2j\ge2$. By
Lemma \ref{lem:nbrvector}, $L_1:=L(\qq,x)$ is the $\qq$-neighbor of $L$ along a hyperbolic plane.
That lemma also shows that $L_1$ is integral and isometric to $L$, so
$\vv L_1=\vv L=\vv L'$ and $\ord_\qq\vv L_1\le1$. Hence the induction hypothesis will apply to the
pair $(L_1,L')$.

It remains to show that $(L_1:L_1\cap L')<(L:L\cap L')$. First, we have $L\cap L'\subset L_x$. Indeed, for
$w\in L\cap L'$ we have $B(x,w)=\pi_\qq^{j}B(y,w)\in\qq^{j}\subset\qq$, since $y,w\in L'$ and
$L'$ is integral. Hence $L\cap L'\subset L_1$, and so $L\cap L'\subset L_1\cap L'$. Second,
this inclusion is strict, as $\pi_\qq^{-1}x=\pi_\qq^{\,j-1}y$ lies in $L'$ and in $L_1$, but not in
$L$, because $x\notin\qq L$. Therefore
\[
(L_1:L_1\cap L')=(L':L_1\cap L')<(L':L\cap L')=(L:L\cap L'),
\]
where the outer equalities hold because $\vv L_1=\vv L'=\vv L$. The induction hypothesis gives a
chain from $L_1$ to $L'$, and together with the step from $L$ to $L_1$ this is the required chain.
\end{proof}

\subsection{Kneser's construction}\label{subsec:kneser}

\begin{theorem}\label{thm:neighbour}
Let $K$ be an integral $\Ok$-lattice of rank $n\ge4$ on a quadratic space $V$ over $F$, and let
$\qq$ be an odd prime of $F$ with $\ord_{\qq}\vv K\le1$. Let $N$ be a lattice on $V$ isometric to
some member of $\spn^{+}(K)$. Then there are lattices $K=K_0,K_1,\dots,K_s$ on $V$ such that, for each
$i$, $K_i$ is the $\qq$-neighbor of $K_{i-1}$ along some hyperbolic plane of $(K_{i-1})_\qq$, and
$K_s$ is isometric to $N$.
\end{theorem}

\begin{proof}
Since $\qq$ is odd, $K_\qq$ is diagonalizable (Lemma \ref{lem:split}). As $\ord_\qq\vv K\le1$, at
most one scale exponent is nonzero, so $K_\qq=J_0\perp J_1$ with $J_0$ unimodular of rank
$\ge n-1\ge3$. By Lemma \pref{lem:odduni}{2}, $J_0$ splits a hyperbolic plane. In particular
$V_{\qq}$ is isotropic.

By Lemma \pref{lem:spinor}{1}, we may replace $N$ by an isometric lattice and assume that
$N\in\spn^{+}(K)$ and that $N_{\PP}=K_{\PP}$ for every prime $\PP\ne\qq$. As $N\in\gen(K)$, we
have $\vv N_\qq=\vv K_\qq$, and both lattices are integral. Lemma \ref{lem:localchain}, applied to
$L:=K_{\qq}$ and $L':=N_{\qq}$ on $U:=V_\qq$, yields lattices
$K_\qq=L_0',L_1',\dots,L_s'=N_\qq$ such that each $L_i'$ is the $\qq$-neighbor of $L_{i-1}'$ along
a hyperbolic plane. Let $K_i$ be the lattice on $V$ with
\[
(K_i)_{\qq}:=L_i',\qquad (K_i)_{\PP}:=K_{\PP}\ \ (\PP\ne\qq).
\]
Each $K_i$ is integral, and by Definition \pref{def:nbr}{2}, $K_i$ is the $\qq$-neighbor of
$K_{i-1}$ along the corresponding hyperbolic plane. Finally $K_s=N$.
\end{proof}

\subsection{The bound for \texorpdfstring{$r$}{r}}

\begin{lemma}\label{lem:rbound}
Let $K$ be an integral $\Ok$-lattice of rank $\ge4$ on a quadratic space $V$ over $F$, and let
$\qq$ be an odd prime of $F$ with $\ord_{\qq}\vv K\le1$. Then every lattice in $\spn^{+}(K)$ is
isometric to a lattice $N'$ with $\qq^{h(\gen K)-1}N'\subset K$, and
$r_\qq(K)\le h(\gen K)-1$.
\end{lemma}

\begin{proof}
Every lattice $K''\in\gen(K)$ is integral with $\vv K''=\vv K$, so $\ord_\qq\vv K''\le1$. As in
the proof of Theorem \ref{thm:neighbour}, $K''_\qq$ splits a hyperbolic plane. Thus $K''$ has
$\qq$-neighbors in the sense of Definition \ref{def:nbr}, and they lie in $\gen(K)$ by
Lemma \pref{lem:nbrbasic}{1}.

Consider the graph $\mathcal G$ whose vertices are the classes in $\gen(K)$, where two classes are
joined if some lattice in one of them is a $\qq$-neighbor of some lattice in the other. Let
$c_1$ and $c_2$ be joined classes and let $K_1\in c_1$. Then some lattice in $c_2$ is a
$\qq$-neighbor of $K_1$. Indeed, by Lemma \pref{lem:nbrbasic}{4} the relation of being a
$\qq$-neighbor is symmetric, so there are $L_1\in c_1$ and $L_2\in c_2$ such that $L_2$ is a
$\qq$-neighbor of $L_1$. Choose $\sigma\in O(V)$ with $\sigma L_1=K_1$. Then $\sigma L_2$ lies in
$c_2$, and it is a $\qq$-neighbor of $K_1$ by Lemma \pref{lem:nbrbasic}{3}. Hence, for
every path of length $s$ in $\mathcal G$ from the class of $K$ to a class $c$, there are lattices
$K=K_0,K_1,\dots,K_s$ such that $K_i$ is a $\qq$-neighbor of $K_{i-1}$ for each $i$, and
$K_s\in c$.

Now let $N\in\spn^{+}(K)$. By Theorem \ref{thm:neighbour}, the class of $N$ lies in the connected
component of the class of $K$. This component has at most $h(\gen K)$ vertices, so there is a path
of length $s\le h(\gen K)-1$ from the class of $K$ to the class of $N$. Let $K=K_0,K_1,\dots,K_s$
be the corresponding lattices. We have $\qq K_i\subset K_{i-1}$ for each $i$ by
Lemma \pref{lem:nbrbasic}{2}, and therefore $\qq^{s}K_s\subset K$. So the lattice $N':=K_s$ is
isometric to $N$ and satisfies $\qq^{h(\gen K)-1}N'\subset\qq^{s}N'\subset K$.

Finally, $r_\qq(K)$ is defined, because $K_\qq$ is isotropic. By the first assertion, the hypothesis of Lemma \ref{lem:spinorcrit} holds with
$s:=h(\gen K)-1$. So for every $r\ge h(\gen K)-1$, every totally positive $\beta\in\Ok$ that is
represented by $\gen(\qq^{r}K)$ is represented by $K$. By Definition \ref{def:rK}, this means that
$r_\qq(K)\le h(\gen K)-1$.
\end{proof}

Note that the graph $\mathcal G$ is formed from classes, and not from proper classes, although
Theorem \ref{thm:neighbour} concerns the proper spinor genus. This suffices, because
Lemma \ref{lem:spinorcrit} only requires that the lattice $N'$ be isometric to $N$.

A bound of this kind is asserted, without proof, in the remark following
\cite[Theorem 3.1]{HI}, there with the number of classes in the proper spinor genus in place of
$h(\gen K)$, and it is used in the form above in the example at the end of \cite[\S4]{HI}.

\begin{proposition}\label{prop:rH}
Let $r$, $\mathcal K$, and $X$ be as in \S\textup{\ref{subsec:setup}}. Then
$r\le\mathcal H_{m-1}(X)-1$.
\end{proposition}

\begin{proof}
As $\mathcal H_{m-1}(X)\ge1$, it suffices to show that $r_\qq(K)\le\mathcal H_{m-1}(X)-1$ for
every lattice $K$ that represents a class in $\mathcal K$. Every lattice in $\gen(K)$ is integral
and positive definite of rank $m-1$, with volume ideal $\vv K$, so its class lies in
$\mathcal K$. Hence $h(\gen K)\le\#\mathcal K\le\mathcal H_{m-1}(X)$, and Lemma \ref{lem:rbound}
gives $r_\qq(K)\le h(\gen K)-1\le\mathcal H_{m-1}(X)-1$.
\end{proof}

\subsection{Explicit bounds}\label{subsec:explicit}

We can now prove Theorem \ref{cor:crude}. We will first make the bound of
Corollary \ref{cor:exponent} explicit for a suitable prime $\qq$ (Theorem \ref{thm:explicitintro}),
and then simplify the resulting constant.

\begin{lemma}\label{lem:leastq}
Let $a$ be a positive integer, and let $p_0$ be the least prime not dividing $a$. Then
$p_0\le\max\{41,1+\tfrac32\log a\}$.
\end{lemma}

\begin{proof}
Suppose $p_0>41$ and put $x:=p_0-1\ge42$. Every prime $p\le x$ divides $a$, so
$\vartheta(x):=\sum_{p\le x}\log p\le\log a$. By \cite[(3.16)]{RS},
$\vartheta(x)>x(1-1/\log x)\ge x(1-1/\log42)>\tfrac23x$. Hence $p_0-1<\tfrac32\log a$.
\end{proof}

Put $n:=m-1$, let $V_{m,F}$ be as in \eqref{eq:V}, and define
\begin{gather}
C_{m,F}:=2\bigl(d_F^{n/2}n^{ln}\bigr)^{ln+1}(6l)^{lnm/2}
\bigl(n^{ln}d_F^{n}V_{m,F}\bigr)^{m/2},\notag\\
J_{m,F}:=33+n\log_2n+\frac nl\log_2d_F+\frac1l\log_2V_{m,F},\notag\\
\Xi_{m,F}:=12l\cdot30^{2lm}\,C_{m,F}\,J_{m,F}^{\,m-2}+\log G_{m,F}+4,\label{eq:Xi}
\end{gather}
with $G_{m,F}$ as in \eqref{eq:G}.

\begin{theorem}\label{thm:explicitintro}
Let $M$ be an integral positive definite $\Ok$-lattice of rank $m\ge5$, put
$D:=\Nm(\vv M)$, and let $A\in\Ok$ be totally positive and represented by $M_\PP$ for every
prime $\PP$. If
\[
\NFQ(A)\ >\ \exp\Bigl(\Xi_{m,F}\,D^{\rho_m}\bigl(\log(4D)\bigr)^{2lm+m-1}\Bigr),
\]
then $A$ is represented by $M$.
\end{theorem}

\begin{proof}
Let $q$ be the least rational prime not dividing $2D$, and let $\qq$ be a prime of $F$ above $q$.
Then $\qq$ is odd, $\qq\nmid\vv M$, and $\Nm\qq\le q^{l}$. Put $\Lambda:=\log(4D)\ge\log4$. By
Lemma \ref{lem:leastq}, $q\le\max\{41,1+\tfrac32\log(2D)\}\le30\Lambda$, as $30\log4>41$.

Let $r$ and $X$ be as in \S\ref{subsec:setup} for this prime $\qq$, and write $\mathcal H$ for
$\mathcal H_{m-1}(X)$. We have $1+r\le\mathcal H$ by Proposition \ref{prop:rH}. Hence, by
Corollary \ref{cor:exponent}, it suffices to show that
\begin{equation}\label{eq:logbound}
\log\Bigl(G_{m,F}\,q^{4l\mathcal H}D^{1+\frac2{m-3}+\frac2{m-4}}\Bigr)\le
\Xi_{m,F}\,D^{\rho_m}\Lambda^{2lm+m-1}.
\end{equation}

Let us bound $\mathcal H$. Let $d_m$ be as in \S\ref{subsec:setup}, so that
$m(1+1/d_m)=\rho_m$. As $\Nm\qq\le q^{l}$, $d_m\le m-3$, and $\Dp\Dodd=D$, the definition of $X$ gives
$X\le V_{m,F}\,q^{4l}D^{2+2/d_m}$. Put
\[
U:=U_n\bigl(V_{m,F}\,q^{4l}D^{2+2/d_m}\bigr)=n^{n}d_F^{n/l}V_{m,F}^{1/l}\,q^{4}D^{(2+2/d_m)/l},
\]
so that $U^{lm/2}=\bigl(n^{ln}d_F^{n}V_{m,F}\bigr)^{m/2}q^{2lm}D^{\rho_m}$.
The function $\mathcal H_n$ of Proposition \ref{prop:classes} is increasing, and $l(n+1)/2=lm/2$.
So the definitions of $\mathcal H_n$ and of $C_{m,F}$ give
\[
\mathcal H\le C_{m,F}\,q^{2lm}D^{\rho_m}\bigl(1+\log_2U\bigr)^{m-2}.
\]
We now bound the logarithmic factor. We have
\[
\log_2U=n\log_2n+\frac nl\log_2d_F+\frac1l\log_2V_{m,F}+4\log_2q+\frac{2+2/d_m}{l}\log_2D .
\]
As $q\le30\Lambda$ and $\log\Lambda\le\Lambda$, we have $\log_2q\le\log_230+\Lambda/\log2$. Moreover
$\log_2D\le\Lambda/\log2$ and $2+2/d_m\le4$. Since $\Lambda\ge1$ and $1+4\log_230+8/\log2<33$, the
definition of $J_{m,F}$ therefore gives $1+\log_2U\le J_{m,F}\Lambda$.

We also have $\log q\le\log30+\log\Lambda\le\log30+\Lambda-1\le3\Lambda$. Together with $q\le30\Lambda$
and the bounds for $\mathcal H$ and for $1+\log_2U$, this gives
\[
4l\,\mathcal H\log q\le12l\cdot30^{2lm}C_{m,F}J_{m,F}^{\,m-2}D^{\rho_m}\Lambda^{2lm+m-1}.
\]
Finally, $1+\frac2{m-3}+\frac2{m-4}\le4$ gives
$\log\bigl(G_{m,F}D^{1+\frac2{m-3}+\frac2{m-4}}\bigr)\le\log G_{m,F}+4\Lambda$, and
$\log G_{m,F}+4\Lambda\le(\log G_{m,F}+4)D^{\rho_m}\Lambda^{2lm+m-1}$ because
$D^{\rho_m}\Lambda^{2lm+m-1}\ge\Lambda\ge1$. Adding the last two estimates and recalling the
definition \eqref{eq:Xi} of $\Xi_{m,F}$, we obtain
\eqref{eq:logbound}.
\end{proof}

\begin{proof}[Proof of Theorem \textup{\ref{cor:crude}}]
Put $Z:=2lm\,d_F\ge10$. By Theorem \ref{thm:explicitintro}, it suffices to show that
$\Xi_{m,F}\le Z^{\,l^{2}m^{2}+4lm^{2}}$. Put
\[
\Xi_0:=12l\cdot30^{2lm}C_{m,F}J_{m,F}^{\,m-2},\qquad\text{so that}\quad\Xi_{m,F}=\Xi_0+\log G_{m,F}+4 .
\]

We first show that $\Xi_{m,F}\le2\Xi_0$. As $C_{m,F}\ge V_{m,F}^{m/2}$ and $J_{m,F}\ge1$, we have
$\Xi_0\ge30^{2lm}V_{m,F}^{m/2}\ge30^{2lm}+V_{m,F}^{m/2}$, where the second inequality holds because
both factors are at least $2$. Now
\[
30^{2lm}\ge l\bigl(12\log2+5\log l+(m+2)\log m\bigr)+4\qquad\text{and}\qquad
V_{m,F}^{m/2}\ge d_F^{(m+l)m/2}\ge(m+2l)\log d_F ,
\]
and the sum of the two right-hand sides is $\log G_{m,F}+4$ by \eqref{eq:G}. So
$\Xi_0\ge\log G_{m,F}+4$, and $\Xi_{m,F}\le2\Xi_0$.

Next we bound the factors of $2\Xi_0$ by powers of $Z$. We have
\[
24l\le Z^{2},\qquad30\le Z^{2},\qquad6l\le Z,\qquad n^{l}d_F\le Z^{l},\qquad
V_{m,F}\le Z^{lm+m+9l},
\]
and the fourth inequality gives $d_F^{n/2}n^{ln}\le Z^{ln}$ and $n^{ln}d_F^{n}\le Z^{ln}$. We also
have $J_{m,F}\le Z^{2}$. Indeed, if $l\ge2$, then $d_F\ge2$, so $Z\ge8m\ge40$ and
\[
J_{m,F}\le33+(2n+2m+9)\log_2Z\le(4m+17)\log_2Z\le\bigl(\tfrac Z2+17\bigr)\log_2Z\le Z^{2}.
\]
If $l=1$, then $F=\Q$, $Z=2m$, and $J_{m,\Q}=43+n\log_2n+(m+2)\log_2m\le4m^{2}$. Inserting these
bounds into the definitions of $C_{m,F}$ and of $\Xi_0$, and using $2\le Z$, we obtain
\[
2\Xi_0\le Z^{E}\qquad\text{with}\quad E:=3+4lm+ln(ln+1)+lnm+\tfrac12(lm+m+9l)m+2(m-2).
\]

It remains to bound $E$. With $n=m-1$ we get
\[
E=l^{2}(m-1)^{2}+\tfrac32lm^{2}+\tfrac12m^{2}+\tfrac{17}2lm+2m-l-1 .
\]
As $m\ge5$, we have $l^{2}(m-1)^{2}=l^{2}m^{2}-l^{2}(2m-1)\le l^{2}m^{2}-\tfrac95lm$, and therefore
\[
E\le l^{2}m^{2}+\tfrac32lm^{2}+\tfrac12m^{2}+\tfrac{67}{10}lm+2m .
\]
Using $l\ge1$ and $m\ge5$ to bound each of the last three terms by a multiple of $lm^{2}$, we
conclude that
\[
E\le l^{2}m^{2}+lm^{2}\bigl(\tfrac32+\tfrac12+\tfrac{67}{50}+\tfrac25\bigr)
\le l^{2}m^{2}+4lm^{2}.
\]
\end{proof}

\section{Real quadratic fields}\label{sec:rq}

In this section we will collect the tools for the proof of Theorem \ref{thm:tuniv}, namely
elements of small norm in prescribed residue classes (\S\ref{subsec:GoN}), local criteria for the
representation of elements of $\Ok_\PP$ (\S\ref{subsec:locallemmas}), and the sublattice of rank five
provided by \cite{KK+} (\S\ref{subsec:rankfive}).

Throughout \S\S\ref{sec:rq}--\ref{sec:pell}, $d>1$ is a squarefree integer and $F=\Q(\sqrt d)$, so
that $l=2$. Integer-valued lattices, integral lattices, and $S$-universality are as defined in
\S\ref{subsec:critintro}.
We number the real embeddings so that $\tau_1(\sqrt d)>0$, identify $F$ with $\tau_1(F)\subset\R$,
and write $\alpha':=\tau_2(\alpha)$ for $\alpha\in F$. Thus $\alpha$ is totally positive if and
only if $\alpha>0$ and $\alpha'>0$, and $\NFQ(\alpha)=\alpha\alpha'$. Put
$\omega:=\frac12(1+\sqrt d)$ if $d\equiv1\pmod 4$, and $\omega:=\sqrt d$ otherwise. Then
$\Ok=\Z+\Z\omega$ and $\omega-\omega'=\sqrt{d_F}$. Put
\begin{equation}\label{eq:Theta}
\Theta_1:=\omega-\lfloor\omega'\rfloor,\qquad\Theta_2:=2\omega-\lfloor2\omega'\rfloor .
\end{equation}
As $\omega'\notin\Q$, the conjugates $\Theta_1'=\omega'-\lfloor\omega'\rfloor$ and
$\Theta_2'=2\omega'-\lfloor2\omega'\rfloor$ lie in $(0,1)$, and $\Theta_i=i\sqrt{d_F}+\Theta_i'$
for $i=1,2$. Hence $\Theta_1,\Theta_2\in\Ok^{+}$ and
\begin{equation}\label{eq:NTheta}
\Theta_2<2\sqrt{d_F}+1,\qquad\NFQ(\Theta_1)<\sqrt{d_F}+1,\qquad\NFQ(\Theta_2)<2\sqrt{d_F}+1 .
\end{equation}

For an odd prime $p$ let $\gamma_p$ be the least positive integer that is a quadratic nonresidue
modulo $p$, and put $\gamma_2:=5$. For $\mu\in\Z_{>0}$ put
\begin{equation}\label{eq:gmu}
g(\mu):=\max\{\gamma_p\mid p\text{ prime},\ p\le2\mu^{2}\},\qquad
\mathcal E_\mu:=\{1,2,\dots,g(\mu)\}\cup\{\Theta_1,\Theta_2\}\subset\Ok^{+}.
\end{equation}
As $\gamma_3=\gamma_5=2$ and $\gamma_7=3$, we have $g(1)=g(2)=5$.

\begin{lemma}\label{lem:gmu}
For every $\mu\in\Z_{>0}$ we have $5\le g(\mu)\le\max\{5,\sqrt2\,\mu\}$.
\end{lemma}

\begin{proof}
The lower bound holds because $\gamma_2=5$ and $2\le2\mu^{2}$. For the upper bound, recall that
$\gamma_3=2$, $\gamma_7=3$, and $\gamma_{23}=5$, and that $\gamma_p<\sqrt p$ for every other odd prime
$p$ \cite[p. 123, remark after Theorem 6-35]{Le}. Hence
$\gamma_p\le\max\{5,\sqrt p\}\le\max\{5,\sqrt2\,\mu\}$ for every odd prime $p\le2\mu^{2}$.
\end{proof}

\subsection{Elements of small norm in a residue class}\label{subsec:GoN}

\begin{lemma}\label{lem:GoN}
Let $\mmm$ be a nonzero ideal of $\Ok$, let $k_\mmm$ be the least positive integer in $\mmm$, and
let $a\in\Ok$. Then $k_\mmm\le\Nm\mmm$, and there is a totally positive $\beta\in a+\mmm$ with
\[
\NFQ(\beta)\le k_\mmm^{2}+\tfrac12\Nm(\mmm)\sqrt{d_F}.
\]
\end{lemma}

\begin{proof}
The integer $\Nm\mmm=\#(\Ok/\mmm)$ annihilates $\Ok/\mmm$, so $\Nm\mmm\in\mmm$ and
$k_\mmm\le\Nm\mmm$. By the Hermite normal form, the $\Z$-module $\mmm\subset\Z+\Z\omega$ has a basis
$k,\ b+c\omega$ with $k,c\in\Z_{>0}$ and $b\in\Z$. Then $\mmm\cap\Z=k\Z$, so that $k=k_\mmm$, and
$\Nm\mmm=kc$. Subtracting from $a$ a suitable integer multiple of $b+c\omega$, we may assume that
$a=x+y\omega$ with $x,y\in\Z$ and $|y|\le c/2$. Then
$|a-a'|=|y|(\omega-\omega')\le\frac c2\sqrt{d_F}$.

Let $j\in\Z$ be minimal such that $\beta:=a+jk$ is totally positive. It exists because $\beta$ and
$\beta'$ tend to $+\infty$ as $j\to+\infty$ and to $-\infty$ as $j\to-\infty$. By minimality,
$\min\{\beta,\beta'\}\le k$, since otherwise $a+(j-1)k$ would be totally positive. As
$\beta-\beta'=a-a'$, it follows that
$\max\{\beta,\beta'\}=\min\{\beta,\beta'\}+|a-a'|\le k+\frac c2\sqrt{d_F}$, and
\[
\NFQ(\beta)=\min\{\beta,\beta'\}\cdot\max\{\beta,\beta'\}\le k\Bigl(k+\frac c2\sqrt{d_F}\Bigr)
=k_\mmm^{2}+\tfrac12\Nm(\mmm)\sqrt{d_F}.
\qedhere
\]
\end{proof}

\subsection{Local lemmas}\label{subsec:locallemmas}

For a prime $\PP$, every class of $F_\PP^{\times}/(F_\PP^{\times})^{2}$ contains an element
$\alpha\in\Ok_\PP$ with $\ord_\PP\alpha\in\{0,1\}$, as $\pi_\PP^{2}$ is a square.

\begin{lemma}\label{lem:squareclasses}
Let $\PP$ be a prime, let $N$ be an $\Ok_\PP$-lattice, and let $t,\alpha\in\Ok_\PP$ be nonzero with
$\ord_\PP\alpha\in\{0,1\}$. Assume that $N$ represents an element of
$t\alpha(\Ok_\PP^{\times})^{2}$. Then $N$ represents every nonzero $\beta\in t\Ok_\PP$ for which
$\beta/t$ lies in the class of $\alpha$ in $F_\PP^{\times}/(F_\PP^{\times})^{2}$.
\end{lemma}

\begin{proof}
Let $\beta$ be such an element, and write $\beta/t=\alpha\lambda^{2}$ with
$\lambda\in F_\PP^{\times}$. As $\ord_\PP(\beta/t)\ge0$ and $\ord_\PP\alpha\le1$, we have
$2\ord_\PP\lambda\ge-1$, and so $\lambda\in\Ok_\PP$. Choose $x\in N$ and $z\in\Ok_\PP^{\times}$ with
$Q(x)=t\alpha z^{2}$. Then $\lambda z^{-1}x\in N$ and $Q(\lambda z^{-1}x)=t\alpha\lambda^{2}=\beta$.
\end{proof}

\begin{lemma}\label{lem:goodprimes}
Let $\PP$ be odd and let $N$ be an $\Ok_\PP$-lattice of rank $5$ with $\ssc N\subset\Ok_\PP$ and
$\ord_\PP\vv N\le2$. Then $Q(N)=\Ok_\PP$.
\end{lemma}

\begin{proof}
By Lemma \ref{lem:split}, $N=\langle\alpha_1,\dots,\alpha_5\rangle$ with
$\ord_\PP\alpha_1\le\dots\le\ord_\PP\alpha_5$. These orders are the scale exponents of
Definition \ref{def:exponents}, so they are nonnegative with sum $\ord_\PP\vv N\le2$. Hence
$\alpha_1,\alpha_2,\alpha_3$ are units, $\langle\alpha_1,\alpha_2,\alpha_3\rangle$ is unimodular of
rank $3$, and $Q(N)\supset\Ok_\PP$ by Lemma \pref{lem:odduni}{2}. Finally $Q(N)\subset\Ok_\PP$,
as $\ssc N\subset\Ok_\PP$.
\end{proof}

The next lemma produces elements of $\Ok$ in prescribed square classes at finitely many primes.

\begin{lemma}\label{lem:pin}
Let $\mathcal S$ be a finite set of primes, and for $\PP\in\mathcal S$ let $\alpha_\PP\in\Ok_\PP$
with $j_\PP:=\ord_\PP\alpha_\PP\in\{0,1\}$. Put
$\mmm:=\prod_{\PP\in\mathcal S}\PP^{\,j_\PP+2e_\PP+1}$. Then there is $a\in\Ok$ such that
$\beta\in\alpha_\PP(\Ok_\PP^{\times})^{2}$ for every $\beta\in a+\mmm$ and every $\PP\in\mathcal S$.
\end{lemma}

\begin{proof}
By the Chinese remainder theorem there is $a\in\Ok$ with
$a-\alpha_\PP\in\PP^{\,j_\PP+2e_\PP+1}\Ok_\PP$ for every $\PP\in\mathcal S$. If $\beta\in a+\mmm$ and
$\PP\in\mathcal S$, then $\beta-\alpha_\PP\in\PP^{\,j_\PP+2e_\PP+1}\Ok_\PP$, so
$\beta/\alpha_\PP\in1+\PP^{2e_\PP+1}\Ok_\PP\subset(\Ok_\PP^{\times})^{2}$ by Lemma \ref{lem:lst}.
\end{proof}

At the dyadic primes we will use the following bounds. Recall that
$\sum_{\pp\mid2}e_\pp f_\pp=2$, as $l=2$.

\begin{lemma}\label{lem:dyadic}
For every prime $\pp\mid2$ let $j_\pp\in\{0,1\}$, and put
$\mathfrak d:=\prod_{\pp\mid2}\pp^{\,j_\pp+2e_\pp+1}$. Then $\Nm(\mathfrak d)\le256$ and
$16\in\mathfrak d$.
\end{lemma}

\begin{proof}
As $\Nm\pp=2^{f_\pp}$, we have $\Nm(\mathfrak d)=2^{s}$ with
\[
s=\sum_{\pp\mid2}f_\pp(j_\pp+2e_\pp+1)\le2\sum_{\pp\mid2}e_\pp f_\pp+2\sum_{\pp\mid2}f_\pp\le8 ,
\]
because $\sum_{\pp\mid2}f_\pp\le\sum_{\pp\mid2}e_\pp f_\pp=2$.
For every $\pp\mid2$ we have $\ord_\pp16=4e_\pp\ge2e_\pp+2\ge j_\pp+2e_\pp+1$, and so
$16\in\mathfrak d$.
\end{proof}

\subsection{Volumes and a sublattice of rank five}\label{subsec:rankfive}

\begin{lemma}\label{lem:hadamard}
Let $L$ be an integer-valued lattice generated by nonzero vectors $w_1,\dots,w_n$. Then
$\Nm(\vv L)\le\prod_{i=1}^{n}\NFQ\bigl(Q(w_i)\bigr)$.
\end{lemma}

\begin{proof}
Choose $I\subset\{1,\dots,n\}$ such that $(w_i)_{i\in I}$ is a basis of $FL$, and put
$L_0:=\sum_{i\in I}\Ok w_i$. Then $\vv L_0$ is generated by
$\delta:=\det\bigl(B(w_i,w_{i'})\bigr)_{i,i'\in I}$. For $k=1,2$, the real matrix
$\bigl(\tau_kB(w_i,w_{i'})\bigr)_{i,i'\in I}$ is positive definite with diagonal entries
$\tau_kQ(w_i)$. So Hadamard's inequality \cite[Theorem 7.8.1]{HJ} gives
$0<\tau_k\delta\le\prod_{i\in I}\tau_kQ(w_i)$, and multiplying over $k$ we get
$\Nm(\vv L_0)\le\prod_{i\in I}\NFQ(Q(w_i))$. As $L_0\subset L$ and $\rank L_0=\rank L$, we have
$\vv L_0=[L:L_0]^{2}\vv L$ (\S\ref{sec:notation}), so $\Nm(\vv L)\le\Nm(\vv L_0)$. Finally
$\NFQ(Q(w_i))\ge1$ for every $i$, as $Q(w_i)$ is a totally positive element of $\Ok$.
\end{proof}

The lattices of \cite{KK+} are our integer-valued lattices, and its classical lattices are our
integral lattices, with
$Q(x)=B(x,x)$. In the notation of \cite{KK+}, $\Delta_D=d_F$ and $\omega_D=\omega$.

\begin{proposition}\label{prop:five}
Let $\mu\in\Z_{>0}$, put $g:=g(\mu)$, assume $\sqrt{d_F}\ge8g\mu^{2}+\frac12$, and let $L$ be a
$\mu\mathcal E_\mu$-universal integral lattice. Then there are an integer $C_1$ and linearly
independent $v_1,\dots,v_5\in L$ with
\[
Q(v_1)=\mu,\quad Q(v_2)=2\mu,\quad Q(v_3)=C_1\mu,\quad Q(v_4)=\mu\Theta_1,\quad
Q(v_5)=\mu\Theta_2 ,
\]
where $2\le C_1\le g$.
\end{proposition}

\begin{proof}
We apply the results of \cite{KK+} with $m=\mu$. Choose $v_1,v_2\in L$ with $Q(v_1)=\mu$ and
$Q(v_2)=2\mu$. They are linearly independent. Indeed, $v_2=\lambda v_1$ with $\lambda\in F$ would
give $\lambda^{2}=2$, and then $F=\Q(\sqrt2)$ and $d_F=8$, whereas $d_F\ge(8g\mu^2+\frac12)^2>1640$. As
$Q(v_1)Q(v_2)=2\mu^{2}\prec d_F/4$, \cite[Lemma 4.1(2)]{KK+} gives $B(v_1,v_2)\in\Z$. Hence
\cite[Corollary 4.11(1)]{KK+} applies to $v_1,v_2$. It yields a prime $p_1\le2\mu^{2}$ and an
integer $C_1$ with $2\le C_1\le\gamma_{p_1}\le g$, where for $p_1=2$ the value $\gamma_2=5$ is
admissible by \cite[Remark 4.7]{KK+}.

Now choose $v_3,v_4,v_5\in L$ with $Q(v_3)=C_1\mu$, $Q(v_4)=\mu\Theta_1$, and $Q(v_5)=\mu\Theta_2$.
They exist because these elements lie in $\mu\mathcal E_\mu$. The proof of
\cite[Theorem 5.4]{KK+} shows that $v_1,\dots,v_5$ are linearly independent. Indeed, that proof
uses its hypothesis that $L$ is $\mu\Ok^{+}$-universal only to obtain vectors $v_1,\dots,v_5$ with
the five values above, where $C_1$ is the integer that \cite[Corollary 4.11(1)]{KK+} attaches to
$v_1,v_2$. Its remaining hypothesis is $d_F>4C_1\mu^{2}\Theta_2$. This holds, since
\eqref{eq:NTheta} and $C_1\le g$ give
\[
4C_1\mu^{2}\Theta_2<4g\mu^{2}\bigl(2\sqrt{d_F}+1\bigr)<d_F ,
\]
where the last inequality holds because
$d_F-8g\mu^{2}\sqrt{d_F}=\sqrt{d_F}\bigl(\sqrt{d_F}-8g\mu^{2}\bigr)\ge\frac12\bigl(8g\mu^{2}+\frac12\bigr)
>4g\mu^{2}$.
\end{proof}

For $\mu=1$ and $\mu=2$, the hypothesis of Proposition \ref{prop:five} reads $d_F\ge1641$ and
$d_F\ge25\,761$, respectively, as $d_F$ is an integer. The same argument shows that the bounds
$2025$ and $32\,400$ in the row ``$\rank L\le4$'' of \cite[Theorem 1.2]{KK+} may be replaced by
$1641$ and $25\,761$. The bounds of \cite{KK+} come from the cruder estimate
$\sqrt{d_F}<9C_1\mu^{2}$ in
\cite[Corollary 3.3(2)]{KK+}.

\section{Criterion sets for real quadratic fields}\label{sec:crit}

In this section we will prove Theorem \ref{cor:critexplicit}. We will deduce it from
Theorem \ref{thm:tuniv}, which gives criterion sets in terms of the following numbers.

\begin{definition}\label{def:T}
For an integer $m_0\ge5$ and a real number $W\ge1$, let $\mathcal T_{F,m_0}(W)$ be the infimum of the
real numbers $c\ge0$ with the following property: if $L$ is an integral lattice with
$5\le\rank L\le m_0$ and $\Nm(\vv L)\le W$, and $A\in\Ok^{+}$ satisfies $\NFQ(A)>c$ and
$A\in Q(L_\PP)$ for every prime $\PP$, then $A\in Q(L)$.
\end{definition}

By Theorem \ref{cor:crude}, applied with $5\le m\le m_0$ and $D\le W$, the number
$\mathcal T_{F,m_0}(W)$ is finite. It is non-decreasing in $m_0$ and in $W$. Moreover,
$c=\mathcal T_{F,m_0}(W)$ itself has the property in Definition \ref{def:T}, because for every $A$
with $\NFQ(A)>\mathcal T_{F,m_0}(W)$ some $c<\NFQ(A)$ has this property. We state the results of
this section in terms of $\mathcal T_{F,m_0}(W)$, so that every upper bound for these numbers can
be inserted into them.

Let $t\in\Z_{>0}$. A set $\mathcal C\subset t\Ok^{+}$ is a \emph{criterion set for
$t\Ok^{+}$-universality of integer-valued lattices} if every $\mathcal C$-universal integer-valued lattice is
$t\Ok^{+}$-universal. It is a
\emph{criterion set for $t\Ok^{+}$-universality of integral lattices} if every
$\mathcal C$-universal integral lattice is $t\Ok^{+}$-universal. Up to passing to classes modulo
$\mathcal U_F^{2}$, these are the $(S,\mathrm{nc})$- and $(S,\mathrm{cl})$-criterion sets of
\cite[Definition 2.1]{KKR} for $S=t\Ok^{+}$. For $t=1$ they are the universality criterion sets of
Definition \ref{def:crit}. Recall the numbers $g(\mu)$ from \eqref{eq:gmu}.

\begin{theorem}\label{thm:tuniv}
Let $F$ be a real quadratic field and $t\in\Z_{>0}$.
\begin{enumerate}
\item\label{thm:tuniv:1} Put $g:=g(t)$, assume $\sqrt{d_F}\ge8gt^{2}+\frac12$, and put
$W:=258\,t^{2}\bigl(9g^{2}t^{10}d_F\bigr)^{7/3}$. Then
\[
\bigl\{A\in t\Ok^{+}\mid\NFQ(A)\le\max\{\mathcal T_{F,6}(W),W\}\bigr\}
\]
is a criterion set for $t\Ok^{+}$-universality of integral lattices.
\item\label{thm:tuniv:2} Put $g:=g(2t)$, assume $\sqrt{d_F}\ge32gt^{2}+\frac12$, and put
$W:=1.1\cdot10^{10}\,t^{2}\bigl(9g^{2}t^{10}d_F\bigr)^{7/3}$. Then
\[
\bigl\{A\in t\Ok^{+}\mid\NFQ(A)\le\max\{\tfrac14\mathcal T_{F,6}(W),W\}\bigr\}
\]
is a criterion set for $t\Ok^{+}$-universality of integer-valued lattices.
\end{enumerate}
\end{theorem}

By Lemma \ref{lem:gmu}, $g(t)\le\max\{5,\sqrt2\,t\}$ in \hyperref[thm:tuniv:1]{(1)} and
$g(2t)\le\max\{5,2\sqrt2\,t\}$ in \hyperref[thm:tuniv:2]{(2)}. For $t=1$ we have $g(1)=g(2)=5$, and the hypotheses of Theorem \ref{thm:tuniv} read
$\sqrt{d_F}\ge40.5$ and $\sqrt{d_F}\ge160.5$, i.e., $d_F\ge1641$ and $d_F\ge25\,761$, as $d_F$ is an
integer. Moreover $258\cdot225^{7/3}<8\cdot10^{7}$ and $1.1\cdot10^{10}\cdot225^{7/3}<3.4\cdot10^{15}$.
The number $\mathcal T_{F,6}(W)$ is non-decreasing in $W$, and a set that contains a universality
criterion set of one of the two kinds is again one of this kind. Thus Theorem \ref{thm:tuniv} has the following
special case.

\begin{corollary}\label{cor:tone}
Let $F$ be a real quadratic field.
\begin{enumerate}
\item\label{cor:tone:1} If $d_F\ge1641$ and $W:=8\cdot10^{7}\,d_F^{7/3}$, then
$\bigl\{A\in\Ok^{+}\mid\NFQ(A)\le\max\{\mathcal T_{F,6}(W),W\}\bigr\}$ is a universality
criterion set for integral lattices.
\item\label{cor:tone:2} If $d_F\ge25\,761$ and $W:=3.4\cdot10^{15}\,d_F^{7/3}$, then
$\bigl\{A\in\Ok^{+}\mid\NFQ(A)\le\max\{\tfrac14\mathcal T_{F,6}(W),W\}\bigr\}$ is a universality
criterion set for integer-valued lattices.
\end{enumerate}
\end{corollary}

\subsection{Proof of Theorem \ref{thm:tuniv}}

We will treat both cases at once. In case (1) put $\mu:=t$, and in case (2) put $\mu:=2t$. In both
cases $g=g(\mu)\ge5$ and $\sqrt{d_F}\ge8g\mu^{2}+\frac12\ge40.5$. Put $V_0:=9g^{2}t^{10}d_F$, put
$V:=V_0$ in case (1) and $V:=2^{10}V_0$ in case (2), and
\begin{equation}\label{eq:Pi}
\Pi:=258\,t^{2}\,V_0\,V^{4/3}.
\end{equation}
The number $V_0$ will bound the norm of the volume of a sublattice $L_5$ of rank five, and $V$ the
norm of the odd part of that volume. The number $\Pi$ will bound the norms of all auxiliary
elements that we need, and the norms of the volumes of the lattices to which Tartakowsky's theorem
will be applied. In case (1) we have $\Pi=258\,t^{2}V_0^{7/3}=W$. In case (2) we have
$4^{6}\,\Pi=4^{6}\cdot2^{40/3}\cdot258\,t^{2}V_0^{7/3}$ and $4^{6}\cdot2^{40/3}\cdot258<1.1\cdot10^{10}$.
Hence
\begin{equation}\label{eq:PiW}
\Pi\le W\ \text{ in case (1)},\qquad 4^{6}\,\Pi\le W\ \text{ in case (2)}.
\end{equation}

Let $\mathcal C$ be the set in the statement, and let $L$ be a $\mathcal C$-universal integer-valued
lattice, integral in case (1). We have to show that $t\Ok^{+}\subset Q(L)$. By \eqref{eq:PiW}, $L$
represents every element of $t\Ok^{+}$ whose norm is at most $\Pi$. We will first find a sublattice
$L_5$ of $L$ of rank five (Step 1), and we will determine the primes at which $L_5$ may fail to
represent all local integers (Step 2). Then we will fix an element $A\in t\Ok^{+}$ and add to $L_5$
one further vector of $L$, chosen in terms of $A$, such that the resulting lattice $L_A$ represents
$A$ at every prime and satisfies $\Nm(\vv L_A)\le\Pi$ (Steps 3 and 4). Finally,
Definition \ref{def:T}, applied to $L_A$ in case (1) and to $(L_A)^{2}$ in case (2), will show that
$L$ represents $A$ if $A\notin\mathcal C$ (Step 5).

\emph{Step 1: a sublattice of rank five.} The elements of $t\mathcal E_\mu$ have norm at most
$t^{2}\max\{g^{2},3\sqrt{d_F}\}\le V_0\le\Pi$, by \eqref{eq:NTheta}, so $L$ represents all of
them. In case (1), Proposition \ref{prop:five} with $\mu=t$ yields an integer $C_1$ with
$2\le C_1\le g$ and linearly independent $v_1,\dots,v_5\in L$ with
\begin{equation}\label{eq:vvalues}
Q(v_1)=t,\quad Q(v_2)=2t,\quad Q(v_3)=C_1t,\quad Q(v_4)=t\Theta_1,\quad Q(v_5)=t\Theta_2 .
\end{equation}
In case (2), the lattice $L^{2}$ of Convention \ref{conv:scaling} is integral, as
$\ssc(L^{2})=2\,\ssc L\subset\nn L\subset\Ok$ by Lemma \ref{lem:normrange}. It is
$\mu\mathcal E_\mu$-universal, as $2t\mathcal E_\mu=\mu\mathcal E_\mu$. Proposition \ref{prop:five},
applied to $L^{2}$ and $\mu=2t$, yields $C_1$ and linearly independent $v_1,\dots,v_5\in L$ with
\eqref{eq:vvalues}. In both cases, put $L_5:=\Ok v_1+\dots+\Ok v_5$, a free lattice with basis
$v_1,\dots,v_5$. By Lemma \ref{lem:hadamard} and \eqref{eq:NTheta},
\[
\Nm(\vv L_5)\le\prod_{i=1}^{5}\NFQ(Q(v_i))=4C_1^{2}t^{10}\NFQ(\Theta_1)\NFQ(\Theta_2)
<4g^{2}t^{10}\bigl(\sqrt{d_F}+1\bigr)\bigl(2\sqrt{d_F}+1\bigr)\le V_0,
\]
where the last step holds because $(\sqrt{d_F}+1)(2\sqrt{d_F}+1)\le\frac94d_F$ for
$\sqrt{d_F}\ge40.5$.

\emph{Step 2: the odd primes of large volume.} As $2B(L,L)\subset\nn L\subset\Ok$, the entries of
the Gram matrix of $v_1,\dots,v_5$ lie in $\frac12\Ok$, and in $\Ok$ in case (1). Hence
$\ord_\PP\vv L_5\ge0$ at every odd prime $\PP$, and at every dyadic prime $\pp$ we have
$\ord_\pp\vv L_5\ge0$ in case (1) and $\ord_\pp\vv L_5\ge-5e_\pp$ in case (2). Let $\mathcal P$ be
the set of odd primes $\PP$ with $\ord_\PP\vv L_5\ge3$, and put
$\mathfrak f:=\prod_{\PP\in\mathcal P}\PP$. As
$\prod_{\pp\mid2}\Nm\pp^{\,5e_\pp}=2^{5\sum_\pp e_\pp f_\pp}=2^{10}$, we get
\begin{equation}\label{eq:fbound}
(\Nm\mathfrak f)^{3}\le\prod_{\PP\ \text{odd}}\Nm\PP^{\,\ord_\PP\vv L_5}
=\Nm(\vv L_5)\prod_{\pp\mid2}\Nm\pp^{-\ord_\pp\vv L_5}\le V .
\end{equation}
Now let $\PP$ be an odd prime with $\PP\notin\mathcal P$. Then $\ssc(L_5)_\PP\subset\Ok_\PP$, as $2$
is a unit at $\PP$, and $\ord_\PP\vv L_5\le2$. So Lemma \ref{lem:goodprimes} gives
\begin{equation}\label{eq:L5good}
Q\bigl((L_5)_\PP\bigr)=\Ok_\PP\qquad\text{for every odd prime }\PP\notin\mathcal P .
\end{equation}

\emph{Step 3: the auxiliary vector.} From now on we fix an element $A\in t\Ok^{+}$. Let
$\mathcal S$ be the finite set that consists of the primes in $\mathcal P$ and of the dyadic primes.
For each $\PP\in\mathcal S$ we choose $\alpha_\PP\in\Ok_\PP$ with
$j_\PP:=\ord_\PP\alpha_\PP\in\{0,1\}$ that lies in the class of $A/t$ in
$F_\PP^{\times}/(F_\PP^{\times})^{2}$. We will find a totally positive $\beta\in\Ok$ of small norm
with
\begin{equation}\label{eq:betaclass}
\beta\in\alpha_\PP(\Ok_\PP^{\times})^{2}\qquad\text{for every }\PP\in\mathcal S .
\end{equation}
Put
\[
\mathfrak d:=\prod_{\pp\mid2}\pp^{\,j_\pp+2e_\pp+1},\qquad
\mathfrak f_A:=\prod_{\PP\in\mathcal P}\PP^{\,j_\PP+1},\qquad\mmm:=\mathfrak d\,\mathfrak f_A .
\]
As $e_\PP=0$ for $\PP\in\mathcal P$, we have
$\mmm=\prod_{\PP\in\mathcal S}\PP^{\,j_\PP+2e_\PP+1}$. So by Lemma \ref{lem:pin} there is $a\in\Ok$
such that every $\beta\in a+\mmm$ satisfies \eqref{eq:betaclass}.

By Lemma \ref{lem:dyadic}, $\Nm(\mathfrak d)\le256$ and $16\in\mathfrak d$. The ideal
$\mathfrak f_A$ divides $\mathfrak f^{2}$, and so $\Nm(\mathfrak f_A)\le(\Nm\mathfrak f)^{2}\le V^{2/3}$
by \eqref{eq:fbound}. Moreover, the integer $\Nm(\mathfrak f_A)$ lies in $\mathfrak f_A$. Hence
$16\,\Nm(\mathfrak f_A)\in\mmm$ and $\Nm(\mmm)\le256\,\Nm(\mathfrak f_A)$. By Lemma \ref{lem:GoN}, there
is a totally positive $\beta\in a+\mmm$ with
\[
\NFQ(\beta)\le256\,\Nm(\mathfrak f_A)^{2}+128\,\Nm(\mathfrak f_A)\sqrt{d_F}
\le256\,V^{4/3}+128\,V^{2/3}\sqrt{d_F}\le258\,V^{4/3}.
\]
For the last step, note that $V\ge225\,d_F$ and $d_F\ge1641$ give
$V^{2/3}\ge36\,d_F^{2/3}\ge64\sqrt{d_F}$.

The element $t\beta$ lies in $t\Ok^{+}$, and
$\NFQ(t\beta)=t^{2}\NFQ(\beta)\le258\,t^{2}V^{4/3}\le\Pi$. So $L$ represents $t\beta$. Choose
$x\in L$ with $Q(x)=t\beta$, and put
\[
L_A:=L_5+\Ok x\subset L .
\]

\emph{Step 4: the lattice $L_A$.} The lattice $L_A$ contains $L_5$ and is generated by the six
vectors $v_1,\dots,v_5,x$, and so $5\le\rank L_A\le6$. By Lemma \ref{lem:hadamard}, applied to these
generators, and by Steps 1 and 3,
\[
\Nm(\vv L_A)\le\prod_{i=1}^{5}\NFQ(Q(v_i))\cdot\NFQ(t\beta)\le V_0\cdot258\,t^{2}V^{4/3}=\Pi .
\]
We claim that $A\in Q\bigl((L_A)_\PP\bigr)$ for every prime $\PP$. If $\PP\notin\mathcal S$, then
$\PP$ is odd and $\PP\notin\mathcal P$, and so $A\in\Ok_\PP=Q\bigl((L_5)_\PP\bigr)$ by
\eqref{eq:L5good}. Let $\PP\in\mathcal S$. By \eqref{eq:betaclass}, the lattice $(L_A)_\PP$
represents the element $Q(x)=t\beta$ of $t\alpha_\PP(\Ok_\PP^{\times})^{2}$. As $A\in t\Ok_\PP$ and
$A/t$ lies in the class of $\alpha_\PP$, Lemma \ref{lem:squareclasses} shows that
$A\in Q\bigl((L_A)_\PP\bigr)$.

\emph{Step 5: conclusion.} If $A\in\mathcal C$, then $A\in Q(L)$ by assumption. Assume that
$A\notin\mathcal C$.

In case (1) we then have $\NFQ(A)>\mathcal T_{F,6}(W)$. The lattice $L_A$ is integral, as
$\ssc L_A\subset\ssc L\subset\Ok$. Moreover $5\le\rank L_A\le6$ and $\Nm(\vv L_A)\le\Pi\le W$, and
$A\in Q\bigl((L_A)_\PP\bigr)$ for every prime $\PP$ by Step 4. By Definition \ref{def:T},
$A\in Q(L_A)\subset Q(L)$.

In case (2) we have $\NFQ(2A)=4\NFQ(A)>\mathcal T_{F,6}(W)$. The lattice $(L_A)^{2}$ is integral, as
in Step 1. It has the same rank as $L_A$, and $\vv\bigl((L_A)^{2}\bigr)=2^{\rank L_A}\vv L_A$, so
$\Nm\bigl(\vv((L_A)^{2})\bigr)=4^{\rank L_A}\Nm(\vv L_A)\le4^{6}\,\Pi\le W$. Moreover
$2A\in2Q\bigl((L_A)_\PP\bigr)=Q\bigl(((L_A)^{2})_\PP\bigr)$ for every prime $\PP$. By
Definition \ref{def:T}, $2A\in Q\bigl((L_A)^{2}\bigr)=2Q(L_A)$, and so $A\in Q(L_A)\subset Q(L)$. This
completes the proof of Theorem \ref{thm:tuniv}.

\subsection{Proof of Theorem \ref{cor:critexplicit}}

A set that contains a universality criterion set of one of the two kinds is again one of this
kind. By
Corollary \ref{cor:tone}, it therefore suffices to show that
\begin{equation}\label{eq:TW}
\max\{\mathcal T_{F,6}(W),W\}\le\exp\bigl((24\,d_F)^{E}\bigr)
\end{equation}
with $E=476$ for $W=8\cdot10^{7}\,d_F^{7/3}$ and $d_F\ge1641$, and with $E=486$ for
$W=3.4\cdot10^{15}\,d_F^{7/3}$ and $d_F\ge25\,761$.

Let $W\ge1$ and $m\in\{5,6\}$. For $l=2$, the bound of Theorem \ref{cor:crude} for rank $m$ is
\[
\exp\Bigl((4m\,d_F)^{12m^{2}}D^{\rho_m}\bigl(\log(4D)\bigr)^{5m-1}\Bigr).
\]
It is non-decreasing in $D$. For $1\le D\le W$ it is at most
$\exp\bigl((24\,d_F)^{432}\,W^{10}(\log(4W))^{29}\bigr)$, because $4m\le24$, $12m^{2}\le432$,
$\rho_m\le\rho_5=10$, $5m-1\le29$, and $\log(4W)>1$. Hence, by Definition \ref{def:T},
\[
\mathcal T_{F,6}(W)\le\exp\Bigl((24\,d_F)^{432}\,W^{10}\bigl(\log(4W)\bigr)^{29}\Bigr).
\]

Now let $W=K\,d_F^{7/3}$ and $d_F\ge d_0$, where $(K,d_0)=(8\cdot10^{7},1641)$ or
$(K,d_0)=(3.4\cdot10^{15},25\,761)$. Put $\lambda:=\log(24\,d_F)$ and $\lambda_0:=\log(24\,d_0)$. Then
$\lambda\ge\lambda_0$, and
\[
\log W=\log K-\tfrac73\log24+\tfrac73\lambda\le a_1\lambda,\qquad
\log(4W)=\log(4K)-\tfrac73\log24+\tfrac73\lambda\le a_2\lambda,
\]
with $a_1:=(\log K-\frac73\log24)/\lambda_0+\frac73$ and
$a_2:=(\log(4K)-\frac73\log24)/\lambda_0+\frac73$, as both numerators are positive. Moreover
$\log(a_2\lambda)\le a_3\lambda$ with $a_3:=\log(a_2\lambda_0)/\lambda_0$, because $a_2\lambda_0>e$
and $x\mapsto\log(a_2x)/x$ is decreasing for $a_2x\ge e$. Therefore
\[
(24\,d_F)^{432}\,W^{10}\bigl(\log(4W)\bigr)^{29}\le(24\,d_F)^{E_0}\qquad\text{with}\quad
E_0:=432+10a_1+29a_3 .
\]
For $(K,d_0)=(8\cdot10^{7},1641)$ we find $a_1<3.353$ and $a_3<0.341$, so $E_0<476$. For
$(K,d_0)=(3.4\cdot10^{15},25\,761)$ we find $a_1<4.460$ and $a_3<0.3082$, so $E_0<486$. In both cases
also $\log W\le a_1\lambda\le e^{\lambda}=24\,d_F$, as $a_1<5$ and $\lambda>10$. This proves
\eqref{eq:TW}. \qed

\medskip
Inserting Theorem \ref{cor:crude} into Theorem \ref{thm:tuniv} in the same way gives explicit
criterion sets for $t\Ok^{+}$-universality.

\section{The fields \texorpdfstring{$\Q(\sqrt{n^{2}-1})$}{Q(sqrt(n\string^2-1))}}\label{sec:pell}

Throughout this section $n\ge2$ is an integer such that $d:=n^{2}-1$ is squarefree, and
$F=\Q(\sqrt d)$. Then $n$ is even, as $8\mid n^{2}-1$ for odd $n$. Hence $d\equiv3\pmod4$,
$\Ok=\Z[\sqrt d]$, $\omega=\sqrt d$, and $d_F=4d$. Put $\eps:=n+\sqrt d$.

In this section we will prove Theorem \ref{thm:pell}. The proof does not use
the results of \S\ref{sec:prelim} and of Sections \ref{sec:lower}--\ref{sec:crit}. It will use two facts about these
fields. First, every totally positive integer is, up to the square of a unit,
of the form $k+k'\eps^{\pm1}$ with $k,k'\in\Z_{\ge0}$ (Lemma \pref{lem:pell}{3}). Second, an element
of $\Ok$ whose conjugates are small is rational (Lemma \pref{lem:pell}{4}), and so vectors that
represent small rational integers span lattices over $\Z$ (Proposition \ref{prop:pellsuff}).

A \emph{$\Z$-lattice} is a free $\Z$-module $T$ of finite rank together with a symmetric bilinear
form $B_T\colon T\times T\to\Q$ such that $Q_T(x):=B_T(x,x)$ lies in $\Z_{>0}$ for every nonzero
$x\in T$, i.e., a positive definite integer-valued lattice over $\Z$. It is \emph{integral} if $B_T(T,T)\subset\Z$, and \emph{universal} if
$\Z_{>0}\subset Q_T(T)$. For $\mathcal C\subset\Z_{>0}$, it is \emph{$\mathcal C$-universal} if
$\mathcal C\subset Q_T(T)$. If $T$ is not universal, its \emph{truant} is the least positive integer
that it does not represent. For $a_1,\dots,a_s\in\Z_{>0}$ we write $\langle a_1,\dots,a_s\rangle$
for $\Z^{s}$ with $Q_T(x)=\sum_ia_ix_i^{2}$.

For a $\Z$-lattice $T$ we write $T\otimes\Ok$ for the
$\Ok$-module $T\otimes_\Z\Ok\subset T\otimes_\Z F$, where $T\otimes_\Z F$ carries the
$F$-bilinear extension of $B_T$. Diagonalizing $B_T$ over $\Q$ shows that $T\otimes_\Z\R$ is
positive definite, and so $T\otimes\Ok$ is positive definite. Moreover
$2B_T(T,T)\subset\Z$, which gives $\nn(T\otimes\Ok)\subset\Ok$. Hence $T\otimes\Ok$ is an
integer-valued lattice, and it is integral if $T$ is.

\begin{lemma}\label{lem:pell}
\begin{enumerate}
\item\label{lem:pell:1} $\eps$ is a totally positive unit, $\mathcal U_F=\{\pm\eps^{s}\mid s\in\Z\}$,
$\mathcal U_F^{+}=\{\eps^{s}\mid s\in\Z\}$, $\mathcal U_F^{2}=\{\eps^{2s}\mid s\in\Z\}$, and
$\eps+\eps^{-1}=2n$.
\item\label{lem:pell:2} For every $j\in\Z$, the elements $\eps^{j},\eps^{j+1}$ form a $\Z$-basis of $\Ok$.
\item\label{lem:pell:3} Every $\alpha\in\Ok^{+}$ satisfies $\alpha=k\eps^{j}+k'\eps^{j+1}$ for some $j\in\Z$ and
$k,k'\in\Z_{\ge0}$. Consequently there are $\eta\in\mathcal U_F$, $\sigma\in\{1,-1\}$, and
$k,k'\in\Z_{\ge0}$ with $\alpha\eta^{2}=k+k'\eps^{\sigma}$.
\item\label{lem:pell:4} If $\beta\in\Ok$ satisfies $|\beta|+|\beta'|<2\sqrt d$, then $\beta\in\Z$.
\item\label{lem:pell:5} If $c,i\in\Z_{>0}$ satisfy $cn<d$, then $c\eps\ne i\gamma^{2}$ and $c\eps^{-1}\ne i\gamma^{2}$
for every $\gamma\in\Ok$.
\end{enumerate}
\end{lemma}

\begin{proof}
(1) We have $\NFQ(\eps)=n^{2}-d=1$, and $\eps$ and $\eps'=n-\sqrt d$ are positive. Let
$u=x+y\sqrt d>1$ be a unit, with $x,y\in\Z$. Then $|u'|=u^{-1}<1$, so $x=\frac12(u+u')>0$ and
$y=(u-u')/(2\sqrt d)>0$. If $y=1$, then $x^{2}=d\pm1\in\{n^{2},n^{2}-2\}$. Here $n^{2}-2$ is not a
square. Indeed, $n^{2}-z^{2}=2$ with $z\in\Z$ is impossible, because $n-z$ and $n+z$ have the same
parity. So $x=n$ and
$u=\eps$. If $y\ge2$, then $x^{2}\ge4d-1$, and $u\ge2\sqrt d+\sqrt{4d-1}>2\sqrt d+1>n+\sqrt d=\eps$,
using $\sqrt d>n-1$. So $\eps$ is the least unit exceeding $1$, and $\mathcal U_F=\{\pm\eps^{s}\}$.
As $(\eps^{s})'=\eps^{-s}$, we get $\mathcal U_F^{+}=\{\eps^{s}\}$ and
$\mathcal U_F^{2}=\{\eps^{2s}\}$. Finally $\eps+\eps^{-1}=\eps+\eps'=2n$.

(2) We have $\Ok=\Z+\Z\sqrt d=\Z+\Z\eps$, and multiplication by $\eps^{j}$ is an automorphism of
the $\Z$-module $\Ok$.

(3) The first assertion is \cite[Lemma 4]{Ki}. If $j=2i$ is even, then $\alpha\eps^{-2i}=k+k'\eps$.
If $j=2i-1$ is odd, then $\alpha\eps^{-2i}=k'+k\eps^{-1}$.

(4) Write $\beta=x+y\sqrt d$ with $x,y\in\Z$. Then $2|y|\sqrt d=|\beta-\beta'|\le|\beta|+|\beta'|
<2\sqrt d$, so $y=0$.

(5) If $\gamma=u+v\sqrt d\in\Ok$ satisfies $i\gamma^{2}=c\eps^{\pm1}=cn\pm c\sqrt d$, then
$i(u^{2}+dv^{2})=cn$ and $2iuv=\pm c\ne0$. So $v\ne0$ and $d\le idv^{2}\le cn$, contrary to the
assumption.
\end{proof}

Kim \cite{Ki} used part (3) to show that the octonary form
$\langle1,1,1,1,\eps,\eps,\eps,\eps\rangle$ is universal over $F$, and Kominers \cite{Kom3}
determined all real quadratic fields over which this form is universal. Part (3) implies that in the
fields $F$ of this section the indecomposable elements of $\Ok^{+}$, i.e., those that are not sums
of two elements of $\Ok^{+}$, are exactly the powers of $\eps$ (compare \cite[\S6.1]{KKR}). More
generally, Hejda--Kala \cite[Theorem 2]{HK} showed that in every real quadratic field each totally
positive integer is a combination, with nonnegative integer coefficients, of two consecutive
indecomposable elements.

\begin{proposition}\label{prop:pellsuff}
Let $\mathcal C\subset\Z_{>0}$ be finite, put $c:=(\max\mathcal C)^{2}$, let $\delta\in\{1,4\}$,
and assume $n>1+\delta c/2$. Assume that every $\mathcal C$-universal $\Z$-lattice, integral if
$\delta=1$, is universal. Let $L$ be a $(\mathcal C\cup\eps\,\mathcal C)$-universal integer-valued
lattice,
integral if $\delta=1$. Then $L$ is universal.
\end{proposition}

\begin{proof}
Put $\theta:=\sqrt\delta\in\{1,2\}$. Then $\theta B(x,y)\in\Ok$ for $x,y\in L$. For $\delta=1$
this holds because $\ssc L\subset\Ok$, and for $\delta=4$ because $2B(x,y)=Q(x+y)-Q(x)-Q(y)$. For
$a\in\mathcal C$ choose $w_a,w'_a\in L$ with $Q(w_a)=a$ and $Q(w'_a)=a\eps$. We will use the
Cauchy--Schwarz inequality $B(x,y)^{2}\preceq Q(x)Q(y)$, valid in both embeddings as $L$ is
positive definite, and Lemma \pref{lem:pell}{4}. We have $ab\le c$ for $a,b\in\mathcal C$,
and $n>1+\delta c/2$ gives $\delta c<2(n-1)$. Hence $\delta c<(n-1)(n+1)=d$ and
$\delta c<2n-1<\eps$.

(a) Let $a,b\in\mathcal C$ and $\beta:=\theta B(w_a,w_b)\in\Ok$. Then $\beta^{2}\preceq\delta ab
\le\delta c$, so $|\beta|+|\beta'|\le2\sqrt{\delta c}<2\sqrt d$, and $\beta\in\Z$.

(b) Let $a,b\in\mathcal C$ and $\beta:=\theta\eps^{-1}B(w'_a,w'_b)\in\Ok$. From
$B(w'_a,w'_b)^{2}\preceq ab\,\eps^{2}$ we get $\beta^{2}\preceq\delta ab\le\delta c$, and again
$\beta\in\Z$.

(c) Let $a,b\in\mathcal C$ and $\beta:=\theta B(w_a,w'_b)\in\Ok$. Then
$\beta^{2}\preceq\delta ab\,\eps$, i.e., $\beta^{2}\le\delta ab\,\eps$ and
$\beta'^{2}\le\delta ab\,\eps^{-1}$. Using $\eps+\eps^{-1}=2n$ and $\delta c<2(n-1)$, we get
\[
|\beta|+|\beta'|\le\sqrt{\delta c}\,\bigl(\eps^{1/2}+\eps^{-1/2}\bigr)
=\sqrt{\delta c\,(\eps+2+\eps^{-1})}=\sqrt{2\delta c\,(n+1)}<2\sqrt{(n-1)(n+1)}=2\sqrt d .
\]
So $\beta\in\Z$, and then $\beta^{2}=\beta'^{2}\le\delta c\,\eps^{-1}<1$ gives $\beta=0$.

Put $Z_1:=\sum_{a\in\mathcal C}\Z w_a$ and $Z_2:=\sum_{a\in\mathcal C}\Z w'_a$, which are free
$\Z$-modules of finite rank. By (a), we have $Q(Z_1)\subset\Z$ and $\theta B(Z_1,Z_1)\subset\Z$.
As $L$ is positive definite, $Z_1$ with the restriction of $B$ is therefore a $\Z$-lattice,
integral if $\delta=1$, which represents every $a\in\mathcal C$. By (b), the same holds for $Z_2$
with the form $\eps^{-1}B$. By hypothesis, both are universal. Thus, for every $h\in\Z_{\ge0}$ there
are $u\in Z_1$ and $u'\in Z_2$ with $Q(u)=h$ and $Q(u')=h\eps$. Moreover $B(Z_1,Z_2)=0$ by (c).

Now let $\alpha\in\Ok^{+}$. By Lemma \pref{lem:pell}{3}, $\alpha=k\eps^{j}+k'\eps^{j+1}$ with
$j\in\Z$ and $k,k'\in\Z_{\ge0}$. One of $j,j+1$ is even, say $2i$, and the other is odd, say
$2i'+1$, so that $\alpha=k_1\eps^{2i}+k_2\eps^{2i'+1}$ with $\{k_1,k_2\}=\{k,k'\}$. Choose
$u\in Z_1$ with $Q(u)=k_1$ and $u'\in Z_2$ with $Q(u')=k_2\eps$, and put
$x:=\eps^{i}u+\eps^{i'}u'\in L$. Then
\[
Q(x)=\eps^{2i}Q(u)+\eps^{2i'}Q(u')+2\eps^{i+i'}B(u,u')=k_1\eps^{2i}+k_2\eps^{2i'+1}=\alpha .
\qedhere
\]
\end{proof}

\begin{proposition}\label{prop:pellnec}
Let $a\in\Z_{>0}$ with $an<d$, let $T$ be a $\Z$-lattice with truant $a$, and put
$L_a:=(T\otimes\Ok)\perp\langle\eps,2\eps,\dots,(a-1)\eps\rangle$. Then $L_a$ is an integer-valued
lattice,
integral if $T$ is, which represents every $\alpha\in\Ok^{+}$ with $\NFQ(\alpha)<a^{2}$ but no
element of $a\,\mathcal U_F^{2}$. The lattice $L_a^{\eps}$ represents every $\alpha\in\Ok^{+}$ with
$\NFQ(\alpha)<a^{2}$ but no element of $a\eps\,\mathcal U_F^{2}$.
\end{proposition}

\begin{proof}
The assertion on $L_a^{\eps}$ follows from that on $L_a$. Indeed, $L_a^{\eps}$ represents $\alpha$
if and only if $L_a$ represents $\eps^{-1}\alpha$, and $\NFQ(\eps^{-1}\alpha)=\NFQ(\alpha)$ and
$\eps^{-1}\Ok^{+}=\Ok^{+}$.

Let $\alpha\in\Ok^{+}$ with $\NFQ(\alpha)<a^{2}$. By Lemma \pref{lem:pell}{3},
$\alpha\eta^{2}=k+k'\eps^{\sigma}$ with $\eta\in\mathcal U_F$, $\sigma\in\{1,-1\}$, and
$k,k'\in\Z_{\ge0}$. Then $\NFQ(\alpha)=\NFQ(k+k'\eps^{\sigma})=k^{2}+2nkk'+k'^{2}$, so
$k,k'\le a-1$. Now $T$ represents $k$, and if $k'\ge1$, the component $\langle k'\eps\rangle$
represents $k'\eps^{\sigma}$, namely $k'\eps\cdot1^{2}$ and $k'\eps\cdot(\eps^{-1})^{2}$. Hence
$L_a$ represents $\alpha\eta^{2}$, and therefore $\alpha$.

Suppose that $a=Q(x)+\sum_{i=1}^{a-1}i\eps y_i^{2}$ with $x\in T\otimes\Ok$ and $y_i\in\Ok$. Each
summand $\gamma$ satisfies $0\preceq\gamma\preceq a$. Since $an<d$ gives $a\le n-1<\sqrt d$, we get
$|\gamma|+|\gamma'|\le2a<2\sqrt d$, so $\gamma\in\{0,1,\dots,a\}$ by Lemma \pref{lem:pell}{4}. A
summand $i\eps y_i^{2}=c'$ with $c'\ge1$ would give $c'\eps^{-1}=iy_i^{2}$, which
Lemma \pref{lem:pell}{5} excludes, as $c'n\le an<d$. So every $y_i$ vanishes, and $Q(x)=a$. Write
$x=u+v\sqrt d$ with $u,v\in T$. Then $Q(x)=Q_T(u)+dQ_T(v)+2B_T(u,v)\sqrt d$, so
$Q_T(u)+dQ_T(v)=a$. As $d>a$, this forces $v=0$ and $Q_T(u)=a$, contrary to the choice of $T$.
Finally, $Q(L_a)\eta^{2}=Q(L_a)$ for $\eta\in\mathcal U_F$, so $L_a$ represents no element of
$a\,\mathcal U_F^{2}$.
\end{proof}

\begin{proof}[Proof of Theorem \textup{\ref{thm:pell}}]
For Theorem \pref{thm:pell}{1} let $\mathcal C:=\mathcal C_{15}$, $\delta:=1$, and
$n_0:=114$. For Theorem \pref{thm:pell}{2} let $\mathcal C:=\mathcal C_{290}$, $\delta:=4$, and
$n_0:=168\,202$. Put $c:=(\max\mathcal C)^{2}$, so that $c=225$, respectively $c=84\,100$. Then
$n\ge n_0$ gives $n>1+\delta c/2$, which is $n>113.5$, respectively $n>168\,201$. By the $15$-Theorem
\cite{Bh}, respectively the $290$-Theorem \cite{BH}, every $\mathcal C$-universal
$\Z$-lattice, integral in the case of Theorem \pref{thm:pell}{1}, is universal. Hence, by
Proposition \ref{prop:pellsuff}, $S_n$ is a universality criterion set for integral lattices, and
$S'_n$ is a universality criterion set for integer-valued lattices.

The classes $s\,\mathcal U_F^{2}$ with $s\in\mathcal C\cup\eps\,\mathcal C$ are pairwise distinct.
Indeed, if $a,b\in\mathcal C$, $i,j\in\{0,1\}$, and $a\eps^{i}\mathcal U_F^{2}=b\eps^{j}\mathcal U_F^{2}$,
then $a\eps^{\,i-j+2s}=b$ for some $s\in\Z$ by Lemma \pref{lem:pell}{1}. As $\eps^{k}\notin\Q$ for
$k\ne0$, this gives $i=j$ and $a=b$.

Now let $a\in\mathcal C$. There is a $\Z$-lattice $T_a$ with truant $a$, integral in the case of
Theorem \pref{thm:pell}{1}. For $a=1$ take $T_1:=\langle2\rangle$. In the case of
Theorem \pref{thm:pell}{1}, the integers $2,3,5,6,7,10,14,15$ are the truants of $\langle1\rangle$,
$\langle1,1\rangle$, $\langle1,2\rangle$, $\langle1,1,3\rangle$, $\langle1,1,1\rangle$,
$\langle1,2,3\rangle$, $\langle1,1,2\rangle$, $\langle1,2,5,5\rangle$, respectively \cite{Bh}. In the
case of Theorem \pref{thm:pell}{2}, every $a\in\mathcal C_{290}$ is the truant of a $\Z$-lattice
\cite{BH}. Since $an\le15n<n^{2}-1=d$ in the first case and $an\le290n<d$ in the second,
Proposition \ref{prop:pellnec} shows that $a\,\mathcal U_F^{2}$ and $a\eps\,\mathcal U_F^{2}$ are
truants, in the sense of \cite[Definition 2.2]{KKR}, of the lattices $L_a$ and $L_a^{\eps}$ built
from $T_a$. These lattices are integral in the case of Theorem \pref{thm:pell}{1}. By
\cite[Proposition 3.1 and Theorem 3.3]{KKR}, every universality criterion set for integer-valued lattices (for integral
lattices, in the case of Theorem \pref{thm:pell}{1}) therefore contains an element of
$a\,\mathcal U_F^{2}$ and an element of $a\eps\,\mathcal U_F^{2}$.

Finally, let $s\in\mathcal C\cup\eps\,\mathcal C$. As the classes above are pairwise distinct, the
set $(\mathcal C\cup\eps\,\mathcal C)\smallsetminus\{s\}$ contains no element of
$s\,\mathcal U_F^{2}$. So it is not a universality criterion set of the respective kind, and
$s$ cannot be omitted.
\end{proof}

\section*{Acknowledgments}
We thank Jakub Kr\'asensk\'y and the members of the
\href{https://www1.karlin.mff.cuni.cz/~kala/web/ufoclan}{UFOCLAN} group for interesting discussions.

\section*{Statement on the use of AI}
The initial ideas for all main results and for their proofs are due to the authors. The text was
written to a large extent with AI tools. AI tools also contributed mathematically: substantially in
Sections \ref{sec:global}, \ref{sec:r}, \ref{sec:rq}, and \ref{sec:crit}, including the explicit constants of
Theorems \ref{cor:crude} and \ref{cor:critexplicit}, and by streamlining arguments elsewhere. We
primarily used Claude Fable 5.2 and Claude Opus 5.5
(Anthropic), with contributions from GPT-6 (OpenAI).

\end{document}